\documentclass[a4paper,10pt]{article}
\usepackage{amssymb}
\usepackage{fullpage}
\usepackage{graphicx}        % standard LaTeX graphics tool
\usepackage{amscd}
\usepackage{amsmath}
\usepackage{mathrsfs}
\usepackage{verbatim}

\newtheorem{theorem}{Theorem}
\newtheorem{lemma}[theorem]{Lemma}
\newtheorem{definition}[theorem]{Definition}
\newtheorem{corollary}[theorem]{Corollary}

\newtheorem{remark}[theorem]{Remark}
\newcommand{\R}{\mathbb{R}}

\newcommand{\ul}{\textbf}

\newcommand{\PPP}{\mbox{\boldmath$\psi$}}

\newcommand{\desda}{\Leftrightarrow}
\newcommand{\ve}{\varepsilon}
\DeclareMathAlphabet\gothic{U}{euf}{m}{n}

\begin{document}

%\begin{frontmatter}

%% Title, authors and addresses

%% use the tnoteref command within \title for footnotes;
%% use the tnotetext command for the associated footnote;
%% use the fnref command within \author or \address for footnotes;
%% use the fntext command for the associated footnote;
%% use the corref command within \author for corresponding author footnotes;
%% use the cortext command for the associated footnote;
%% use the ead command for the email address,
%% and the form \ead[url] for the home page:
%%
%% \title{Title\tnoteref{label1}}
%% \tnotetext[label1]{}
%% \author{Name\corref{cor1}\fnref{label2}}
%% \ead{email address}
%% \ead[url]{home page}
%% \fntext[label2]{}
%% \cortext[cor1]{}
%% \address{Address\fnref{label3}}
%% \fntext[label3]{}

\title{Evolution Equations on Gabor Transforms and their Applications}

%% use optional labels to link authors explicitly to addresses:
%% \author[label1,label2]{<author name>}
%% \address[label1]{<address>}
%% \address[label2]{<address>}

\author{Remco Duits and Hartmut F\"{u}hr and Bart Janssen and \\
Mark Bruurmijn and Luc Florack and Hans van Assen
}

\maketitle
%\address{Remco Duits, Department of Mathematics and Computer Science \& Department of Biomedical Engineering, Eindhoven University of Technology, 5600 MB Eindhoven, The Netherlands, \\
%\emph{R.Duits@tue.nl} \\
%Hartmut F\"{u}hr, Lehrstuhl A f\"ur Mathematik, RWTH Aachen University, 52056 Aachen, Germany, \emph{fuehr@MathA.rwth-aachen.de} \\
%Bart Janssen, Department of Mathematics and Computer Science, Eindhoven University of Technology, 5600 MB Eindhoven, The Netherlands, \\
%\emph{B.J.Janssen@tue.nl} \\
%Mark Bruurmijn, Department of Biomedical Engineering, Eindhoven University of Technology, 5600 MB Eindhoven, The Netherlands, \emph{L.C.M.Bruurmijn@student.tue.nl} \\
%Luc Florack, Department of Mathematics and Computer Science \& Department of Biomedical Engineering, Eindhoven University of Technology, 5600 MB Eindhoven, The Netherlands, \\
%\emph{L.M.J.Florack@tue.nl} \\
%Hans van Assen, Department of Biomedical Engineering, Eindhoven University of Technology, 5600 MB Eindhoven, The Netherlands, \emph{H.C.v.Assen@tue.nl} \\
%}
\begin{abstract}
%% Text of abstract
We introduce a systematic approach to the design, implementation
and analysis of left-invariant evolution schemes acting on Gabor transform,
primarily for applications in signal and image analysis. Within this approach we relate operators on signals to
operators on Gabor transforms. In order to obtain a translation and modulation invariant operator on the space of signals, the corresponding operator on the reproducing kernel space of Gabor transforms must be left invariant, i.e. it should commute with the left regular action of the reduced Heisenberg group $H_r$. By using the left-invariant vector fields on $H_r$ in the generators of our evolution equations on Gabor transforms, we naturally employ the essential group structure on the domain of a Gabor transform. Here we distinguish between two tasks.
Firstly, we consider non-linear adaptive left-invariant convection (reassignment) to sharpen Gabor transforms, while maintaining the original signal. Secondly, we consider signal enhancement via left-invariant diffusion on the corresponding Gabor transform.
We provide numerical experiments and analytical evidence for our methods and we consider an explicit medical imaging application.
\end{abstract}

%\begin{keyword}
%% keywords here, in the form: keyword \sep keyword
\textbf{keywords: }
Evolution equations, Heisenberg group , Differential reassignment , Left-invariant vector fields , Diffusion on Lie groups , Gabor transforms , Medical imaging.
%% MSC codes here, in the form: \MSC code , code
%35R03 , 22E25 ,  35Q68
%% or \MSC[2008] code \sep code (2000 is the default)
%\end{keyword}

%\end{frontmatter}

%%
%% Start line numbering here if you want
%%
% \linenumbers

%% main text
%\section{}
%\label{}

%% The Appendices part is started with the command \appendix;
%% appendix sections are then done as normal sections
%% \appendix

%% \section{}
%% \label{}

%% References
%%
%% Following citation commands can be used in the body text:
%% Usage of \cite is as follows:
%%   \cite{key}         ==>>  [#]
%%   \cite[chap. 2]{key} ==>> [#, chap. 2]
%%

\section{Introduction \label{ch:1}}

The Gabor transform of a signal $f \in \mathbb{L}_2(\mathbb{R}^d)$ is a
function $\mathcal{G}_{\psi}[f]: \R^{d} \times \R^d \to \mathbb{C}$ that
can be roughly understood as a musical score of $f$, with
$\mathcal{G}_{\psi}[f] (p,q)$ describing
the contribution of frequency $q$ to the behaviour of $f$ near $p$
\cite{Gabor1946,Helstrom1966}. This interpretation is
necessarily of limited precision, due to the various uncertainty
principles, but it has nonetheless turned out to be a very rich
source of mathematical theory as well as practical signal processing
algorithms.

\begin{figure}
\centerline{\includegraphics[width=0.7\hsize]{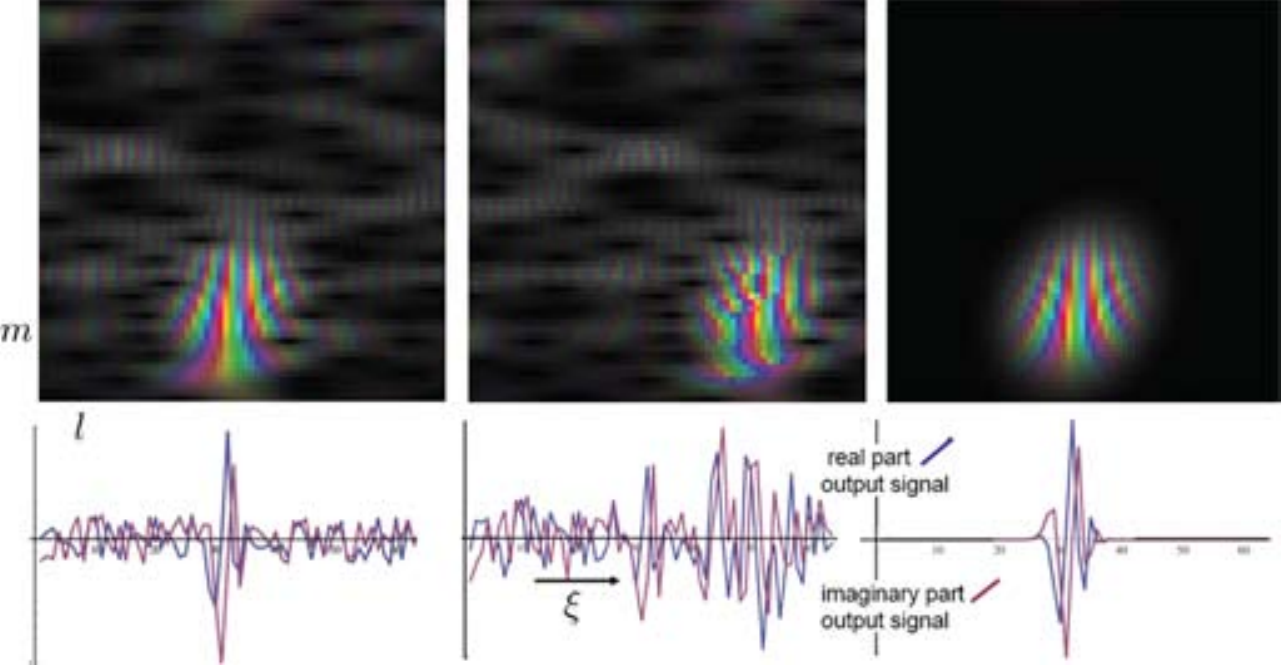}}
\caption{
Top row from left to right, (1) the Gabor transform of original signal $f$, (2) processed Gabor transform $\Phi_{t}(\mathcal{W}_{\psi}f)$ where $\Phi_{t}$ denotes a phase invariant shift (for more elaborate adaptive convection/reassignment operators see Section \ref{ch:reass} where we operationalize the theory in \cite{Daudet}) using a discrete Heisenberg group, where $l$ represents discrete spatial shift and $m$ denotes discrete local frequency, (3) processed Gabor transform $\Phi_{t}(\mathcal{W}_{\psi}f)$ where $\Phi_{t}$ denotes a \emph{phase covariant } diffusion operator on Gabor transforms with stopping time $t>0$. For details on phase covariant diffusions on Gabor transforms, see \cite[ch:7]{DuitsFuehrJanssen} and \cite[ch:6]{Janssen2009b}. Note that phase-covariance is preferable over phase invariance. For example, restoration of the old phase in the phase invariant shift creates noisy artificial patterns (middle image) in the phase of the transported strong responses in the Gabor domain. Bottom row, from left to right: (1) Original complex-valued signal $f$, (2) output signal $\Upsilon_{\psi}f=\mathcal{W}_{\psi}^{*}\Phi_{t}\mathcal{W}_{\psi}f$ where $\Phi_t$ denotes a \emph{phase-invariant} spatial shift (due to phase invariance the output signal looks bad and clearly phase invariant spatial shifts in the Gabor domain do not correspond to spatial shifts in the signal domain), (3) Output signal $\Upsilon_{\psi}f=\mathcal{W}_{\psi}^{*}\Phi_{t}\mathcal{W}_{\psi}f$ where $\Phi_t$ denotes \emph{phase-covariant} adaptive diffusion in the Gabor domain with stopping time $t>0$.
% Bottom row: explanation color-coding of the phase in the corresponding Gabor transforms in the top row.
}
\label{fig:translation}
\end{figure}

The use of a window function for the Gabor transform results in a
smooth, and to some extent blurred, time-frequency
representation. %;
%though keep in mind that by the uncertainty principle,
%there is no such thing as a ``true time-frequency
%representation''.
For purposes of signal analysis, say for the extraction of instantaneous
frequencies, various authors tried to improve the
resolution of the Gabor transform, literally in order to sharpen the
time-frequency picture of the signal; this type
of procedure is often called ``reassignment'' in the literature.
For instance, Kodera et al.~\cite{Kodera1976} studied techniques for the
enhancement of the spectrogram, i.e. the
squared modulus of the short-time Fourier transform. Since the phase of
the Gabor transform is neglected,
the original signal is not easily recovered from the reassigned
spectrogram. Since then, various authors developed
reassignment methods that were intended to allow (approximate) signal
recovery \cite{Auger1995,Chassande,Daudet}. We claim that a proper treatment
of phase may be understood as {\em phase-covariance},
rather than {\em phase-invariance}, as advocated previously. An illustration of this
claim is contained in Figure \ref{fig:translation}.

We adapt the group theoretical approach developed for the Euclidean
motion groups in the recent works
\cite{DuitsRTHESIS,FrankenPhDThesis,DuitsAMS1,DuitsAMS2,DuitsR2006AMS,Duits2005IJCV,DuitsIJCV2010,Prckovska,Rodrigues},
thus illustrating the scope of the methods devised for general Lie groups in
\cite{DuitsR2006SS2} in signal and image processing.
Reassignment will be
seen to be a special case of left-invariant convection. A useful source of ideas specific to Gabor analysis and reassignment was the paper~\cite{Daudet}.

\subsection{Structure of the article}

This article provides a systematic approach to the design, implementation and analysis of left-invariant evolution
schemes acting on Gabor transform, primarily for applications in signal and image analysis.
The article is structured as follows:
\begin{itemize}
 \item {\bf Introduction and preliminaries:} Sections \ref{ch:1}, \ref{ch:Hr}. In Section \ref{ch:Hr} we consider time-frequency analysis and its inherit connection to the Heisenberg group and subsequently we introduce relevant structures on the Heisenberg group.
 \item {\bf Evolution equations on the Heisenberg group and on related manifolds:} Sections~\ref{ch:evolutionintro}, \ref{ch:phasespace}~and~\ref{ch:evPS}.
 Section \ref{ch:evolutionintro} we set up the left-invariant evolution equations on Gabor transforms. We explain the
rationale behind these convection-diffusion schemes, and we comment on their interpretation in
differential-geometric terms. Sections \ref{ch:phasespace} and \ref{ch:evPS} are concerned with a transfer of the schemes from the full
Heisenberg group to phase space,
resulting in a dimension reduction that is beneficial for
implementation.
\item {\bf Convection:} In Section \ref{ch:reass} we consider convection (reassignment) as an important special case.
 For a suitable choice of Gaussian window, it is possible to exploit
Cauchy-Riemann equations for the analysis of the algorithms, and the design of more efficient
alternatives. For example, in Section \ref{ch:reass} we deduce from these Cauchy-Riemann relations that differential reassignment according to Daudet et al.~\!\cite{Daudet}, boils down to a convection system on Gabor transforms that is equivalent to erosion on the modulus, while preserving the phase.
\item {\bf Discretization and Implementation:} Sections \ref{ch:discrete} and \ref{ch:discreteVF}. In order to derive various suitable algorithms for differential reassignment and diffusion we consider discrete Gabor transforms in Section~\ref{ch:discrete}. As these Gabor transforms are defined on discrete Heisenberg groups, we need left-invariant shift operators for the generators in our left-invariant evolutions on discrete
Heisenberg groups. These discrete left-invariant shifts are derived in Section~\ref{ch:discreteVF}. They are crucial in the
algorithms for left-invariant evolutions on Gabor transforms, as straightforward finite difference approximations of the continuous framework produce considerable errors due to phase oscillations in the Gabor domain.
\item {\bf Implementation and analysis of reassignment:} Sections~\ref{ch:upwind}~and~\ref{ch:chirp}. In Section~\ref{ch:upwind} we employ the results from
the previous Sections in four algorithms for
discrete differential reassignment. We compare these four numerical algorithms by applying them to reassignment of a chirp signal. We provide evidence that it actually works as reassignment (via numerical experiments, subsection \ref{ch:evaluate}) and indeed yields a concentration about the expected curve in the time-frequency plane (Section \ref{ch:chirp}). We show this by deriving the corresponding analytic solutions of reassigned Gabor transforms of arbitrary chirp signals in Section~\ref{ch:chirp}.
\item {\bf Diffusion:}  In Section~\ref{ch:covdiff} we consider signal enhancement via left-invariant diffusion on Gabor transforms. Here we do not apply thresholds on Gabor coefficients. Instead we use both spatial and frequency context of Gabor-atoms in locally adaptive flows in the Gabor domain. We include a basic experiment of enhancement of a noisy 1D-chirp signal. This experiments is intended as a preliminary feasibility study to show possible benefits for various imaging- and signal applications. The benefits may be comparable to our directly related\footnote{Replace $H(2d+1)$ and the Schr\"{o}dinger representation with $SE(d)$ and its left-regular representation on $\mathbb{L}_{2}(\R^{d})$, $d=2,3$.} previous works on enhancement of (multiply crossing) elongated structures in 2D- and 3D medical images via nonlinear adaptive evolutions on invertible orientation scores and/or diffusion-weighted MRI images, \cite{DuitsRTHESIS,FrankenPhDThesis,DuitsAMS1,DuitsAMS2,DuitsR2006AMS,Duits2005IJCV,DuitsIJCV2010,Prckovska,Rodrigues}.
\item {\bf 2D Imaging Applications:} In Section \ref{ch:2D} we investigate extensions of our algorithms to left-invariant evolutions on Gabor transforms of 2-dimensional greyscale images.
    We apply experiments of differential reassignment and texture enhancement on basic images. Finally, we consider a cardiac imaging application where cardiac wall deformations can be directly computed from robust frequency field estimations from Gabor transforms of MRI-tagging images. Our approach by Gabor transforms is inspired by \cite{Arts} and it is relatively simple compared to existing approaches on cardiac deformation (or strain/velocity) estimations, cf.~\!\cite{Duits2011QAM,FlorackDEF,VanAssen2b,vanAssen2007CBfMII}.
\end{itemize}

\section{Gabor transforms and the reduced Heisenberg group \label{ch:Hr}}

Throughout the paper, we fix integers $d \in \mathbb{N}$ and $n \in
\mathbb{Z} \setminus \{ 0 \}$.
The continuous Gabor-transform $\mathcal{G}_{\psi}[f]:\R^{d} \times
\R^{d} \to \mathbb{C}$ of a square integrable signal $f:\R^{d} \to
\mathbb{C}$ is commonly defined as
\begin{equation}\label{austria}
\mathcal{G}_{\psi}[f](p,q)=\int \limits_{\R^{d}} f(\xi)
\overline{\psi(\xi -p)} e^{-2\pi n i\, (\xi-p) \cdot q }\, {\rm d}\xi ~,
\end{equation}
where $\psi \in \mathbb{L}_2(\mathbb{R}^d)$ is a suitable window function. For
window functions centered around zero both
in space and frequency, the Gabor coefficient
$\mathcal{G}_{\psi}[f](p,q)$ expresses the contribution of the frequency
$ n q$ to the behaviour of $f$ near $p$.

This interpretation is suggested by the Parseval formula associated to
the Gabor transform, which reads
\begin{equation} \label{Cpsi}
\int_{\R^d} \int_{\R^d} |\mathcal{G}_{\psi}[f](p,q)|^{2} \; {\rm d}p\,
{\rm d}q %= %\int_{\R^d} \int_{\R^d} |\mathcal{W}_{\psi}[f](p,q,0)|^{2}
= C_\psi  \int_{\R^d} |f(p)|^{2} \; {\rm d}p,~\mbox{ where }
C_\psi =  \frac{1}{n} \|\psi\|^{2}_{\mathbb{L}_{2}(\R^d)}
\end{equation}
for all $f,\psi \in \mathbb{L}_{2}(\R^d)$.
This property can be rephrased as an inversion formula:
\begin{equation}\label{rec}
f(\xi) =\frac{1}{C_\psi}  \int \limits_{\R^d} \int \limits_{\R^d}
\mathcal{G}_{\psi} [f](p,q)\, e^{i 2\pi n(\xi-p) \cdot q} \psi(\xi-p) \, {\rm
d}p  {\rm d}q  ~,
\end{equation}
to be read in the weak sense. The inversion  formula is commonly
understood as the decomposition of $f$ into building
blocks, indexed by a time and a frequency parameter; most applications
of Gabor analysis are based on this heuristic interpretation.
For many such applications, the phase of the Gabor transform is of
secondary importance (see, e.g., the characterization
of function spaces via Gabor coefficient decay \cite{Groechenig2001}).
However, since the Gabor transform uses highly oscillatory
complex-valued functions, its phase information is often crucial,
a fact that has been specifically acknowledged in the context of
reassignment for Gabor transforms \cite{Daudet}.
%and more generally for ``a deeper understanding of the mathematical structure of time frequency shifts'' \cite[ch.9]{Groechenig2001}
%in general
%Now $H_r/\Theta \equiv \R2$ and in principle (if no operators are
%applied to the Gabor transform) the integration over $s$ in %(\ref{rec})
%is redundant as the integrand is independent on $s>0$.

For this aspect of Gabor transform, as for many others\footnote{Regarding the phase factor that arises in the composition of time frequency shifts and the $H_{r}$ group-structure in the Gabor domain we quote ``\emph{This phase factor is absolutely essential for a deeper understanding of the mathematical structure of time frequency shifts, and it is the very reason for a non-commutative group in the analysis.}'' from \cite[Ch:9]{Groechenig2001}.}, the
group-theoretic viewpoint becomes particularly beneficial.
The underlying group is the {\em reduced Heisenberg group} $H_r$. As a
set, $H_r = \mathbb{R}^{2d} \times \mathbb{R}/\mathbb{Z}$, with the
group product
\[
(p,q,s+\mathbb{Z})(p',q',s'+\mathbb{Z})= (p+p', q+q', s+s'+
\frac{1}{2}(q\cdot p'-p \cdot q')+\mathbb{Z})~.
\]
This makes $H_r$ a connected (nonabelian) nilpotent Lie group. The Lie
algebra is spanned by vectors
$A_1,\ldots,A_{2d+1}$ with Lie brackets $[A_i,A_{i+d}] = -A_{2d+1}$, and all other brackets vanishing.

$H_r$ acts on $\mathbb{L}_{2}(\R^d)$ via the
{\em Schr\"odinger representations} $\mathcal{U}^{n}: H_r \to
\mathcal{B}(\mathbb{L}_{2}(\R))$,
\begin{equation}\label{rep}
\mathcal{U}^{n}_{g=(p,q,s+\mathbb{Z})}\psi(\xi)=e^{2\pi i n (s+q \xi
-\frac{pq}{2})}\psi (\xi-p), \qquad \psi \in \mathbb{L}_{2}(\R).
\end{equation}
The associated matrix coefficients are defined as
\begin{equation}
 \label{GT_matrix_coeff}
  \mathcal{W}_{\psi}^{n}f(p,q,s+\mathbb{Z}) =
(\mathcal{U}^{n}_{(p,q,s+\mathbb{Z})}\psi,f)_{\mathbb{L}_{2}(\R^d)}.
\end{equation} In the following, we will often omit the superscript $n$ from $U$ and $\mathcal{W}_\psi$, implicitly
assuming that we use the same choice of $n$ as in the definition of $\mathcal{G}_\psi$. Then a simple comparison of (\ref{GT_matrix_coeff}) with
(\ref{austria}) reveals that
\begin{equation}\label{rel}
\mathcal{G}_{\psi}[f](p,q)=\mathcal{W}_{\psi} f(p,q,s=-\frac{pq}{2}).
\end{equation}
Since
%\[
$ \mathcal{W}_{\psi} f(p,q,s+\mathbb{Z}) = e^{2 \pi i n s}
\mathcal{W}_{\psi} f(p,q,0+\mathbb{Z})$~,
%\]
the phase variable $s$ does not affect the modulus, and (\ref{Cpsi})
can be rephrased as
\begin{equation}\label{Cpsi_matrix_coeff} \int_{0}^1 \int_{\R^d}
\int_{\R^d} |\mathcal{W}_{\psi}[f](p,q,s +\mathbb{Z})|^{2} \; {\rm d}p\,
{\rm d}q {\rm d} s
= C_\psi  \int_{\R^d} |f(p)|^{2} \; {\rm d}p.
\end{equation} Just as before, this induces a weak-sense inversion
formula, which reads
\[
 f = \frac{1}{C_{\psi}} \int_{0}^1
\int_{\R^d} \int_{\R^d} \mathcal{W}_{\psi}[f](p,q,s +\mathbb{Z})
\mathcal{U}^{n}_{(p,q,s+\mathbb{Z})} \psi \; {\rm d}p\, {\rm d}q {\rm d} s~.
\] As a byproduct of (\ref{Cpsi_matrix_coeff}), we note that the
Schr\"odinger representation
is irreducible.
Furthermore, the orthogonal projection $\mathbb{P}_{\psi}$  of
$\mathbb{L}_{2}(H_{r})$ onto the range
$\mathcal{R}(\mathcal{W}_{\psi})$ turns out to be right convolution
with a suitable (reproducing) kernel function,
\[
 (\mathbb{P}_{\psi}U)(h) = U \ast K (h) = \int_{H_r} U(g) K(g^{-1}h) dg~,
\]
with $dg$ denoting the left Haar measure (which is just the Lebesgue measure on $\mathbb{R}^{2d} \times \mathbb{R}/\mathbb{Z}$) and
$
K(p,q,s) =
\frac{1}{C_\psi} \mathcal{W}_{\psi} \psi (p,q,s) = \frac{1}{C_\psi} (U_{(p,q,s)} \psi, \psi)$.
 %\]

The chief reason for choosing the somewhat more redundant function
$\mathcal{W}_{\psi} f$ over
$\mathcal{G}_{\psi}[f]$ is that $\mathcal{W}_\psi$ translates
time-frequency shifts acting on the
signal $f$ to shifts in the argument. If $\mathcal{L}$ and $\mathcal{R}$
denote the left and right regular
representation, i.e., for all $g,h \in H_r$ and $F \in \mathbb{L}_{2}(H_r)$,
\[
 (\mathcal{L}_g F)(h) = F(g^{-1}h)~,~(\mathcal{R}_g F)(h) = F(hg)~,
\] then $\mathcal{W}_{\psi}$ intertwines $\mathcal{U}$ and
$\mathcal{L}$,
\begin{equation} \label{eqn:Gabor_intertwine}
 \mathcal{W}_{\psi} \circ \mathcal{U}_g^n = \mathcal{L}_g \circ
\mathcal{W}_{\psi} ~.
\end{equation}
Thus the group parameter $s$ in $H_r$ keeps track of the phase
shifts induced by the noncommutativity
of time-frequency shifts. By contrast, {\em right shifts} on the Gabor
transform corresponds to changing the window:
\begin{equation}
 \mathcal{R}_g (W_\psi^n (h)) = (\mathcal{U}_{hg} \psi, f) =
\mathcal{W}_{\mathcal{U}_g \psi} f(h)~.
\end{equation}

\section{Left Invariant Evolutions on Gabor Transforms \label{ch:evolutionintro}}

We relate operators $\Phi: \mathcal{R}(\mathcal{W}_{\psi}) \to
\mathbb{L}_{2}(H_{r})$ on Gabor transforms, which actually use and
change the relevant phase information of a Gabor transform, in a
well-posed manner to operators $\Upsilon_{\psi}: \mathbb{L}_{2}(\R^d)
\to \mathbb{L}_{2}(\R^d)$ on signals via
\begin{equation}\label{net}
\hspace{-0.05cm}\mbox{}
\begin{array}{l}
(\Upsilon_{\psi} f)(\xi) =     (\mathcal{W}_{\psi}^{*}  \circ \Phi \circ
\mathcal{W}_{\psi}f)(\xi) \\
 =\frac{1}{C_{\psi}}\int \limits_{[0,1]} \int \limits_{\R^d} \int
\limits_{\R^d} (\Phi(\mathcal{W}_{\psi} f))(p,q,s)\, e^{i 2\pi
n[(\xi,q)+s -(1/2)(p,q)]} \psi(\xi-p) \; {\rm d}p \, {\rm d}q \, {\rm d}s.
\end{array}
\end{equation}
Our aim is to design operators $\Upsilon_\psi$ that address signal
processing problems such as denoising or detection.

\subsection{Design principles}

We now formulate a few desirable properties of $\Upsilon_\psi$, and
sufficient conditions for $\Phi$ to guarantee that
$\Upsilon_\psi$ meets these requirements.
\begin{enumerate}
 \item {\em Covariance with respect to time-frequency-shifts:} The
operator $\Upsilon_{\psi}$ should commute
with time-frequency shifts;
\[
\Upsilon_{\psi} \circ \mathcal{U}_{g}= \mathcal{U}_{g} \circ \Upsilon_{\psi}
\]
for all $g \in H(2d+1)$. This requires a proper treatment of
the phase.

One easy way of guaranteeing covariance of $\Upsilon_{\psi}$ is to
ensure left invariance of $\Phi$:
\begin{equation}\label{Eqcomm}
\Phi \circ \mathcal{L}_{g}=\mathcal{L}_{g} \circ \Phi
\end{equation}
for all $g \in H(2d+1)$. If $\Phi$
commutes with $\mathcal{L}_g$, for all $g \in H_r$, it follows from
(\ref{eqn:Gabor_intertwine}) that
\begin{equation*}
 \Upsilon_{\psi} \circ \mathcal{U}^n_g  = \mathcal{W}_{\psi}^{*}
\circ \Phi \circ \mathcal{W}_{\psi} \circ \mathcal{U}^n_g
 = \mathcal{W}_\psi^* \circ \Phi \circ \mathcal{L}_g \circ
\mathcal{W}_{\psi}
 =  \mathcal{U}^n_g \circ \Upsilon_{\psi} ~.
\end{equation*}
Generally speaking, left invariance of $\Phi$ is not a \emph{necessary}
condition for invariance of $\Upsilon_\psi$:
Note that $\mathcal{W}_\psi^* = \mathcal{W}_\psi^* \circ
\mathbb{P}_{\psi}$. Thus if $\Phi$ is
left-invariant, and $A : \mathbb{L}_2(H_r) \to
\mathcal{R}(\mathcal{W}_{\psi}^n)^\bot$ an arbitrary
operator, then $\Phi + A$ cannot be expected to be left-invariant, but
the resulting operator on the signal
side will be the same as for $\Phi$, thus covariant with respect to
time-frequency shifts.

The authors of \cite{Daudet} studied reassignment
procedures that leave the phase invariant, whereas we shall put emphasis on phase-covariance.
Note however that the two properties are not mutually
exclusive; convection along equiphase lines fulfills both. (See also the discussion in
subsection \ref{subsect:Horizontal}.)
\item {\em Covariance with respect to rotation and translations :}
\begin{equation}
 \Upsilon_{\psi} \circ \mathcal{U}^{SE(d)}_{g}=  \mathcal{U}^{SE(d)}_{g} \circ \Upsilon_{\psi}
\end{equation}
for all $g \in SE(d)=\R^{d} \rtimes SO(d)$ with unitary representation
$\mathcal{U}^{SE(d)}: SE(d) \to B(\mathbb{L}_{2}(\R^{2}))$ given by
$(\mathcal{U}^{SE(d)}_{(x,R)}f)(\xi)= f(R^{-1}(\xi-x))$, for almost every $(x,R) \in SE(d)$.
Rigid body motions on signals and Gabor transforms relate via
\begin{equation} \label{defV}
\begin{array}{ll}
(\mathcal{W}_{\psi} \mathcal{U}^{SE(d)}_{(\ul{x},R)} f)(p,q,s)& =
(\mathcal{V}^{SE(d)}_{(x,R)}\mathcal{W}_{\, \mathcal{U}_{(0,R)}\psi}f)(p,q,s) := \\ &
(\mathcal{W}_{\, \mathcal{U}_{(0,R)}\psi}f)(R^{-1}(p-x),R^{-1}q,s+\frac{pq}{2}),
\end{array}
\end{equation}
for all $f \in \mathbb{L}_{2}(\R^{d})$ and for all $(x,R) \in SE(d)$ and for all $(p,q,s) \in H(2d+1)$ and therefore we require the kernel to be isotropic (besides $\Phi \circ \mathcal{L}_{(x,0,0)}= \mathcal{L}_{(x,0,0)} \circ \Phi$ included in Eq.~(\ref{Eqcomm})) and we require
\begin{equation}\label{rotcovrequire}
\Phi \circ \mathcal{V}^{SE(d)}_{(0,R)}= \mathcal{V}^{SE(d)}_{(0,R)} \circ \Phi
\end{equation}
for all $R \in SO(d)$.
\item {\em Nonlinearity:} The requirement that $\Upsilon_{\psi}$ commute
with $\mathcal{U}^n$ immediately
rules out linear operators $\Phi$. Recall that $\mathcal{U}^n$ is irreducible,
and by Schur's lemma \cite{Dieudonne}, any linear
intertwining operator is a scalar multiple of the identity operator.
\item  By contrast to left invariance, right invariance of $\Phi$ is
undesirable. By a similar argument as for left-invariance, it would
provide that $\Upsilon_{\psi} = \Upsilon_{\mathcal{U}^n_g \psi}$.
\end{enumerate}
We stress that one cannot expect that the processed Gabor transform
$\Phi(\mathcal{W}_\psi f)$ is again the
Gabor transform of a function constructed by the same kernel $\psi$, i.e. we do not expect that
$\Phi(\mathcal{R}(\mathcal{W}_{\psi}^n))
\subset \mathcal{R}(\mathcal{W}_{\psi}^n)$.

\subsection{Invariant differential operators on $H_r$}

The basic building blocks for the evolution equations are the
left-invariant differential operators on $H_r$ of
degree one.
These operators are conveniently obtained by differentiating the {\em
right} regular representation, restricted
to one-parameter subgroups through the generators $\{ A_1,\ldots,A_{2d+1} \}
= \{\partial_{p_{1}}, \ldots,\partial_{p_{d}}, \partial_{q_{1}}, \ldots,\partial_{q_{d}}, \partial_{s}\} \subset T_{e}(H_{r})
$,
\begin{equation} \label{rightrep}
{\rm d}\mathcal{R}(A_{i})U(g) = \lim \limits_{\epsilon \to 0}
\frac{U(ge^{ \epsilon A_i})-U(g)}{\epsilon} \textrm{ for all }g \in H_r
\textrm{ and smooth }U \in \mathcal{C}^{\infty}(H_{r}),
\end{equation}
The resulting differential operators $\{{\rm
d}\mathcal{R}(A_{1}),\ldots,{\rm
d}\mathcal{R}(A_{2d+1})\}=:\{\mathcal{A}_{1},\ldots,\mathcal{A}_{2d+1}\}$
denote the left-invariant vector fields on $H_{r}$, and brief computation of (\ref{rightrep}) yields:
\[
%\[
%\begin{array}{l}
\mathcal{A}_{i}=\partial_{p_{i}} + \frac{q_{i}}{2} \partial_{s}, \ \
\mathcal{A}_{d+i}=\partial_{q_{i}} - \frac{p_{i}}{2} \partial_{s}, \ \
\mathcal{A}_{2d+1}=\partial_{s}, \qquad \textrm{ for }i=1,\ldots,d. \ %& %\\
,
%\end{array}
\]
The differential operators obey the same commutation relations as their
Lie algebra counterparts $A_1,\ldots,A_{2d+1}$
%{\small
\begin{equation} \label{commutators}
[\mathcal{A}_{i},\mathcal{A}_{d+i}]:=\mathcal{A}_{i}\mathcal{A}_{d+i}-\mathcal{A}_{d+i}\mathcal{A}_{i}=-\mathcal{A}_{2d+1},
\end{equation}
and all other commutators are zero. I.e. $\rm d \mathcal{R}$ is a Lie algebra isomorphism.

\subsection{Setting up the equations}

For the effective operator $\Phi$, we will choose left-invariant
evolution operators with stopping time $t>0$. To stress the dependence
on the stopping time we shall write $\Phi_t$ rather than $\Phi$.
Typically, such operators are defined by
$\Phi_{t}(\mathcal{W}_{\psi}f)(p,q,s)=W(p,q,s,t)$ where $W$ is the
solution of %} \\[8pt]
\begin{equation} \label{theeqs}
\boxed{
\begin{array}{l}
\left\{
\begin{array}{ll}
\partial_{t}W(p,q,s,t) &=
Q(|\mathcal{W}_{\psi}f|,\mathcal{A}_{1},\ldots, \mathcal{A}_{2d}) W
(p,q,s,t) , \\
W(p,q,s,0) &=\mathcal{W}_{\psi}f(p,q,s).
\end{array}
\right.
\end{array}
}
\end{equation}
where we note that the left-invariant vector fields
$\{\mathcal{A}_{i}\}_{i=1}^{2d+1}$ on $H_{r}$ are given by
\[
\begin{array}{l}
\mathcal{A}_{i}=\partial_{p_{i}} + \frac{q_{i}}{2} \partial_{s},
\mathcal{A}_{d+i}=\partial_{q_{i}} - \frac{p_{i}}{2} \partial_{s},
\mathcal{A}_{2d+1}=\partial_{s}, \qquad \textrm{ for }i=1,\ldots,d, \ %& %\\
,
\end{array}
\]
with left-invariant quadratic differential form
{\small
\begin{equation}\label{therightchoice}
%\hspace{-0.2cm}\mbox{}
Q(|\mathcal{W}_{\psi}f|,\mathcal{A}_{1},\ldots,\mathcal{A}_{2d})=
-\sum \limits_{i=1}^{2d} a_{i}(|\mathcal{W}_{\psi}f|)(p,q)
\mathcal{A}_{i} + \sum \limits_{i=1}^{2d}\sum
\limits_{j=1}^{2d}\mathcal{A}_{i} \; D_{ij}(|\mathcal{W}_{\psi}f|)(p,q)
\; \mathcal{A}_{j}.
\end{equation}
}
Here $a_{i}(|\mathcal{W}_{\psi} f|)$ and $D_{ij}(|\mathcal{W}_{\psi} f|)$  are functions such that $(p,q) \mapsto
a_{i}(|\mathcal{W}_{\psi}f|)(p,q) \in \R$ and
$(p,q) \mapsto D_{ij}(|\mathcal{W}_{\psi}f|)(p,q) \in \R$ are smooth and either $D=0$ (pure convection) or
$D^{T}=D >0$ holds pointwise (with $D=[D_{ij}]$) for all $i=1,\ldots,2d$, $j=1,\ldots
2d$.
Moreover, in order to guarantee left-invariance, the mappings $a_i : \mathcal{W}_\psi f \mapsto a_{i}(|\mathcal{W}_{\psi} f|)$
need to fulfill the covariance relation
\begin{equation} \label{vi}
a_{i}(|\mathcal{L}_{h}\mathcal{W}_{\psi}f|)(g)=
a_{i}(|\mathcal{W}_{\psi}f|)(p-p',q-q'),
\end{equation}
for all $f \in \mathbb{L}_{2}(\R)$, and all $g=(p,q,s+\mathbb{Z}), h=(p',q',s'+\mathbb{Z}) \in H_{r}$.

For $a_1=\ldots=a_{2d+1}=0$, the equation is a diffusion equation, whereas if
$D=0$, the equation
describes a convection.
We note that existence, uniqueness and square-integrability of the solutions (and thus
well-definedness of $\Upsilon$) are issues that
will have to be decided separately for each particular choice of $a:=(a_{1}, \ldots, a_{2d})^{T}$
and $D$. In general existence and uniqueness are guaranteed,
provided that the convection vector and the diffusion-matrix are smoothly depending on the initial condition, see
Appendix~\ref{ch:unique}. Occasionally,
we shall consider the case where the convection vector $(a_{i})_{i=1}^{2d}$ and the diffusion-matrix $D$ are updated with the absolute value of the current solution $|W(\cdot,t)|$ at time $t>0$,
rather than having them prescribed by the modulus of the initial condition $|W(\cdot,0)|=|\mathcal{W}_{\psi}(f)|=|\mathcal{G}_{\psi}f|$ at time $t=0$, i.e.
\[
\begin{array}{lcl}
a_{i}(|W(\cdot,0)|)  & \textrm{ is sometimes replaced by } & a_{i}(|W(\cdot,t)|) \textrm{ and/or }\\
D_{ij}(|W(\cdot,0)|) & \textrm{ is sometimes replaced by } & D_{ij}(|W(\cdot,t)|)
\end{array}
\]
In case of such replacement the PDE gets non-linear and unique (weak) solutions are not a priori guaranteed. For example in the cases of differential re-assignment we shall consider in Chapter~\ref{ch:reass}, we will restrict ourselves to \emph{viscosity solutions} of the corresponding Hamilton-Jacobi evolution systems, \cite{Evans,Manfredi}.

In order to guarantee rotation covariance we set column vector $a:=(a^{1}, \ldots, a^{2d})^{T}$ and $D=[D_{ij}]$ with row-index $i$ and column-index $j$ and we require
\begin{equation} \label{rotcov}
\begin{array}{l}
(a( \mathcal{V}_{(\ul{0},R)}^{SE(d)}W(\cdot,t)))(g)= \ul{R} \; (a(W(\cdot,t)))(\ul{R}^{-1}g)\ , \\
(D( \mathcal{V}_{(\ul{0},R)}^{SE(d)}W(\cdot,t)))(g)= \ul{R}^{-1} \; (D(W(\cdot,t)))(\ul{R}^{-1}g) \; \ul{R}\ .
\end{array}
\end{equation}
for all $\ul{R}=R \otimes R \in SO(2d)$, $R \in SO(d)$, $g \in H(2d+1)$, $U \in \mathbb{L}_{2}(H(2d+1))$, where we recall (\ref{defV}).

This definition of $\Phi_t$, for each $t>0$ fixed, satisfies the criteria we set up above:
\begin{enumerate}
 \item Since the evolution equation is left-invariant (and provided uniqueness
of the solutions), it follows that
$\Phi_t$ is left-invariant. Thus the associated $\Upsilon_\psi$ is
invariant under time-frequency shifts.
\item The rotated left-invariant gradient transforms as follows
\[
\left(\underline{\mathcal{A}}_{g} \mathcal{V}_{(\ul{0},R)}^{SE(d)} W(\cdot,t)\right)(g)=
\ul{R} \; \left(\underline{\mathcal{A}}_{\ul{R}^{-1}g} W(\cdot,t)\right)(\ul{R}^{-1}g), \textrm{ with }\underline{\mathcal{A}}=(\mathcal{A}_1, \ldots, \mathcal{A}_{2d})^{T}.
\]
Thereby (the generator of) our diffusion operator $\Phi_{t}$ is rotation covariant, i.e. \[
Q(|W(\cdot,t)|,\underline{\mathcal{A}}) \circ \mathcal{V}_{(\ul{0},R)}^{SE(d)}= \mathcal{V}_{(\ul{0},R)}^{SE(d)} \circ Q(|W(\cdot,t)|,\underline{\mathcal{A}}) \textrm{ for all }R \in SO(d),
\]
if Eq.~(\ref{rotcov}) and Eq.~(\ref{vi}) hold. For example, if $a=0$ and $D$ would be constant, then by Eq.~(\ref{rotcov}) and Schur's lemma one has {\small $D=\textrm{diag}\{D^{11},\ldots, D^{11}, D^{d+1,d+1}, \ldots, D^{d+1,d+1}\}$} yielding the Kohn's Laplacian $\Delta_{K}= D_{11}\sum_{i=1}^{d} \mathcal{A}_{i}^{2}+
D_{d+1,d+1} \sum_{i=d+1}^{2d} \mathcal{A}_{i}^{2}$, cf.~\cite{Gaveau}, and indeed $\Delta_{K} \circ \mathcal{V}_{(\ul{0},R)}^{SE(d)}=\mathcal{V}_{(\ul{0},R)}^{SE(d)} \circ \Delta_{K}$.
 \item In order to ensure non-linearity, not all of the functions $a_i$,
$D_{ij}$ should be constant, i.e. the schemes should be
\emph{adaptive convection} and/or \emph{adaptive diffusion}, via
adaptive choices of convection vectors $(a_1,\ldots,a_{2d})^T$ and/or conductivity matrix $D$. We will use ideas similar to our
previous work on adaptive diffusions on invertible orientation scores
\cite{Fran08}, \cite{DuitsAMS1}, \cite{Duits2005IJCV}, \cite{DuitsAMS2}.
%(where we employed evolution equations for the 2D-Euclidean motion
%group).
We use the absolute value of the (evolving) Gabor transform to adapt the diffusion and convection
in order to avoid oscillations.
 \item The two-sided invariant differential operators of degree one
correspond to the center of the Lie algebra,
which is precisely the span of $A_{2d+1}$. Both in the cases of
diffusion and convection, we consistently removed the
$\mathcal{A}_{2d+1} = \partial_s$-direction, and we removed the
$s$-dependence in the coefficients $a_{i}(|\mathcal{W}_{\psi}f|)(p,q)$,
$D_{ij}(|\mathcal{W}_{\psi}f|)(p,q)$ of the generator
$Q(|\mathcal{W}_{\psi}f|,\mathcal{A}_{1},\ldots, \mathcal{A}_{2d})$ by
taking the absolute value $|\mathcal{W}_{\psi}f|$, which is independent
of $s$. %in the adaptation.%\footnote{This was also essential in the
%context of non-linear left-invariant diffusions on invertible
%orientation scores, %\cite{DuitsAMS2}, \cite{Fran08}.}
A more complete discussion of the role of the $s$-variable is contained
in the following subsection.
\end{enumerate}

\subsection{Convection and Diffusion along Horizontal Curves}
\label{subsect:Horizontal}

So far our motivation for (\ref{theeqs}) has been group theoretical. There is one issue we did not address yet,
namely the omission of $\partial_{s}=\mathcal{A}_{2d+1}$ in (\ref{theeqs}). Here we first motivate this omission and then consider the differential geometrical consequence that (adaptive) convection and diffusion takes place along so-called horizontal curves.

The reason for the removal of the $\mathcal{A}_{2d+1}$ direction in our diffusions and convections is simply that this direction leads to a scalar multiplication operator mapping the space of Gabor transform to itself, since
%\[
$\partial_{s} \mathcal{W}_{\psi}f = - 2\pi i n \, \mathcal{W}_{\psi}f$ .
%\]
Moreover, we adaptively steer the convections and diffusions by the modulus of a Gabor transform $|\mathcal{W}_{\psi}f(p,q,s)|=|\mathcal{G}_{\psi}f(p,q)|$, which is independent of $s$, and clearly a vector field $(p,q,s) \mapsto F(p,q)\partial_{s}$ is left-invariant iff $F$ is constant. Consequently it does not make sense to include the separate $\partial_{s}$ in our convection-diffusion equations, as it can only yield a scalar multiplication. Indeed, for all constant $\alpha>0, \beta \in \R$ we have
\[
\begin{array}{l} \,
[\partial_{s},Q(|\mathcal{W}_{\psi}f|, \mathcal{A}_{1},\ldots, \mathcal{A}_{2d})]=0 \textrm{ and } \partial_{s} \mathcal{W}_{\psi}f = - 2\pi i n \, \mathcal{W}_{\psi}f  \Rightarrow \\[6pt]
e^{t \left((\alpha \partial_{s}^2 +\beta \partial_{s}) + Q(|\mathcal{W}_{\psi}f|, \mathcal{A}_{1},\ldots, \mathcal{A}_{2d})\right) }=
e^{-t \alpha(2\pi n)^2 -t \beta 2\pi i n}\, e^{t Q(|\mathcal{W}_{\psi}f|, \mathcal{A}_{1},\ldots, \mathcal{A}_{2d})}.
\end{array}
\]
In other words $\partial_{s}$ is a redundant direction in each tangent space $T_{g}(H_{r})$, $g \in H_{r}$. This however does not imply that it is a redundant direction in the group manifold $H_{r}$ itself, since clearly the $s$-axis represents the relevant phase and stores the non-commutative nature between position and frequency, \cite[ch:1]{DuitsFuehrJanssen}.

The omission of the redundant direction $\partial_{s}$ in $T(H_{r})$ has an important geometrical consequence.
Akin to our framework of linear evolutions on orientation scores, cf. \cite{DuitsAMS2,Fran08}, this means that we enforce \emph{horizontal} diffusion and convection. In other words, the generator of our evolutions will only include derivations within the horizontal part of the Lie algebra spanned by $\{\mathcal{A}_{1}, \mathcal{A}_{2}\}$. On the Lie group $H_{r}$
this means that transport and diffusion only takes place along so-called \emph{horizontal curves} in $H_{r}$ which are curves $t \mapsto (p(t),q(t),s(t)) \in H_{r}$, with $s(t) \in (0,1)$, along which
\begin{equation} \label{SST}
s(t)= \frac{1}{2}\int \limits_{0}^{t} \sum \limits_{i=1}^{d} q_{i}(\tau)p_{i}'(\tau)-p_{i}(\tau)q_{i}'(\tau) {\rm d}\tau \ ,
\end{equation}
see Theorem \ref{th:hor}.
This gives a nice geometric interpretation to the phase variable $s(t)$, as by the Stokes theorem it represents the net surface area between a straight line connection between $(p(0),q(0),s(0))$ and $(p(t),q(t),s(t))$ and the actual horizontal curve connection
\mbox{$[0,t]\ni \tau \mapsto (p(\tau),q(\tau),s(\tau))$}. For example, horizontal diffusion with diagonal $D$ is the forward Kolmogorov equation of Brownian motion $t \mapsto (P(t),Q(t))$ in position and frequency and Eq.~(\ref{SST}) associates a random variable $S(t)$ (measuring the net surface area) to the implicit smoothing in the phase direction due to the commutator $[\mathcal{A}_{1},\mathcal{A}_{2}]=\mathcal{A}_{3}=\partial_{s}$, cf.~\cite{Gaveau,DuitsFuehrJanssen,DuitsAMS1}.

In order to explain why the omission of the redundant direction $\partial_{s}$ from the tangent bundle $T(H_{r})$ implies a restriction to horizontal curves, we consider the dual frame associated to our frame of reference  $\{\mathcal{A}_{1},\ldots,\mathcal{A}_{2d+1}\}$. We will denote this dual frame by $\{{\rm d}\mathcal{A}^{1},\ldots, {\rm d}\mathcal{A}^{2d+1}\}$ and it is uniquely determined by $\langle {\rm d}\mathcal{A}^{i}, \mathcal{A}_{j} \rangle=\delta^{i}_{j}$, $i,j=1,2,3$ where $\delta^{i}_{j}$ denotes the Kronecker delta. A brief computation yields
\begin{equation}\label{dualbasis}
\begin{array}{l}
\left.{\rm d}\mathcal{A}^{i} \right|_{g=(p,q,s)}= {\rm d}p^{i} \ , \
\left. {\rm d}\mathcal{A}^{d+i} \right|_{g=(p,q,s)}= {\rm d}q^{i} \ , \ \ i=1,\ldots,d \\
\left. {\rm d}\mathcal{A}^{2d+1} \right|_{g=(p,q,s)}= {\rm d}s + \frac{1}{2} (p \cdot {\rm d}q-q \cdot {\rm d}p),
\end{array}
\end{equation}
Consequently a smooth curve $t \mapsto \gamma(t)=(p(t),q(t),s(t))$ is horizontal iff
\[
\left.\langle {\rm d}\mathcal{A}^{2d+1} \right|_{\gamma(s)}, \gamma'(s) \rangle=0 \desda s'(t)=\frac{1}{2}(q(t) \cdot p'(t)- p(t) \cdot q'(t)).
\]
%$
%[\mathcal{A}_{i1},\mathcal{A}_{i+d}]=-\mathcal{A}_{2d+1}$, $i=1,\ldots,d$.
\begin{theorem}\label{th:hor}
Let $f \in \mathbb{L}_{2}(\R)$ be a signal and $\mathcal{W}_{\psi}f$ be its Gabor transform associated to the Schwartz function $\psi$. If we just consider convection and no diffusion (i.e. $D=0$) then the solution of (\ref{theeqs}) is given by
\[
W(g,t)= \mathcal{W}_{\psi}f( \gamma^{g}_{f}(t))\ , \qquad g=(p,q,s) \in H_r,
\]
where the characteristic horizontal curve $t \mapsto \gamma^{g_0}_{f}(t)=(p(t),q(t),s(t))$ for each $g_0=(p_0,q_0,s_0) \in H_r$ is given by the unique solution of the following ODE:
\[
\left\{
\begin{array}{ll}
\dot{p}(t)= -a^{1}(|\mathcal{W}_{\psi}f|)(p(t),q(t)), & p(0)=p_0,\\
\dot{q}(t)= -a^{2}(|\mathcal{W}_{\psi}f|)(p(t),q(t)), & q(0)=q_0, \\
\dot{s}(t)= \frac{q(t)}{2} \dot{p}(t)- \frac{p(t)}{2} \dot{q}(t), &s(0)=s_{0}, \\
\end{array}
\right.
\]
Consequently, the operator $\mathcal{W}_{\psi}f \mapsto W(\cdot,t)$
is phase covariant (the phase moves along with the characteristic curves of transport):
\begin{equation} \label{phasecovariant}
\textrm{arg}\{W(g,t)\}= \textrm{arg}\{\mathcal{W}_{\psi}f\}(\gamma_{f}^{g}) \textrm{ for all }t>0.
\end{equation}
\end{theorem}
\textbf{Proof }
First we shall show that $g \, \gamma^{e}_{\mathcal{U}_{g^{-1}}f}(t)=\gamma^{g}_{f}(t)$ for all $g \in H_{r}$ and all $t>0$ and all $f \in \mathbb{L}_{2}(\R)$. To this end we note that both solutions are horizontal curves, i.e.
\[
\langle \left.{\rm d} \mathcal{A}^{3}\right|_{\gamma^{g}_{f}(t)}, \dot{\gamma}^{g}_{f}(t) \rangle = \langle \left.{\rm d}\mathcal{A}^{3}\right|_{g \, \gamma^{e}_{\mathcal{U}_{g^{-1}}f}(t)}, g \, \dot{\gamma}^{e}_{\;\mathcal{U}_{g^{-1}}f}\,(t) \rangle=0,
\]
where ${\rm d}\mathcal{A}^{3}={\rm d}s+ 2(p {\rm d}q-q {\rm d}p)$. So it is sufficient to check whether the first two components of the curves coincide. Let $g=(p,q,s) \in H_{r}$ and define $p_{e}:\R^{+} \to \R$, $q_{e}:\R^{+} \to \R$ and $p_{g}:\R^{+} \to \R$, $q_{g}:\R^{+} \to \R$ by
\[
\begin{array}{ll}
p_{e}(t):= \langle {\rm d}p, g \, \dot{\gamma}^{e}_{\;\mathcal{U}_{g^{-1}}f}\,(t) \rangle \ , &  p_{g}(t):= \langle {\rm d}p, \dot{\gamma}^{g}_{\;f}\,(t) \rangle \ , \\
q_{e}(t):= \langle {\rm d}q, g \, \dot{\gamma}^{e}_{\;\mathcal{U}_{g^{-1}}f}\,(t) \rangle  \ , & q_{g}(t):= \langle {\rm d}q, \dot{\gamma}^{g}_{\;f}\,(t) \rangle \ ,
\end{array}
\]
then it remains to be shown that $p_{e}+p=p_{g}$ and $p_{e}+q=q_{g}$. To this end we compute
\[
\begin{array}{l}
\frac{dp_{e}}{dt}(t) = a^{1}(|\mathcal{W}_{\psi}\mathcal{U}_{g^{-1}}f|)(p_{e}(t),q_{e}(t)) = a^{1}(|\mathcal{L}_{g^{-1}}\mathcal{W}_{\psi}f|)(p_{e}(t),q_{e}(t)) = a^{1}(|\mathcal{W}_{\psi}f|)(p_{e}(t)+p, q_{e}(t)+q),
\end{array}
\]
so that we see that $(p_{e}+p,q_{e}+q)$ satisfies the following ODE system:
\[
\left\{
\begin{array}{ll}
\frac{d}{dt} (p+p_{e})(t) &=a^{1}(|\mathcal{W}_{\psi}f|)\, (p_{e}(t)+p, q_{e}(t)+q), \qquad t>0, \\
\frac{d}{dt} (q+q_{e})(t) &=a^{2}(|\mathcal{W}_{\psi}f|)\, (q_{e}(t)+q, q_{e}(t)+q), \qquad t>0, \\
p+p_{e}= p, \\
q+q_{e}=q.
\end{array}
\right.
\]
This initial value problem has a unique smooth solution, so indeed $p_{g}=p+p_{e}$ and $q_g=q+q_e$.
Finally, we have by means of the chain-rule for differentiation:
\[
\begin{array}{ll}
\frac{d}{dt}(\mathcal{W}_{\psi}f)(\gamma^{g}_{f}(t))&= \sum \limits_{i=1}^{2}
\langle \left. {\rm d}\mathcal{A}^{i} \right|_{\gamma^{g}_{f}(t)},\dot{\gamma}^{g}_{f}(t)\rangle \; \left. \left(\mathcal{A}_{i}\right|_{\gamma^{g}_{f}(t)} \mathcal{W}_{\psi}f\right)(\gamma^{g}_{f}(t))  \\
&= \dot{p}_{g}(t) \; \left. \mathcal{A}_{1}\right|_{\gamma^{g}_{f}(t)} \mathcal{W}_{\psi}f(\gamma^{g}_{f}(t)) +
 \dot{q}_{g}(t) \; \left. \mathcal{A}_{2}\right|_{\gamma^{g}_{f}(t)} \mathcal{W}_{\psi}f(\gamma^{g}_{f}(t)) \\
  &= -a^{1}(|\mathcal{W}_{\psi}f|)(p_{g}(t),q_{g}(t)) \left.\mathcal{A}_{1}\right|_{\gamma^{g}_{f}(t)} \mathcal{W}_{\psi}(\gamma^{g}_{f}(t))\\ &\qquad -a^{2}(|\mathcal{W}_{\psi}f|)(p_{g}(t),q_{g}(t))  \left.\mathcal{A}_{2} \right|_{\gamma^{g}_{f}(t)} \mathcal{W}_{\psi}(\gamma^{g}_{f}(t))\ .
\end{array}
\]
from which the result follows $ \hfill \Box$ \\ \\
%\end{proof}
Also for the (degenerate) diffusion case with $D=D^{T}=[D_{ij}]_{i,j=1,\ldots,d}>0$, the omission of the \mbox{$2d\!+\!1^{\textrm{th}}$} direction $\partial_{s}=\mathcal{A}_{2d+1}$ implies that diffusion takes place along horizontal curves. Moreover, the omission does not affect the smoothness and uniqueness of the solutions of (\ref{theeqs}), since the initial condition is infinitely differentiable (if $\psi$ is a Schwarz function) and the H\"{o}rmander condition \cite{Hoermander}, \cite{DuitsR2006SS2} is by (\ref{commutators}) still satisfied. %since clearly

The removal of the $\partial_{s}$ direction from the tangent space does not imply that one can entirely ignore the $s$-axis in the domain of a (processed) Gabor transform. The domain of a (processed) Gabor transform $\Phi_{t}(\mathcal{W}_{\psi}f)$ should \emph{not}\footnote{As we explain in \cite[App. B and App. C ]{DuitsFuehrJanssen} the Gabor domain is a principal fiber bundle $P_{T}=(H_r,\mathbb{T},\pi,\mathcal{R})$ equipped with the Cartan connection form $\omega_{g}(X_g)= \langle {\rm d}s + \frac{1}{2}(p {\rm d}q -q {\rm d}p) , X_{g}\rangle$, or equivalently, it is a contact manifold, cf. \cite[p.6]{Bryantbook}, \cite[App. B, def. B.14]{DuitsFuehrJanssen} , $(H_{r},{ \rm d}\mathcal{A}^{2d+1}))$.}
 be considered as $\R^{2d} \equiv H_{r}/\Theta$. Simply, because $[\partial_{p},\partial_{q}]=0$ whereas we should have (\ref{commutators}).
%It is not even entirely appropriate to %consider it as the reduced Heisenberg group $H_{r}=H_r/\{0\}\times \{0\}\times \mathbb{Z}$, because of the horizontality constraint.
For further differential geometrical details see the appendices of \cite{DuitsFuehrJanssen}, analogous to the differential geometry on orientation scores, \cite{DuitsAMS2}, \cite[App. D , App. C.1 ]{DuitsFuehrJanssen}.

\section{Towards Phase Space and Back \label{ch:phasespace}}

As pointed out in the introduction it is very important to keep track of the phase variable $s>0$. The first concern that arises here is whether this results in slower algorithms. In this section we will show that this is not the case. As we will explain next, one can use an invertible mapping $\mathcal{S}$ from the space $\mathcal{H}_{n}$ of Gabor transforms to phase space (the space of Gabor transforms restricted to the plane $s=\frac{pq}{2}$). As a result by means of conjugation with $\mathcal{S}$ we can map our diffusions
on $\mathcal{H}_{n}\subset \mathbb{L}_{2}(\R^2 \times [0,1])$ uniquely to diffusions on $\mathbb{L}_{2}(\R^2)$ simply by conjugation with $\mathcal{S}$. From a geometrical point of view it is easier to consider the diffusions on $\mathcal{H}_{n}\subset \mathbb{L}_{2}(\R^{2d} \times [0,1])$ than on $\mathbb{L}_{2}(\R^{2d})$, even though all our numerical PDE-Algorithms take place in phase space in order to gain speed.
\begin{definition}
Let $\mathcal{H}_{n}$ denote the space of all complex-valued functions $F$ on $H_{r}$ such that
$F(p,q,s+\mathbb{Z})=e^{-2 \pi i ns}F(p,q,0)$ and $F(\cdot,\cdot,s+\mathbb{Z}) \in \mathbb{L}_{2}(\R^{2d})$ for all $s\in \mathbb{R}$. Clearly $\mathcal{W}_{\psi}f \in \mathcal{H}_{n}$ for all $f,\psi \in \mathcal{H}_{n}$.
\end{definition}
In fact $\mathcal{H}_{n}$ is the closure of the space $\{\mathcal{W}^n_{\psi}f \; |\; \psi,f \in \mathbb{L}_{2}(\R)\}$ in $\mathbb{L}_{2}(H_{r})$.  The space $\mathcal{H}_{n}$ is bi-invariant, since:
\begin{equation}\label{important}
\begin{array}{ll}
\mathcal{W}_{\psi}^{n} \circ \mathcal{U}^{n}_{g}=\mathcal{L}_{g} \circ \mathcal{W}_{\psi}^{n}  \textrm{ and }
\mathcal{W}^{n}_{\mathcal{U}^{n}_{g}\psi} = \mathcal{R}_{g} \circ \mathcal{W}^{n}_{\psi},
\end{array}
\end{equation}
where $\mathcal{R}$ denotes the right regular representation on $\mathbb{L}_{2}(H_{r})$ and $\mathcal{L}$ denotes the left regular representation of $H_{r}$ on $\mathbb{L}_{2}(H_{r})$. We can identify $\mathcal{H}_{n}$ with $\mathbb{L}_{2}(\R^{2d})$ by means of the following operator $\mathcal{S}:\mathcal{H}_{n} \to \mathbb{L}_{2}(\R^{2d}) $ given by
\begin{equation} \label{S}
(\mathcal{S}F)(p,q)= F(p,q, \frac{pq}{2}+\mathbb{Z})=e^{i \pi n p q} F(p,q,0+\mathbb{Z}).
\end{equation}
Clearly, this operator is invertible and its inverse is given by
\begin{equation} \label{Smin1}
(\mathcal{S}^{-1}F)(p,q,s+\mathbb{Z})= e^{-2 \pi i sn} e^{-i \pi n pq} F(p,q)
\end{equation}
The operator $\mathcal{S}$ simply corresponds to taking the section $s(p,q)=-\frac{pq}{2}$ in the left cosets $H_r/\Theta$ where $\Theta=\{(0,0,s+\mathbb{Z})\; |\; s \in \R\}$ of $H_r$. Furthermore we recall the common Gabor transform $\mathcal{G}_{\psi}^{n}$ given by (\ref{austria}) and its relation (\ref{rel}) to the full Gabor transform, which we can now write as $\mathcal{G}_{\psi}^{n}= \mathcal{S} \circ \mathcal{W}_{\psi}^{n}$.
\begin{theorem}\label{th:one}
Let the operator $\Phi$ map the closure $\mathcal{H}_{n}$, $n \in \mathbb{Z}$, of the space of Gabor transforms into itself, i.e. $\Phi: \mathcal{H}_{n} \to \mathcal{H}_{n}$. Define the left and right-regular rep's of $H_{r}$ on $\mathcal{H}_{n}$ by restriction
\begin{equation}\label{defreptilde}
\mathcal{R}^{(n)}_{g}=\left. \mathcal{R}_{g}\right|_{\mathcal{H}_{n}} \textrm{ and }
\mathcal{L}^{(n)}_{g}=\left. \mathcal{L}_{g}\right|_{\mathcal{H}_{n}}\, \quad \textrm{ for all }g \in H_{r}.
\end{equation}
Define the corresponding left and right-regular rep's of $H_{r}$ on phase space by
\[
\begin{array}{ll}
\tilde{\mathcal{R}}^{(n)}_g := \mathcal{S} \circ \mathcal{R}^{(n)}_g \circ \mathcal{S}^{-1}, &  %\\
\tilde{\mathcal{L}}^{(n)}_g := \mathcal{S} \circ \mathcal{L}^{(n)}_g \circ \mathcal{S}^{-1}.
\end{array}
\]
For explicit formulas see \cite[p.9]{DuitsFuehrJanssen}. Let
$\tilde{\Phi}:= \mathcal{S}\circ \Phi \circ \mathcal{S}^{-1}$ be the corresponding operator on $\mathbb{L}_{2}(\R^{2d})$ and
\[
\Upsilon_{\psi} = (\mathcal{W}_{\psi}^{n})^{*} \circ \Phi \circ \mathcal{W}_{\psi}^{n} = (\mathcal{S} \mathcal{W}_{\psi}^{n})^{-1} \circ \tilde{\Phi} \circ \mathcal{S} \mathcal{W}_{\psi}^{n}= (\mathcal{G}_{\psi}^{n})^{*} \circ \tilde{\Phi} \circ \mathcal{G}_{\psi}^{n}.
\]
Then one has the following correspondence:
\begin{equation}\label{correspondence}
\Upsilon_{\psi} \circ \mathcal{U}^{n}=\mathcal{U}^{n} \circ \Upsilon_{\psi} \Leftarrow \Phi \circ \mathcal{L}^{n} =\mathcal{L}^{n} \circ \Phi \desda \tilde{\Phi} \circ \tilde{\mathcal{L}}^{n} = \tilde{\mathcal{L}}^{n} \circ \tilde{\Phi}.
\end{equation}
If moreover $\Phi(\mathcal{R}(\mathcal{W}_{\psi})) \subset \mathcal{R}(\mathcal{W}_{\psi})$ then the left implication may be replaced by an equivalence. If $\Phi$ does not satisfy this property then one may replace $\Phi \to \mathcal{W}_{\psi}\mathcal{W}_{\psi}^{*}\Phi$ in (\ref{correspondence}) to obtain full equivalence. Note that $\Upsilon_{\psi}=\mathcal{W}_{\psi}^{*}\Phi \mathcal{W}_{\psi}=
\mathcal{W}_{\psi}^{*}(\mathcal{W}_{\psi}\mathcal{W}_{\psi}^{*}\Phi) \mathcal{W}_{\psi}$.
\end{theorem}
\textbf{Proof }
For details see our technical report \cite[Thm 2.2]{DuitsFuehrJanssen}.
%\end{proof}
\section{Left-invariant Evolutions on Phase Space \label{ch:evPS}}

For the next three chapters, for the sake of simplicity, we fix $d=1$ and consider left-invariant evolutions on Gabor transforms of 1D-signals. We will return to the case $d=2$ in Section \ref{ch:2D} where besides of some extra bookkeeping the rotation-covariance (2nd design principle) comes into play.
%, in order to obtain manageable expressions.

We want to apply Theorem \ref{th:one} to our left-invariant evolutions (\ref{theeqs}) to obtain the left-invariant diffusions on phase space (where we reduce 1 dimension in the domain). To this end we first compute the left-invariant vector fields
$\{\tilde{\mathcal{A}}_{i}\}:=\{\mathcal{S} \mathcal{A}_{i} \mathcal{S}^{-1}\}_{i=1}^{3}$
on phase space.
The left-invariant vector fields on phase space are
{\small
\begin{equation}\label{phasespacegen}
\hspace{-0.13cm}\mbox{}
\begin{array}{l}
\tilde{\mathcal{A}}_{1}U(p',q')=\mathcal{S} \mathcal{A}_{1}\mathcal{S}^{-1}U(p',q')= ((\partial_{p'}\!-\! 2n \pi i q' ) U)(p',q'), \\
\tilde{\mathcal{A}}_{2}U(p',q')=\mathcal{S} \mathcal{A}_{2}\mathcal{S}^{-1}U(p',q')=(\partial_{q'}U)(p',q'), \\
\tilde{\mathcal{A}}_{3}U(p',q')=\mathcal{S} \mathcal{A}_{3}\mathcal{S}^{-1}U(p',q')= -2 i n \pi U(p',q')\ ,
\end{array}
\end{equation}
}
for all $(p,q)\in \R$ and all locally defined smooth functions $U: \Omega_{(p,q)} \subset \R^2\to \mathbb{C}$.
%\bibliographystyle{ieeetr}
%\begin{small}
%\bibliography{strings,literature,root,crossref}
%\end{small}
%\subsection{Left-Invariant Evolution Equations on Phase Space}

Now that we have computed the left-invariant vector fields on phase space, we can express our left-invariant evolution equations (\ref{theeqs}) on phase space
\begin{equation} \label{theeqsPS}
\boxed{
\begin{array}{l}
\left\{
\begin{array}{ll}
\partial_{t}\tilde{W}(p,q,t) &= \tilde{Q}(|\mathcal{G}_{\psi}f|,\tilde{\mathcal{A}}_{1}, \tilde{\mathcal{A}}_{2}) \tilde{W} (p,q,t) , \\
\tilde{W}(p,q,0) &=\mathcal{G}_{\psi}f(p,q).
\end{array}
\right.
\end{array}
}
\end{equation}
with left-invariant quadratic differential form
{\small
\begin{equation}\label{therightchoicePS}
\tilde{Q}(|\mathcal{G}_{\psi}f|,\tilde{\mathcal{A}}_{1}, \tilde{\mathcal{A}}_{2})=
-\sum \limits_{i=1}^{2} a_{i}(|\mathcal{G}_{\psi}f|)(p,q) \tilde{\mathcal{A}}_{i} + \sum \limits_{i=1}^{2}\sum \limits_{j=1}^{2} \tilde{\mathcal{A}}_{i} \; D_{ij}(|\mathcal{G}_{\psi}f|)(p,q) \; \tilde{\mathcal{A}}_{j}.
\end{equation}
}
Similar to the group case, the $a_{i}$ and $D_{ij}$ are functions such that $(p,q) \mapsto a_{i}(|\mathcal{G}_{\psi}f|)(p,q) \in \R$ and
$(p,q) \mapsto a_{i}(|\mathcal{G}_{\psi}f|)(p,q) \in \R$ are smooth and either $D=0$ (pure convection) or $D^{T}=D>0$ (with $D=[D_{ij}]$ $i,j=1,\ldots,2d$), so H\"{o}rmander's condition \cite{Hoermander} (which guarantees smooth solutions $\tilde{W}$, provided the initial condition $\tilde{W}(\cdot,\cdot,0)$ is smooth) is satisfied because of (\ref{commutators}).
\begin{theorem}\label{th:correspondence}
The unique solution $\tilde{W}$ of (\ref{theeqsPS}) is obtained from the unique solution $W$ of (\ref{theeqs}) by means of
\[
\tilde{W}(p,q,t)= (\mathcal{S}\, W(\cdot,\cdot,\cdot,t))(p,q)\ , \textrm{ for all }t\geq 0 \textrm{ and for all }(p,q) \in \R^2,
\]
with in particular $\tilde{W}(p,q,0)=\mathcal{G}_{\psi}(p,q)=(\mathcal{S} \mathcal{W}_{\psi})(p,q)=(\mathcal{S}W(\cdot,\cdot,\cdot,0))(p,q)$ .
\end{theorem}
\textbf{Proof }
This follows by the fact that the evolutions (\ref{theeqs}) leave the function space $\mathcal{H}_{n}$ invariant and the fact that the evolutions (\ref{theeqsPS}) leave the space $\mathbb{L}_{2}(\R^2)$ invariant, so that we can apply direct conjugation with the invertible operator $\mathcal{S}$ to relate the unique solutions, where we have
\begin{equation} \label{conj}
\begin{array}{ll}
\tilde{W}(p,q,t) &= (e^{t \tilde{Q}(|\mathcal{G}_{\psi}f|,\tilde{\mathcal{A}}_{1}, \tilde{\mathcal{A}}_{2})} \mathcal{G}_{\psi}f)(p,q)\\ & = (e^{t \tilde{Q}(|\mathcal{G}_{\psi}f|, \mathcal{S} \mathcal{A}_{1}  \mathcal{S}^{-1},  \mathcal{S} \mathcal{A}_{2} \mathcal{S}^{-1})} \mathcal{S}\mathcal{W}_{\psi}f)(p,q) \\
 &=  (e^{\mathcal{S}\, \circ \, t \, Q(|\mathcal{W}_{\psi}f|,  \mathcal{A}_{1}  ,  \mathcal{A}_{2} )\, \circ \,\mathcal{S}^{-1}} \mathcal{S}\mathcal{W}_{\psi}f)(p,q) \\&= (\mathcal{S}\, \circ \, e^{ t \, Q(|\mathcal{W}_{\psi}f|,  \mathcal{A}_{1}  ,  \mathcal{A}_{2} )} \, \circ \,\mathcal{S}^{-1} \mathcal{S}\circ \mathcal{W}_{\psi}f)(p,q) \\ &= (\mathcal{S}\, W(\cdot,\cdot,\cdot,t))(p,q)
\end{array}
\end{equation}
for all $t>0$ on densely defined domains. For every $\psi \in \mathbb{L}_{2}(\R)\cap S(\R)$, the space of Gabor transforms is a reproducing kernel space with a bounded and smooth reproducing kernel, so that $\mathcal{W}_{\psi}f$ (and thereby $|\mathcal{W}_{\psi}f|=|\mathcal{G}_{\psi}f|$) is uniformly bounded and continuous and equality (\ref{conj}) holds for all $p,q \in \mathbb{R}^2$. $\hfill \Box$
%\end{proof}
\subsection{The Cauchy Riemann Equations on Gabor Transforms and the Underlying Differential Geometry}

As previously observed in \cite{Daudet}, the Gabor transforms associated to Gaussian windows obey Cauchy-Riemann equations which are particularly useful for the analysis of convection schemes, as well as for the design of more efficient algorithms.
More precisely, if the window is a Gaussian $\psi(\xi)=\psi_{a}(\xi):=e^{-\pi n \frac{(\xi-c)^2}{a^2}}$ and $f$ is some arbitrary signal in $\mathbb{L}_{2}(\R)$ then we have
\begin{equation} \label{CR1q}
\begin{array}{l}
(a^{-1}\mathcal{A}_{2}+ i a \mathcal{A}_{1})\mathcal{W}_{\psi}(f)=0 \desda
(a^{-1}\mathcal{A}_{2}+ i a \mathcal{A}_{1})\log \mathcal{W}_{\psi}(f)=0
\end{array}
\end{equation}
where we included a general scaling $a>0$. On phase space this boils down to
\begin{equation} \label{CR1}
\begin{array}{l}
(a^{-1}\tilde{\mathcal{A}}_{2}+ i a \tilde{\mathcal{A}}_{1})\mathcal{G}_{\psi}(f)=0
%\desda
%(a^{-1}\tilde{\mathcal{A}}_{2}+ i a \tilde{\mathcal{A}}_{1})\log \mathcal{G}_{\psi}(f)=0
\end{array}
\end{equation}
since $\mathcal{G}_{\psi}(f)=\mathcal{S}\mathcal{W}_{\psi}(f)$ and
$\mathcal{A}_{i}=\mathcal{S}^{-1}\tilde{\mathcal{A}}_{i}\mathcal{S}$ for $i=1,2,3$.

For the case $a=1$, equation (\ref{CR1}) was noted in \cite{Daudet}. Regarding
general scale $a>0$ we note that
\[
\mathcal{G}_{\psi_{a}}f(p,q)= \sqrt{a}\; \mathcal{G}_{\mathcal{D}_{\frac{1}{a}}\psi}(f)(p,q)=\sqrt{a}\mathcal{G}_{\psi}\mathcal{D}_{\frac{1}{a}}f(\frac{p}{a},a q)
\]
with $\psi=\psi_{a=1}$ with unitary dilation operator $\mathcal{D}_{a}:\mathbb{L}_{2}(\R) \to \mathbb{L}_{2}(\R)$ given by
\begin{equation}\label{dil}
\mathcal{D}_{a}(\psi)(x)=a^{-\frac{1}{2}}f(x/a), \qquad a>0.
 \end{equation}
As a direct consequence of Eq.~\!(\ref{CR1}), respectively (\ref{CR1q}), we have
\begin{equation} \label{fnalCRwitha}
\begin{array}{lll}
|\tilde{U}^{a}|\partial_{q}\tilde{\Omega}^{a}=- a^2 \partial_{p}|\tilde{U}^{a}| &\textrm{ and } &|\tilde{U}^{a}|\partial_{p}\tilde{\Omega}^{a}= a^{-2}\partial_{q}|\tilde{U}^{a}|\; + \;  2\pi q . \\
\mathcal{A}_{2}\Omega^{a}=-a^{2} \frac{\mathcal{A}_{1}|U^{a}|}{|U^{a}|} &\textrm{ and } & \mathcal{A}_{1}\Omega^{a}=a^{-2} \frac{\mathcal{A}_{2}|U^{a}|}{|U^{a}|}\ . \end{array}
\end{equation}
%\mbox{
where $\tilde{U}^{a}= \mathcal{G}_{\psi_{a}}(f)$, $U^{a}\!= \mathcal{W}_{\psi_{a}}(f)\!$, $\tilde{\Omega}^{a}\!=\! \arg \{\mathcal{G}_{\psi_{a}}(f)\}$ and $\Omega^{a}\!=\! \arg \{\mathcal{W}_{\psi_{a}}(f)\}$.

The proper differential-geometric context for the analysis of the evolution equations (and in particular the involved Cauchy-Riemann equations) that we study in this paper is provided by sub-Riemannian geometry.
%In order to get a differential geometrical understanding of our parabolic evolutions on Gabor transforms and the involved Cauchy-Riemann relations we need some additional structure.
As already mentioned
in Subsection \ref{subsect:Horizontal} we omit $\mathcal{A}_{3}$ from the tangent bundle $T(H_{r})$ and consider the sub-Riemannian manifold $(H_{r}, {\rm d}\mathcal{A}^{3})$,  recall (\ref{dualbasis}), as the domain
of the evolved Gabor transforms. Akin to our previous work on parabolic evolutions on orientation scores (defined on the sub-Riemannian manifold $(SE(2),-\sin \theta {\rm d}x+\cos \theta {\rm dy})$) cf.\cite{DuitsAMS2}, we need a left-invariant first fundamental form (i.e. metric tensor) on this sub-Riemannian manifold in order to analyze our parabolic evolutions from a geometric viewpoint.
\begin{lemma}
The only left-invariant metric tensors $\gothic{G}: H_{r} \times T(H_{r}) \times T(H_{r}) \to \mathbb{C}$ on the sub-Riemannian manifold $(H_{r}, {\rm d}\mathcal{A}^{3})$ are given by
\[
\gothic{G}_{g}= \sum \limits_{(i,j) \in \{1,2\}^2}  g_{ij} \left.{\rm d}\mathcal{A}^{i}\right|_{g} \otimes \left.{\rm d}\mathcal{A}^{j}\right|_{g}
\]
\end{lemma}
\textbf{Proof } Let $\gothic{G}: H_{r} \times T(H_{r}) \times T(H_{r}) \to \mathbb{C}$ be a left-invariant metric tensor on the sub-Riemannian manifold $(H_{r}, {\rm d}\mathcal{A}^{3})$. Then since
the tangent space of $(H_{r}, {\rm d}\mathcal{A}^{3})$ at $g \in H_{r}$ is spanned by $\{\left.\mathcal{A}_{1}\right|_{g},\left.\mathcal{A}_{2}\right|_{g}\}$ we have
\[
\gothic{G}_{g}= \sum \limits_{(i,j) \in \{1,2\}^2}  g_{ij}(g) \left.{\rm d}\mathcal{A}^{i}\right|_{g} \otimes \left.{\rm d}\mathcal{A}^{j}\right|_{g}
\]
for some $g_{ij}(g) \in \mathbb{C}$. Now $\gothic{G}$ is left-invariant, meaning
$
\gothic{G}_{gh}((L_{g})_{*}X_h,(L_{g})_{*}Y_h)=\gothic{G}_{h}(X_h,Y_h)$
for all vector fields $X,Y$ on $H_{r}$, and since our basis of left-invariant vector fields satisfies
$(L_{g})_{*}\left.\mathcal{A}_{i}\right|_{h}=\left.\mathcal{A}_{i}\right|_{gh}$ and $\langle {\rm d}\mathcal{A}^{i}, \mathcal{A}_{j}\rangle=\delta^{i}_{j}$
we deduce $g_{ij}(gh)=g_{ij}(h)=g_{ij}(e)$ for all $g,h \in H_{r}$.$\hfill \Box$ \\
\\
Akin to the related quadratic forms in the generator of our parabolic evolutions we restrict ourselves to the case
where the metric tensor is diagonal
\begin{equation}\label{metric}
\gothic{G}_{\beta}= \sum \limits_{(i,j) \in \{1,2\}^2} g_{ij} {\rm d}\mathcal{A}^{i} \otimes {\rm d}\mathcal{A}^{j} =\beta^{4}{\rm d}\mathcal{A}^{1} \otimes {\rm d}\mathcal{A}^{1} +
{\rm d}\mathcal{A}^{2} \otimes {\rm d}\mathcal{A}^{2}.
\end{equation}
Here the fundamental positive parameter $\beta^{-1}$ has physical dimension length, so that this first fundamental form is consistent
with respect to physical dimensions. Intuitively, the parameter $\beta$ sets a global balance between changes in frequency space and changes in position space within the induced metric
\[
d(g,h)= \inf \limits_{{\scriptsize \begin{array}{c}
\gamma \in C([0,1],H_{3}), \\
\gamma(0)=g, \gamma(1)=h,
\\
\langle \left.{\rm d}\mathcal{A}^{3}\right|_{\gamma},\dot{\gamma}\rangle=0 \end{array}
}} \int_{0}^{1} \sqrt{\left.\gothic{G}_{\beta}\right|_{\gamma(t)}(\dot{\gamma}(t),\dot{\gamma}(t))}\,{\rm d}t.
\]
Note that the metric tensor $\gothic{G}_{\beta}$ is bijectively related to the linear operator $G:\gothic{H} \to \gothic{H}'$, where $\gothic{H}=\textrm{span}\{\mathcal{A}_{1},\mathcal{A}_{2}\}$ denotes the horizontal part of the tangent space, that maps
$\mathcal{A}_{1}$ to $\beta^{4}{\rm d}\mathcal{A}^{1}$ and $\mathcal{A}_{2}$ to ${\rm d}\mathcal{A}^{2}$.
The inverse operator of $G$ is bijectively related to
\[
\gothic{G}^{-1}_{\beta}= \sum \limits_{(i,j) \in \{1,2\}^2} g^{ij} \mathcal{A}_{i} \otimes \mathcal{A}_{j} = \beta^{-4}\mathcal{A}_{1} \otimes \mathcal{A}_{1} +
\mathcal{A}_{2} \otimes \mathcal{A}_{2}\ .
\]
The fundamental parameter $\beta$ is inevitable when dealing with flows in the Gabor domain, since position $p$ and frequency $q$ have different physical dimension. Later, in Section \ref{ch:covdiff}, we will need this metric in the design of adaptive diffusion on Gabor transforms.
In this section we primarily use it for geometric understanding of the Cauchy-Riemann relations on Gabor transforms, which we will employ in our convection schemes in the subsequent section.

\begin{figure}[h]
\centerline{
\includegraphics[width=0.45\hsize]{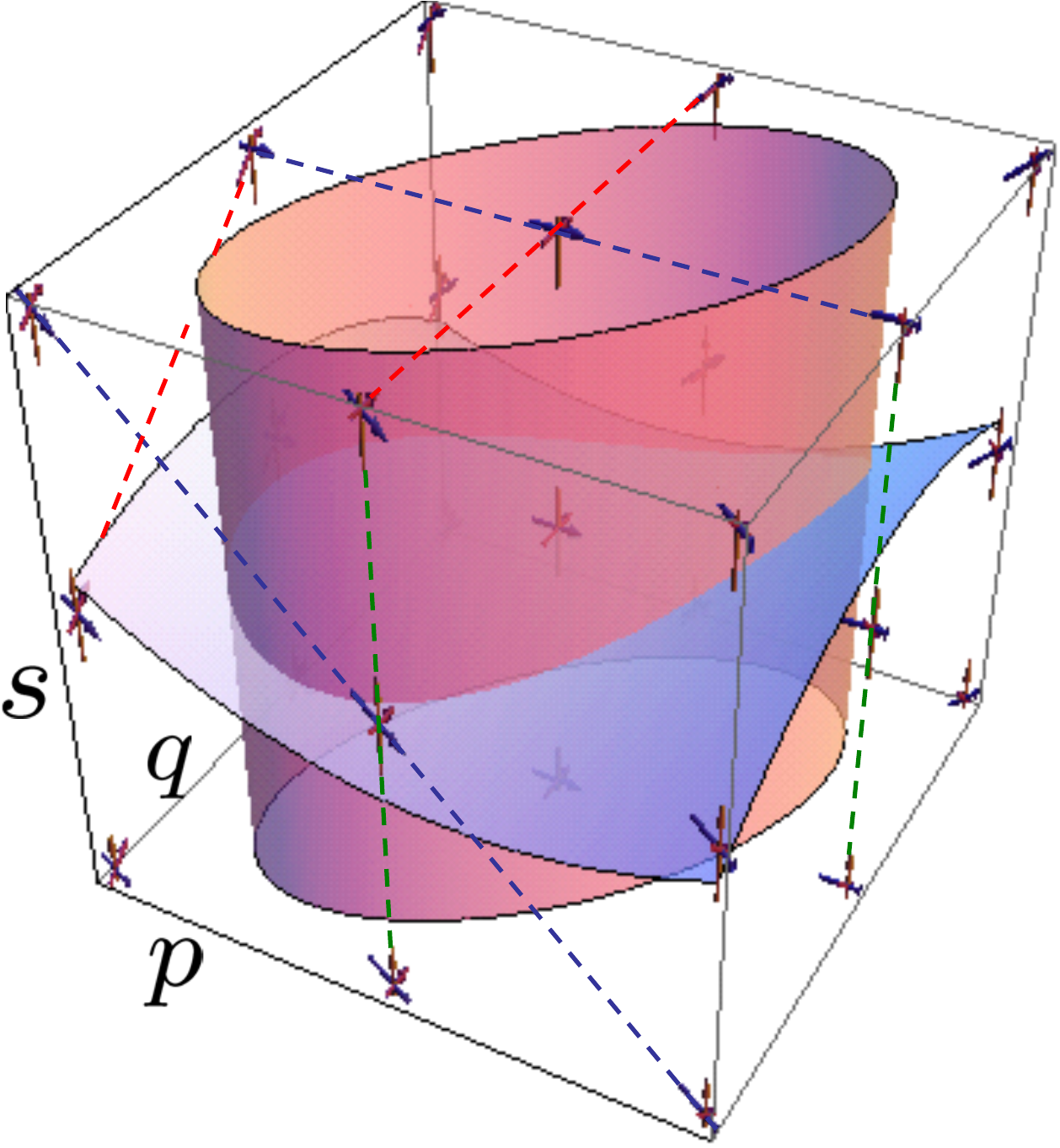}
}
\caption{Equi-amplitude plane (red) and equi-phase plane in the Gabor transform of a chirp signal, that we shall compute exactly in Section \ref{ch:chirp}. The left-invariant horizontal gradients of these surfaces are locally $\gothic{G}_{\beta}$-orthogonal to each other, cf.~Eq.~(\ref{CRMetric}), with
$\gothic{G}_{\beta}= \beta^{4}{\rm d}\mathcal{A}^{1} \otimes {\rm d}\mathcal{A}^{1} +  {\rm d}\mathcal{A}^{2} \otimes {\rm d}\mathcal{A}^{2}$ with $\beta=a^{-1}$. The left-invariant vector fields $\{\mathcal{A}_{1},\mathcal{A}_{2},\mathcal{A}_{3}\}=\{\partial_{p}+\frac{q}{2}\partial s, \partial_{q}-\frac{p}{2}\partial s, \partial_{s}\}$ form a local frame of reference which is indicated by the arrows.
Some exponential curves (the auto-parallel curves w.r.t. Cartan connection)
are indicated by dashed lines.
 }\label{fig:CR}
\end{figure}

In potential theory and fluid dynamics the Cauchy-Riemann equations for complex-valued analytic functions, impose orthogonality between flowlines and equipotential lines. A similar geometrical interpretation can be deduced from the Cauchy-Riemann relations (\ref{fnalCRwitha}) on Gabor transforms :
%that hold between local phase and local amplitude can now be written in geometrical form:
\begin{lemma} \label{lemma:CR}
Let $U:=\mathcal{W}_{\psi_{a}}f=|U| e^{i \Omega}$ be the Gabor transform of a signal $f \in \mathbb{L}_{2}(\R)$. Then
\begin{equation} \label{CRMetric}
\gothic{G}^{-1}_{\beta=\frac{1}{a}}({\rm d} \log |U|, \mathbb{P}_{\gothic{H}^{*}}{\rm d}\Omega)=
\gothic{G}^{-1}_{\beta=\frac{1}{a}}({\rm d}|U|, \mathbb{P}_{\gothic{H}^{*}}{\rm d}\Omega)
=0,
\end{equation}
where the left-invariant gradient of modulus and phase equal
\[
{\rm d}\Omega= \sum \limits_{i=1}^{3} \mathcal{A}_{i}\Omega \; {\rm d}\mathcal{A}^{i}, \ \ {\rm d}|U|= \sum \limits_{i=1}^{2} \mathcal{A}_{i}|U| \; {\rm d}\mathcal{A}^{i},
 \]
whose horizontal part equals $\mathbb{P}_{\gothic{H}^{*}}{\rm d}\Omega=\sum \limits_{i=1}^{2} \mathcal{A}_{i}\Omega \; {\rm d}\mathcal{A}^{i}$, $\mathbb{P}_{\gothic{H}^{*}}{\rm d}|U|\equiv{\rm d}|U|$.
\end{lemma}
\textbf{Proof }By the second line in (\ref{fnalCRwitha}) we have
\[
a^{-2} \gothic{G}_{\beta=a^{-1}}^{-1}({\rm d}\log|U|, \mathbb{P}_{\gothic{H}^{*}} {\rm d}\Omega)= a^{2}
\frac{\mathcal{A}_{1}|U| \mathcal{A}_{1}\Omega}{|U|}  +  a^{-2}\frac{\mathcal{A}_{2}|U| \mathcal{A}_{2}\Omega}{|U|}=0,
\]
from which the result follows.$\hfill \Box$
\begin{corollary}
Let $g_0 \in H_{r}$. Let $\mathcal{W}_{\psi_{a}}f=|U| e^{i \Omega}$ with in particular
$\mathcal{W}_{\psi_{a}}f(g_0)=|U(g_0)| e^{i \Omega(g_0)}$. Then the horizontal part $\mathbb{P}_{\gothic{H}^{*}}\left. {\rm d}\Omega \right|_{g_0}$ of the normal covector $\left. {\rm d}\Omega \right|_{g_0}$ to the equi-phase surface $\{(p,q,s) \in H_{r}\; |\; \Omega(p,q,s)=\Omega(g_0)\}$  is $\gothic{G}_{\beta}$-orthogonal to the normal covector $\left. {\rm d}|U| \right|_{g_0}$
to the equi-amplitude surface $\{(p,q,s) \in H_{r}\; |\; |U|(p,q,s)=|U|(g_{0})\}$.
\end{corollary}
%Equality (\ref{CRMetric}) sets an important analogy with potential theory/fluid dynamics where equi-potential lines of the harmonic potential are orthogonal to flow lines (isolines of the Harmonic conjugate) due to the Cauchy Riemann relations. In our technical report
%\cite[App.F]{DuitsFuehrJanssen}
%we use this to show that the convection-algorithms of Section \ref{ch:reass} take place along horizontal curves whose tangent vectors are $\mathcal{G}_{\beta}$-orthogonal to the tangent vectors of horizontal curves aligned with the principal axis of our horizontal diffusion-algorithms (again fitting in our general framework (\ref{theeqs}) but then with $a_{1}=a_{2}=0$ in stead of $D=0$ in  (\ref{therightchoice})) explained in \cite[ch:7]{DuitsFuehrJanssen}.
For a visualization of the Cauchy-Riemannian geometry see Fig.\ref{fig:CR}, where we also include the exponential curves along which our diffusion and convection (Eq.(\ref{theeqs})) take place.
\begin{remark}
Akin to our framework of left-invariant evolutions on orientation scores \cite{DuitsAMS2} one express the left-invariant evolutions on Gabor transforms (Eq.(\ref{theeqs})) in covariant derivatives, so that transport and diffusion takes place along
the covariantly constant curves (auto-parallels) w.r.t. Cartan Connection on the sub-Riemannian manifold $(H_{r},{\rm d}\mathcal{A}^{3})$.
A brief computation (for analogous details see \cite{DuitsAMS2})
shows that the auto-parallels $t \mapsto \gamma(t)$ w.r.t. the Cartan connection $\nabla$ coincide with the horizontal exponential curves. Auto-parallels are by definition curves that satisfy
\begin{equation} \label{basic}
\nabla_{\dot{\gamma}}\dot{\gamma}=0 \desda \ddot{\gamma}^{i} - \sum_{k,j} c^{i}_{kj} \dot{\gamma}^{k}\dot{\gamma}^{j}=\ddot{\gamma}^{i}=0,
\end{equation}
where in case of the Cartan-connection the Christoffel-symbols coincide with minus the anti-symmetric Lie algebra structure constants
and with \mbox{$\dot{\gamma}^{i}=\langle \left.{{\rm d}\mathcal{A}^{i}}\right|_{\gamma},\dot{\gamma} \rangle$}. So indeed Eq.~\!(\ref{basic}) holds iff
$\dot{\gamma}^{i}= c^{i} \in \R$, $i=1,2$, i.e. {\small $\gamma(t)= \gamma(0)\, \exp(t \sum \limits_{i=1}^{2} c^{i} A_{i})$} for all $t \in \R$.
\end{remark}

%hier

\section{Convection operators on Gabor Transforms that are both phase-covariant and phase-invariant \label{ch:reass}}

In differential reassignment, cf. \cite{Daudet,Chassande} the practical goal
is to sharpen Gabor distributions towards lines (minimal energy curves \cite[App.D]{DuitsFuehrJanssen}) in $H_{r}$, while maintaining the signal as much as possible.

We would like to achieve this by left-invariant convection on Gabor transforms $U:=\mathcal{W}_{\psi}(f)$. This means one should set $D=0$ in Eq.~\!(\ref{theeqs}) and (\ref{theeqsPS}) while considering the mapping $U \mapsto W(\cdot,t)$ for a suitably chosen fixed time $t>0$. Let us denote this mapping by $\Phi_{t}: \mathcal{H}_{n} \to \mathcal{H}_{n}$ given by $\Phi_{t}(U)=W(\cdot,t)$.
Such a mapping is called \emph{phase invariant} if
\[
\arg \left( (\Phi_{t}(U))(p,q,s) \right)= \arg \left( U(p,q,s) \right),
\]
for all $(p,q,s) \in H_{r}$ and all $U \in \mathcal{H}_{n}$, allowing us to write $\Phi_{t}(|U|e^{i\Omega})=e^{i\Omega} \cdot \Phi_{t}^{\textrm{net}}(|U|)$ where $\Phi_{t}^{\textrm{net}}$ is the effective operator on the modulus.
Such a mapping is called \emph{phase covariant} if the phase moves along with the flow (characteristic curves of transport), i.e. if Eq.~(\ref{phasecovariant}) is satisfied. Somewhat contrary to intuition, the two properties are not exclusive.

Our convection operators $\Phi_{t}$ (obtained by setting $D=0$ in Eq.~\!(\ref{theeqs})) are both phase covariant and phase invariant iff their generator is. In order to achieve both phase invariance and phase covariance one should construct the generator such that the flow is along equi-phase planes of the initial Gabor transform $\mathcal{W}_{\psi}f$. As we restricted ourselves to horizontal convection there is only one direction in the horizontal part of the tangent space we can use.
\begin{lemma}\label{lemma:0}
let $\Omega_{g_0}$ be an open set in $H_{r}$ and
let $U : \Omega_{g_0} \to \mathbb{C}$ be differentiable.
The only horizontal direction in the tangent bundle above $\Omega_{g_0}$ that preserves the phase of $U$ is given by
\[
\textrm{span}\{ -\mathcal{A}_{2}\Omega \mathcal{A}_{1} + \mathcal{A}_{1}\Omega \mathcal{A}_{2}\},
\]
with $\Omega=\arg \{U\}$.
\end{lemma}
\textbf{Proof }The horizontal part of the tangent space is spanned by $\{\mathcal{A}_{1},\mathcal{A}_{2}\}$, the horizontal part of the phase gradient is given by $\mathbb{P}_{\gothic{H}^{*}} {\rm d}\Omega=\mathcal{A}_{1}{\Omega} {\rm d}\mathcal{A}^{1}+\mathcal{A}_{2}{\Omega} {\rm d}\mathcal{A}^{2}$. Solving for
\[
\langle \mathbb{P}_{\gothic{H}^{*}} {\rm d}\Omega , \alpha^{1} \mathcal{A}_{1} + \alpha^{2} \mathcal{A}_{2} \rangle =0
\]
yields $(\alpha^{1},\alpha^{2})= \lambda (-\mathcal{A}_{2}\Omega,\mathcal{A}_{1}\Omega)$, $\lambda \in \mathbb{C}$. $\hfill \Box$

As a result it is natural to consider the following class of convection operators.
\begin{lemma}
The horizontal, left-invariant, convection generators $\mathcal{C}:\mathcal{H}_{n} \to \mathcal{H}_{n}$  given by
\[
\mathcal{C}(U) = \mathcal{M}(|U|)(-\mathcal{A}_{2}\Omega \mathcal{A}_{1}U + \mathcal{A}_{1}\Omega \mathcal{A}_{2}U), \qquad \textrm{ where }\Omega=\arg \{U\},
\]
and where $\mathcal{M}(|U|)$ a multiplication operator
naturally associated to a bounded monotonically increasing differentiable function $\mu:[0,\max(|U|)] \to [0,\mu(\max(|U|))] \subset \R$
with $\mu(0)=0$, i.e. $(\mathcal{M}(|U|)V)(p,q)= \mu(|U|(p,q))\; V(p,q)$ for all $V \in \mathcal{H}_{n}, (p,q)\in \R^2$,
are well-defined and both phase covariant and phase invariant.
\end{lemma}
\textbf{Proof }Phase covariance and phase invariance follows by Lemma \ref{lemma:0}.
The operators are well-posed as the absolute value of Gabor transform is almost everywhere smooth (if $\psi$ is a Schwarz function) bounded and moreover $\mathcal{C}$ can be considered as
an unbounded operator from $\mathcal{H}_{n}$ into $\mathcal{H}_{n}$, as the bi-invariant space $\mathcal{H}_{n}$ is invariant under bounded multiplication operators that do not depend on the phase $z=e^{2\pi i s}$. $\hfill \Box$ \\
\\
For Gaussian kernels $\psi_{a}(\xi)=e^{-a^{-2}\xi^{2} n \pi}$ we may apply the Cauchy Riemann relations (\ref{CR1}) which simplifies for the special case $\mathcal{M}(|U|)=|U|$ to
\begin{equation} \label{phaseinvariantconvgen}
\mathcal{C}(e^{i \Omega} |U|)= ( a^{2} (\partial_{p}|U|)^2 + a^{-2}(\partial_{q}|U|)^2 )\;e^{i \Omega} .
\end{equation}
Now consider the following phase-invariant adaptive convection equation on $H_{r}$,
\begin{equation} \label{conpde}
\boxed{
\left\{
\begin{array}{l}
\partial_{t}W(g,t)=-\mathcal{C}(W(\cdot,t))(g), \\
W(g,0)=U(g)
\end{array}
\right.
}
\end{equation}
with either
\begin{equation} \label{conchoice}
\boxed{
\begin{array}{l}
1. \ \ \mathcal{C}(W(\cdot,t))= \mathcal{M}(|U|)\, (-\mathcal{A}_{2} \Omega, \mathcal{A}_{1} \Omega) \cdot (\mathcal{A}_{1}W(\cdot,t), \mathcal{A}_{2} W(\cdot,t)) \textrm{ or } \\[8pt]
2. \ \ \mathcal{C}(W(\cdot,t))= e^{i \Omega} \left(a^{2} \frac{(\partial_{p} |W(\cdot,t)|)^2}{|W(\cdot,t)|} + a^{-2} \frac{(\partial_{q}|W(\cdot,t)|)^2}{|W(\cdot,t)|}\right).
\end{array}
}
\end{equation}
In the first choice we stress that $\arg (W(\cdot,t))=\arg(W(\cdot,0))=\Omega$, since transport only takes place along iso-phase surfaces. Initially, in case $\mathcal{M}(|U|)=1$ the two approaches are the same since at $t=0$ the Cauchy Riemann relations (\ref{fnalCRwitha}) hold, but as time increases the Cauchy-Riemann equations are violated (this directly follows by the preservation of phase and non-preservation of amplitude).
Consequently, generalizing the single step convection schemes in \cite{Daudet,Chassande} to a continuous time axis produces two options:
 %a problem that also arises in the single step convection schemes in \cite{Daudet,Chassande}.
%Nevertheless, both approaches (\ref{conchoice}) in (\ref{conpde}) make sense from a practical point of view. For small convection times %$t>0$ the difference turns out to be hardly visible, but after a while we have observed differences.

%Here we stress that by taking the modulus, we omit the $s$-dependence and we automatically apply our algorithm in phase space.
\begin{enumerate}
\item With respect to the first choice in (\ref{conchoice}) in (\ref{conpde}) (which is much more cumbersome to implement) we follow the authors in \cite{Daudet} and consider the equivalent equation on phase space:
\begin{equation} \label{reassigmentphaseinvphasespace}
%\boxed{
\left\{
\begin{array}{l}
\partial_{t}\tilde{W}(p,q,t)=-\tilde{\mathcal{C}}(\tilde{W}(\cdot,t))(p,q), \\
\tilde{W}(p,q,0)=\mathcal{G}_{\psi}f(p,q)=:\tilde{U}(p,q)=e^{i \tilde{\Omega}(p,q)} |\tilde{U}(p,q)|= e^{i \tilde{\Omega}(p,q)} |U|(p,q)
\end{array}
\right.
%}
\end{equation}
with $
\tilde{C}(\tilde{W}(\cdot,t))= \mathcal{M}(|U|)\left(-\tilde{\mathcal{A}}_{2} \tilde{\Omega} \tilde{\mathcal{A}}_{1}\tilde{W}(\cdot,t) + (\partial_{q}\tilde{\Omega}-2\pi q) \tilde{\mathcal{A}}_{2}\tilde{W}(\cdot,t)\right)$,
where we recall $\mathcal{G}_{\psi}=\mathcal{S}\mathcal{W}_{\psi}$ and
$\mathcal{A}_{i}=\mathcal{S}^{-1}\tilde{\mathcal{A}}_{i}\mathcal{S}$ for $i=1,2,3$. Note that the authors in \cite{Daudet} consider the case $\mathcal{M}=1$. %However the case $\mathcal{M}=1$ and the earlier %mentioned case $\mathcal{M}(|U|)=|U|$ are equivalent :
%\[
%\frac{\partial}{\partial t}|U|= a^{2} \frac{(\partial_{p}|U|)^2}{|U|}+ a^{-2} \frac{(\partial_{q}|U|)^2}{|U|} \desda
%\frac{\partial }{\partial t} \log |U| = a^{2}(\partial_{p} \log |U|)^2+a^{-2}(\partial_{q} \log |U|)^2.
%\]
%%Although the approach in \cite{Daudet} is highly plausible, the authors did not
%%
In addition to \cite{Daudet} we provide in Section \ref{ch:upwind} an explicit computational finite difference scheme acting on discrete subgroups of $H_{r}$, which is non-trivial due to the oscillations in the Gabor domain. %and which is equivalent to our numerical schemes for erosions on diffusion weighted MRI %\cite{Fran2009,CreusenSSVM2011}.
\item
The second choice in Eq.~(\ref{conchoice}) within Eq.~(\ref{conpde}) is just a phase-invariant inverse Hamilton Jacobi equation on $H_{r}$, with a Gabor transform as initial solution. Rather than computing the viscosity solution cf.~\cite{Evans} of this non-linear PDE, we may
as well store the phase and apply an inverse Hamilton Jacobi system on $\R^2$ with the amplitude $|U|$ as initial condition and
multiply with the stored phase factor afterwards.
More precisely, the viscosity solution of the Hamilton Jacobi system on the modulus
is
given by a basic inverse convolution over the
$(\textrm{max},+)$ algebra, \cite{Burgeth}, (also known as erosion operator in image analysis)
\begin{equation} \label{erosion}
\tilde{W}(p,q,t)=(\tilde{\Phi}_{t}(U))(p,q,t) =(K_{t} \ominus |U|)(p,q) e^{i \Omega(p,q,t)}\ ,
\end{equation}
with kernel
\begin{equation} \label{erosionkernel}
K_{t}(p,q)= \frac{a^{-2}p^{2}+a^{2} q^2}{4t}
\end{equation}
%(describing the
%growth of balls in $\R^2$),
where
\[
(f \ominus g)(p,q)= \inf \limits_{(p',q') \in \R^2} \left[g(p',q')+f(p-p',q-q')\right].
\]
Here the homomorphism between erosion and inverse diffusion
is given by the Cramer transform $C= \gothic{F} \circ \log \circ \mathcal{L}$,
\cite{Burgeth}, \cite{Akian}, that is a concatenation of the multi-variate Laplace transform, logarithm and Fenchel transform %and whose
%The Fenchel
%transform maps a convex function $c:\R^{2} \to \overline{\R}$ onto
%\begin{equation}\label{Fenchel}
%\ul{x} \mapsto [\gothic{F}c](\ul{x})=\sup \{\ul{y}\cdot \ul{x}-c(\ul{y})\; |\; \ul{y} \in \R^2 \}.
%\end{equation}
%isomorphic property is
so that
\[
\begin{array}{l}
\mathcal{C}(f*g) =\gothic{F} \log \mathcal{L}(f*g)=\gothic{F}(\log \mathcal{L}f + \log \mathcal{L}g) = \mathcal{C}f \oplus \mathcal{C}g,
\end{array}
\]
with convolution on the $(\textrm{max},+)$-algebra given by
$f \oplus g(\ul{x}) =\sup \limits_{\ul{y} \in \R^d} [f(\ul{x}-\ul{y}) + g(\ul{y})]$.
\end{enumerate}

%\section{Implementation}
\section{Discrete Gabor Transforms \label{ch:discrete}}

In order to derive various suitable algorithms
for differential reassignment and diffusion we consider discrete Gabor transforms.
We show that Gabor transforms are defined on a (finite) group quotient within the discrete Heisenberg group. In the subsequent section
we shall construct left-invariant shift operators for the generators in our left-invariant evolutions on this quotient group.

Let the discrete signal is given by $\ul{f}=\{\ul{f}[n]\}_{n=0}^{N-1}:=\{f\left( \frac{n}{N}\right)\}_{n=0}^{N-1} \in \R^{N}$. Let the discrete kernel be given by a sampled Gaussian kernel
\begin{equation}\label{Gaussian}
\PPP=\{\PPP[n]\}_{n=-(N-1)}^{N-1}:=\{e^{- \frac{(|n| -\lfloor \frac{N-1}{2} \rfloor)^{2} \pi }{N^2a^2}}\}_{n=-(N-1)}^{N-1} \in \mathbb{R}^{N},
\end{equation}
with $\frac{\pi}{a^2}=\frac{1}{2\sigma^2}$ where $\sigma^{2}$ is the variance of the Gaussian.
The discrete Gabor transform of $\ul{f}$ is then given by
\begin{equation} \label{discretegabor}
\begin{array}{ll}
(\mathcal{W}_{\PPP}^{D}\ul{f})[l,m,k]&:=e^{-2\pi i (\frac{k}{Q}-\frac{ml L}{2M})} \frac{1}{N}\sum \limits_{n=0}^{N-1}\overline{\PPP[n-l L]}\ul{f}[n]\, e^{-\frac{2\pi i nm}{M}}
\end{array}
\end{equation}
with $L,K,N,M,Q \in \mathbb{N}$ and
\begin{equation} \label{LNK}
k=0,1,\ldots, Q\!-\!1 \textrm{ and }l=0,\ldots, K\!-\!1, \ m=0,\ldots, M\!-\!1, \textrm{ with } L=\frac{N}{K},
\end{equation}
and integer oversampling $P=M/L \in \mathbb{Z}$.
Note that we follow the notational conventions of the review paper \cite{Hlawatsch}.
One has
\begin{equation}\label{relation}
\frac{1}{M}=\frac{K}{P} \frac{1}{N}.
\end{equation}
It is important that the discrete kernel is $N$ periodic since $N=KL$ implies
\[
\forall_{\ul{f} \in \ell_{2}(I)} \forall_{l,m,k} \mathcal{W}_{\PPP}^{D}\ul{f}[l+K,m,k]=\mathcal{W}_{\PPP}^{D}\ul{f}[l,m,k] \desda \forall_{n=0,\ldots N} \PPP[n-N]=\PPP[n],
\]
where $I=\{0,\ldots, N-1\}$. Moreover, we note that the kernel chosen in (\ref{Gaussian}) is even.

For Riemann-integrable $f$ with support within $[0,1]$ and $\psi$ even with support within $[-1,1]$, say
\begin{equation} \label{Gaussiancont}
\psi(\xi)= e^{-\frac{\pi ||\xi|-\frac{1}{2}|^{2}}{a^2}} 1_{[-1,1]}(\xi),
\end{equation}
we have
\[
\begin{array}{ll}
(W_{\PPP}^{D}\ul{f})[l,m,k]
 &=\frac{1}{N}e^{-2\pi i (\frac{k}{Q}-\frac{ml L}{2M})} \sum \limits_{n=0}^{N-1} e^{-\pi a^{-2}     \frac{(|n-lL|-\lfloor \frac{N-1}{2}\rfloor)^2}{N^2}   }
f\left(\frac{n }{N} \right) e^{-\frac{2\pi i nm }{M}} \\
 &=e^{-2\pi i (\frac{k}{Q}-\frac{ml}{2P})}\, \frac{1}{N} \sum \limits_{n=0}^{N-1} e^{-\pi a^{-2}\left(|\frac{n}{N}-\frac{l}{K}|-\frac{1}{N}\lfloor \frac{N-1}{2}\rfloor\right)^2}
f\left(\frac{n }{N} \right) e^{-\frac{2\pi (K/P) i nm }{N}}\\
 & \to  e^{-2\pi i (\frac{k}{Q}-\frac{1}{2} \frac{mK}{P} \frac{l}{K})}
 \int \limits_{0}^{1} f(\xi) e^{-\frac{\pi \left(\left|\xi-\frac{l}{K}\right|-\frac{1}{2}\right)^2}{a^2}} e^{-2\pi i \xi \left(\frac{mK}{P}\right)}\, {\rm d}\xi.
\end{array}
\]
Consequently, we obtain the pointwise limit (in the reproducing kernel space of Gabor transforms)
\begin{equation} \label{convresult}
(W_{\PPP}^{D}\ul{f})[l,m,k] \to \mathcal{W}^{n=1}_{\psi}f(p=\frac{l}{K},q=\frac{m K}{P},s=\frac{k}{Q}) \textrm{ as }N \to \infty,
\end{equation}
where we keep both $P$ and $K$ fixed so that only $M \to \infty$ as $N \to \infty$ and with with scaled Gaussian kernel $\psi(\xi)=e^{-\pi a^{-2} \xi^2} 1_{[-1,1]}(\xi)$.
To this end we recall that the continuous Gabor transform was given by
\[
\mathcal{W}^{n=1}_{\psi}f(p,q,s)=e^{-2\pi i (s-\frac{pq}{2})} \int_{\R} \overline{\psi(\xi-p)} f(\xi) e^{-2\pi i \xi q}\, {\rm d}\xi.
\]

\subsection{Diagonalization of the Gabor transform \label{ch:zak}}

In our algorithms, we follow \cite{GuidoJanssen} and \cite{Hlawatsch} and use the diagonalization of the discrete Gabor transform by means of the discrete Zak-transform.
%In fact the discrete Zak-transform is given by
%\section{Diagonalization of the discrete Gabor transform by discrete Zak-transform}
The finite frame operator $\gothic{F}:\ell_{2}(I) \to \ell_{2}(I)$ equals
\[
[\gothic{F}\ul{f}][n] = \sum \limits_{l=0}^{K-1} \sum \limits_{m=0}^{M-1} (\PPP_{lm},\ul{f})\PPP_{lm}[n], \qquad n \in I=\{0,\ldots,N-1\},
\]
with $\PPP_{lm}= \mathcal{U}_{[l,m,k=-\frac{Qlm}{2P}]}\PPP$.
%By its Riesz-bounds the self-adjoint operator $\gothic{F}$ is bounded and invertible.
Operator $\gothic{F}=\gothic{F}^*$ is coercive and has the following orthonormal eigenvectors:
\[
\ul{u}_{nk}[n']=\frac{1}{\sqrt{K}}v[n'-n]\, e^{\frac{2\pi i k}{N}(n'-n)}, \textrm{ with }v(n)=\sum \limits_{l=-\infty}^{\infty} \delta[n-lL]\
\]
\textrm{ for } $n \in \{0,\ldots,L-1\}, k \in \{0,\ldots,K-1\}$, where $N=KL$,
and it is diagonalized by:
\[
\gothic{F}= (Z^{D})^{-1} \circ \Lambda \circ Z^{D}\ ,
\]
with Discrete Zak transform given by $[Z^{D}\ul{f}][n,k]=(\ul{u}_{nk},\ul{f})_{\ell_{2}(I)}$, i.e.
$\gothic{F}\ul{f}=\sum \limits_{n=0}^{L-1} \sum \limits_{k=0}^{K-1} \lambda_{nk}(\ul{u}_{nk},\ul{f})\ul{u}_{nk}$, with eigenvalues $\lambda_{nk}= L \sum \limits_{p=0}^{P-1} |Z\PPP^{D}[n,k-p \frac{N}{M}]|^2$ and integer oversampling factor $P=M/L$.
%\begin{remark}
%It is well-known that (continuous) Gabor transform relates to inverse Fourier transform on the reduced Heisenberg group $H_{r}$, see Appendix \ref{ch:FT}. Similarly, discrete Gabor transform relates to inverse Fourier transform on $\gothic{h}_{r}$. The considerations in this section therefore boil down to diagonalization of inverse Fourier transform on the \emph{discrete} Heisenberg group. In the continuous setting (or in the infinite discrete setting) one could apply a similar diagonalization of the frame operator, but in this setting
%the self-adjoint frame-operator is not compact and is bound to work with generalized eigenvectors in a Gelfand-triple.
%\end{remark}

\section{Discrete Left-invariant vector fields \label{ch:discreteVF}}

Let the discrete signal is given by $\ul{f}=\{\ul{f}[n]\}_{n=0}^{N-1}:=\{f\left( \frac{n}{N}\right)\}_{n=0}^{N-1} \in \R^{N}$. Then
similar to the continuous case the discrete Gabor transform of $\ul{f}$ can be written
\begin{equation}\label{Gaborfinite}
[\mathcal{W}_{\PPP}^{D}\ul{f}][l,m,k]= (\mathcal{U}_{[l,m,k]}\PPP,\ul{f})_{\ell_{2}(I)},
\end{equation}
where $I=\{0,\ldots, N-1\}$ and $(\ul{a},\ul{b})= \frac{1}{N} \sum \limits_{i=0}^{N-1} \overline{a}_{i} b_{i}$ and where
\begin{equation} \label{groupact}
\mathcal{U}_{[l,m,k]}\PPP[n]= e^{2\pi i (\frac{k}{Q}-\frac{ml }{2 P})} e^{\frac{2\pi i nm}{M}} \PPP[n-lL],
\end{equation}
$l=0,\ldots,K-1, m=0,\ldots, M-1, k=0,\ldots, Q-1$. Next we will show that (under minor additional conditions) Eq.~(\ref{groupact}) gives rise to a group representation of a finite dimensional Heisenberg group $\gothic{h}_{r}$, obtained by taking the quotient of the discrete Heisenberg group $h_r$ with a normal subgroup.
\begin{definition}
Assume $\frac{Q}{2P} \in \mathbb{N}$ then the group $h_{r}$ is the set
$\mathbb{Z}^{3}/(\{0\}\times \{0\} \times Q \mathbb{Z})$ endowed with group product
\begin{equation} \label{discreteprod}
[l,m,k][l',m',k']=[l+l' , m+m' , k+k' +\frac{Q}{2P}(ml'-m'l) \, \textrm{\emph{Mod}} Q].
\end{equation}
\end{definition}
\begin{lemma}
Assume $\frac{Q}{2P} \in \mathbb{N}$ and $\frac{K}{P} \in \mathbb{N}$ and
$L$ even, $N$ even. Then the subgroup
\[
[K\mathbb{Z}, M \mathbb{Z}, Q \mathbb{Z}]:= \{[lK,mM,kQ]\; |\; l,m,k \in \mathbb{Z}\}
\]
is a normal subgroup of $h_{r}$. Thereby, the quotient $\gothic{h}_{r}:= h_{r}/([K\mathbb{Z}, M \mathbb{Z}, Q \mathbb{Z}])$ is a group with product
\begin{equation} \label{discreteprod}
[l,m,k][l',m',k']=[l+l' \; \textrm{Mod}K, m+m' \; \textrm{Mod}M, k+k' +\frac{Q}{2P}(ml'-m'l) \; \textrm{Mod}(Q)].
\end{equation}
Eq.~(\ref{groupact}) defines a group representation on this group
and thereby $\gothic{h}_{r}$ is the domain of discrete Gabor transforms endowed with the group product given by Eq.~(\ref{discreteprod}).
\end{lemma}
\textbf{Proof }
Direct computation yields
\[
\begin{array}{l} \
[l,m,k][l'K, m'M,k'Q]=[l+l'K, m+m'M, k+k'Q+\frac{Q}{2P}(ml'K -m'Ml) +\textrm{Mod} Q]= \\ \
[l'K, m'M,k'Q][l,m,k]=[l+l'K, m+m'M, k+k'Q-\frac{Q}{2P}(ml'K -m'Ml) +\textrm{Mod} Q]
\end{array}
\]
since $M/P=L \in \mathbb{N}$ and we assumed $K/P \in \mathbb{N}$. Consequently,
$H:=[K\mathbb{Z}, M \mathbb{Z}, Q \mathbb{Z}]$ is a normal subgroup and thereby (since
$g_{1}H g_{2}H]=g_{1}g_{2} H$) the quotient
$\gothic{h}_{r}:= h_{r}/([K\mathbb{Z}, M \mathbb{Z}, Q \mathbb{Z}])$ is a group
with well-defined group product (\ref{discreteprod}).
The remainder now follows by direct verification of
\[
\mathcal{U}_{[l,m,k]}\mathcal{U}_{[l',m',k']}= \mathcal{U}_{[l,m,k][l',m',k']}
\]
and $\mathcal{U}_{[l,m,k]}=\mathcal{U}_{[l+K,m+M,k+Q]}$. $\hfill \Box$ \\
%\begin{lemma}
\\
In view of Eq.~(\ref{convresult}) we define a monomorphism between $h_{r}$
and $H_{r}$ as follows.
\begin{lemma}
Define the mapping $\phi: h_r \to H_{r}$ by
\[
\phi[l,m,k]=\left(\frac{l}{K}, \frac{mK}{P},\frac{k}{Q}\right)
\]
which sets a monomorphism between $h_r$%the discrete group $h_r=\{[l,m,k]\; |\; l,m,k \in \mathbb{Z}\}$
%%which is equipped with group product
%%\[
%%[l,m,k][l',m',k']=[l+l' , m+m' , k+k' + \frac{Q}{2P}(ml'-m'l) ].
%%\]
and the continuous Heisenberg group $H_r$.
\end{lemma}
\textbf{Proof } Straightforward computation yields
\[
\begin{array}{ll}
\phi[l,m,k]\phi[l',m',k'] &=\left(\frac{l}{K}, \frac{mK}{P},\frac{k}{Q}\right)\left(\frac{l}{K}, \frac{mK}{P},\frac{k}{Q}\right)\\
 &=\left(\frac{l+l'}{K}, \frac{(m+m')K}{P}, \frac{k+k' + \frac{Q}{2P}(ml'-lm')}{Q} \right)=\phi[[l,m,k][l',m',k']].
\end{array}
\]
from which the result follows. $\hfill \Box$ \\
% sets a monomorphism between $h_r$%the discrete group $h_r=\{[l,m,k]\; |\; l,m,k \in
  %on the convergence result (\ref{convresult}) we define the mapping $\phi: h_r \to H_{3}$ by
%\[
%\phi[l,m,k]=\left(\frac{l}{K}, \frac{mK}{P},\frac{k}{Q}\right)
%\]
%which sets a monomorphism between $h_r$%the discrete group $h_r=\{[l,m,k]\; |\; l,m,k \in \mathbb{Z}\}$
%%which is equipped with group product
%%\[
%%[l,m,k][l',m',k']=[l+l' , m+m' , k+k' + \frac{Q}{2P}(ml'-m'l) ].
%%\]
%and the continuous Heisenberg group $H_3$:
%\[
%\begin{array}{ll}
%\phi[l,m,k]\phi[l',m',k'] &=\left(\frac{l}{K}, \frac{mK}{P},\frac{k}{Q}\right)\left(\frac{l}{K}, \frac{mK}{P},\frac{k}{Q}\right)\\
% &=\left(\frac{l+l'}{K}, \frac{(m+m')K}{P}, \frac{k+k' + \frac{Q}{2P}(ml'-lm')}{Q} \right)=\phi[[l,m,k][l',m',k']].
%\end{array}
%\]
%Here we
The mapping $\phi$ maps the discrete variables on a uniform grid in the continuous domain:
\[
\begin{array}{lll}
s \in [0,1) \leftrightarrow k \in [0,Q) \cap \mathbb{Z}, &
p \in [0,1) \leftrightarrow l \in [0,K) \cap \mathbb{Z}, &
q \in [0,N) \leftrightarrow m \in [0,M) \cap \mathbb{Z}.
\end{array}
\]
%The group $h_r$ is not isomorphic to a subgroup of $H_{r}$ since we have periodicity in both $l$ and $m$. However, if again $\frac{N}{2M} \in \mathbb{N}$, then it equals the quotient
%$\gothic{h}_r:=h_{r}/[K\mathbb{Z}, M\mathbb{Z}, Q \mathbb{Z}]$ of the discrete group $h_r$ with the normal subgroup $[K\mathbb{Z}, M\mathbb{Z}, Q \mathbb{Z}]$, %$\gothic{g}$ generated by the elements $(p=\frac{1}{K},0,0)$ and %$(0,q= %\frac{K}{P},0)$ with its normal subgroup $\gothic{h}=??=\{(l %\frac{1}{K},m\frac{K}{P}, k_{1} \frac{1}{K}+k_{2}\frac{1}{2P}) \; | \; %l,m,k_1,k_{2} \in \mathbb{Z}\}$
%where even $\frac{N}{2M}$ and $N$ implies $[K \mathbb{Z}, M \mathbb{Z},Q \mathbb{Z}]=[K \mathbb{Z}, M \mathbb{Z},Q \mathbb{Z}-\frac{Q}{2P}MK \mathbb{Z}]$. % (recall that $MK=PN$).

On the quotient-group $\gothic{h}_r$ we define the forward left-invariant vector fields on discrete Gabor-transforms as follows (where we again use (\ref{relation}) and (\ref{LNK})):
{\small
\begin{equation} \label{forwarddiff}
\begin{array}{lll}
(\mathcal{A}_{1}^{D^{+}} W_{\PPP}^{D}\ul{f})[l,m,k] &= K\, ({\rm d}\mathcal{R}^{D^{+}}[1,0,0] W_{\PPP}^{D}\ul{f})[l,m,k] &= \frac{W_{\PPP}^{D}\ul{f}([l,m,k][1,0,0])-W_{\PPP}^{D}\ul{f}[l,m,k]}{K^{-1}} \\&=
\frac{e^{-\frac{\pi i m }{P}}W_{\PPP}^{D}\ul{f}[l+1,m,k]-W_{\PPP}^{D}\ul{f}[l,m,k]}{K^{-1}}&= \frac{e^{-\frac{\pi i m L }{M}}W_{\PPP}^{D}\ul{f}[l+1,m,k]-W_{\PPP}^{D}\ul{f}[l,m,k]}{K^{-1}} \\[7pt]
(\mathcal{A}_{2}^{D^{+}} W_{\PPP}^{D}\ul{f})[l,m,k] &= \frac{M}{N}\,({\rm d}\mathcal{R}^{D^{+}}[0,1,0] W_{\PPP}^{D}\ul{f})[l,m,k] &= \frac{e^{+\frac{\pi i l }{P}}W_{\PPP}^{D}\ul{f}[l,m+1,k]-W_{\PPP}^{D}\ul{f}[l,m,k]}{K\, P^{-1}},  \\
&\ &= \frac{e^{+\frac{\pi i l L }{M}}W_{\PPP}^{D}\ul{f}[l,m+1,k]-W_{\PPP}^{D}\ul{f}[l,m,k]}{N\, M^{-1}},
 \\[7pt]
(\mathcal{A}_{3}^{D^{+}} W_{\PPP}^{D}\ul{f})[l,m,k] &=Q({\rm d}\mathcal{R}^{D^{+}}[0,0,1] W_{\PPP}^{D}\ul{f})[l,m,k] &= \frac{W_{\PPP}^{D}\ul{f}[l,m,k+1]-W_{\PPP}^{D}\ul{f}[l,m,k]}{Q^{-1}} \\
 & \ &= Q( e^{\frac{-2\pi i}{Q}}-1)W_{\PPP}^{D}\ul{f}[l,m,k]),
\end{array}
\end{equation}
}
and the backward discrete left-invariant vector fields
\begin{equation}\label{backwarddiff}
\begin{array}{ll}
(\mathcal{A}_{1}^{D^{-}} W_{\PPP}^{D}\ul{f})[l,m,k] &= ({\rm d}\mathcal{R}^{D^{-}}[1,0,0] W_{\PPP}^{D}\ul{f})[l,m,k]
%= \frac{W_{\PPP}^{D}\ul{f}([l,m,k])-W_{\PPP}^{D}\ul{f}[l,m,k][-1,0,0]}{K^{-1}} \\
%\\ &
=
\frac{W_{\PPP}^{D}\ul{f}[l,m,k]-e^{+\frac{\pi i m L}{M}}W_{\PPP}^{D}\ul{f}[l-1,m,k]}{K^{-1}}, \\
(\mathcal{A}_{2}^{D^{-}} W_{\PPP}^{D}\ul{f})[l,m,k] &= ({\rm d}\mathcal{R}^{D^{+}}[0,1,0] W_{\PPP}^{D}\ul{f})[l,m,k]
%\\
%&
= \frac{W_{\PPP}^{D}\ul{f}[l,m,k]-e^{-\frac{\pi i l L}{M}}W_{\PPP}^{D}\ul{f}[l,m-1,k]}{N\, M^{-1}}, \\
(\mathcal{A}_{3}^{D^{-}} W_{\PPP}^{D}\ul{f})[l,m,k] &=({\rm d}\mathcal{R}^{D^{+}}[0,0,1] W_{\PPP}^{D}\ul{f})[l,m,k]
%\\ &
= \frac{W_{\PPP}^{D}\ul{f}[l,m,k]-W_{\PPP}^{D}\ul{f}[l,m,k-1]}{Q^{-1}}= Q(1- e^{\frac{2\pi i}{Q}})W_{\PPP}^{D}\ul{f}[l,m,k]\ .
\end{array}
\end{equation}
\begin{remark}
With respect to the step-sizes in (\ref{forwarddiff}) and (\ref{backwarddiff}) we have set $p=\frac{l}{K}$, $q= \frac{m}{M}N$, $\xi=\frac{n}{N}$, $s=\frac{k}{Q}$, so that the actual discrete steps are $\Delta p= K^{-1}$, $\Delta q= N\, M^{-1}$ and $\Delta s= Q^{-1}$.
This discretization is chosen such that both the continuous Gabor transform and the continuous
left-invariant vector fields follow from their discrete counterparts by $N\to \infty$, e.g. recall Eq.\!~(\ref{convresult}).
\end{remark}
 Akin to the continuous case we use the following discrete version $S^{D}$ of the operator $\mathcal{S}$ that maps a Gabor transform $W_{\PPP}^{D}$ onto its
phase space representation $G_{\PPP}^{D}$:
\[
G_{\PPP}^{D}\ul{f}[l,m]:=(S^{D} W^{D}_{\PPP})\ul{f}[l,m]=W^{D}_{\PPP}\ul{f}[l,m,-\frac{Qlm}{2P}], \qquad P=M/L.
\]
The inverse is given by $W_{\PPP}^{D}\ul{f}[l,m,k] =((S^{D})^{-1}G_{\PPP}^{D}\ul{f})[l,m,k]= e^{-2\pi i \left(\frac{k}{Q}+ \frac{lmL}{2M} \right)}G_{\PPP}^{D}\ul{f}[l,m]$.

Again we can use the conjugation with $S^{D}$ to map the left-invariant discrete vector fields $\{\mathcal{A}_{i}^{D^{\pm}}\}_{i=1}^{3}$ to the corresponding discrete vector fields on the discrete phase space: $\tilde{\mathcal{A}}^{D^{\pm}}_{i}= (S^{D})\circ  \mathcal{A}^{D^{\pm}}_{i} \circ (S^{D})^{-1}$.
A brief computation yields the following forward left-invariant differences
\begin{equation}\label{DiscretegenFS1}
\begin{array}{ll}
(\tilde{\mathcal{A}}_{1}^{D^{+}} G_{\PPP}^{D}\ul{f})[l,m]&=\frac{e^{\frac{-2\pi i Lm}{M}} (G_{\PPP}^{D}\ul{f})[l+1,m]- G^{D}_{\PPP}\ul{f}[l,m]}{K^{-1}} \\
(\tilde{\mathcal{A}}_{2}^{D^{+}} G^{D}_{\PPP}\ul{f})[l,m] &= M \, N^{-1}\,(G^{D}_{\PPP}\ul{f}[l,m+1]-G^{D}_{\PPP}\ul{f}[l,m]) \\
(\tilde{\mathcal{A}}_{3}^{D^{+}} G^{D}_{\PPP}\ul{f})[l,m] &=Q( e^{\frac{-2\pi i}{Q}}-1)G^{D}_{\PPP}\ul{f}[l,m]\\
\end{array}
\end{equation}
and the following backward left-invariant differences:
\begin{equation} \label{DiscretegenFS2}
\begin{array}{ll}
(\tilde{\mathcal{A}}_{1}^{D^{-}} G^{D}_{\PPP}\ul{f})[l,m]&=\frac{ (G^{D}_{\PPP}\ul{f})[l,m]- e^{\frac{2\pi i Lm}{M}} G^{D}_{\PPP}\ul{f}[l-1,m]}{K^{-1}} \\
(\tilde{\mathcal{A}}_{2}^{D^{-}} G^{D}_{\PPP}\ul{f})[l,m] &= M\, N^{-1}(G^{D}_{\PPP}\ul{f}[l,m]-G^{D}_{\PPP}\ul{f}[l,m-1]) \\
(\tilde{\mathcal{A}}_{3}^{D^{-}} G^{D}_{\PPP}\ul{f})[l,m] &= Q(1- e^{\frac{2\pi i}{Q}})G^{D}_{\PPP}\ul{f}[l,m]\ . \\
\end{array}
\end{equation}
%\[
%
%\]
The discrete operators are defined on the discrete quotient group $\gothic{h}_{r}$ and do not involve approximations in the setting of discrete Gabor transforms. They are first order approximation
of the corresponding continuous operators on $H_{r}$ as we motivate next. For $f$ compactly supported on $[0,1]$ and both $f$ and $\psi$ Riemann-integrable on $\R$:
\begin{equation} \label{een}
\begin{array}{ll}
\tilde{\mathcal{A}_{1}} \mathcal{G}_{\psi}f(p=\frac{l}{K},q=\frac{mK}{P}) &=(\partial_{p}-2\pi q)(e^{2\pi ipq} \int \limits_{\R} \overline{\psi(\xi-p)} f(\xi) e^{-2\pi i \xi q} \,{\rm d}\xi)(\frac{l}{K},\frac{mK}{P})  \\
 &=- e^{\frac{2\pi il}{P}} \int \limits_{\R} \overline{\psi}'(\xi-\frac{l}{K}) f(\xi) e^{-\frac{2\pi i nm K}{NP}}\, {\rm d}\xi \\
 &= O(\frac{1}{N}) - \frac{1}{N}e^{\frac{2\pi i\, lm}{P}}
 \sum \limits_{n=0}^{N-1} \overline{\psi}'\left(\frac{n}{N}-\frac{l}{K}
 \right) f\left( \frac{n}{N}\right)  e^{-\frac{2\pi i nm}{M}}.
\end{array}
\end{equation}
Moreover, we have
\[
[G^{D}_{\PPP}\ul{f}](l,m)= \frac{1}{N} e^{\frac{2\pi i ml}{P}} \sum \limits_{n=0}^{N-1} e^{-\frac{2\pi i nm}{M}} \overline{\psi}\left(\frac{n}{N}-\frac{l}{K} \right) f\left(\frac{n}{N} \right)
\]
so that straightforward computation yields
\begin{equation}\label{twee}
\begin{array}{ll}
\tilde{\mathcal{A}}_{1}^{D^{+}} G^{D}_{\PPP} \ul{f}[l,m] &= \frac{1}{N} e^{2\pi i \frac{lm}{P}} \sum \limits_{n=0}^{N-1}
\frac{\overline{\psi}\left(\frac{n}{N}-\frac{l+1}{K} \right)- \overline{\psi}\left( \frac{n}{N}-\frac{l}{K}\right)}{K^{-1}}\;  f\left(\frac{n}{N} \right)\, e^{-\frac{2\pi i n m}{N}}
\\
 &= O\left(\frac{1}{K} \right)O\left(\frac{1}{N} \right) -\frac{1}{N}\, e^{\frac{2\pi i\, lm}{P}}
 \sum \limits_{n=0}^{N-1} \overline{\psi}'\left(\frac{n}{N}-\frac{l}{K}
 \right) f\left( \frac{n}{N}\right)  e^{-\frac{2\pi i nm}{M}}
% \\
% &= O\left(\frac{1}{K}\right) + O\left(\frac{1}{N}\right) + %\tilde{\mathcal{A}_{1}} \mathcal{G}_{\psi}f(p=\frac{l}{K},q=\frac{mK}{P}).
\end{array}
\end{equation}
So from (\ref{een}) and (\ref{twee}) we deduce that
\[
\tilde{\mathcal{A}}_{1}^{D^{+}} G^{D}_{\PPP} \ul{f}[l,m] = O\left(\frac{1}{N}\right) + \tilde{\mathcal{A}_{1}} \mathcal{G}_{\psi}f(p=\frac{l}{K},q=\frac{mK}{P}).
\]
So clearly the discrete left-invariant vector fields acting on the discrete Gabor-transforms converge to the continuous vector fields acting on the continuous Gabor transforms pointwise as $N \to \infty$.

%\[
%\begin{array}{ll}
%\mathcal{A}_{1}^{D^{+}} W_{\PPP}^{D}(l,m,k) &=(1/L)(\frac{W_{\PPP}^{D}[l+1,m,k]-W_{\PPP}^{D}[l,m,k]}{L}) \\
% &+ 2 \frac{m}{MQ}\left(W_{\PPP}^{D}\ul{f}(l,m,k+1)-W_{\PPP}^{D}\ul{f}(l,m,k)\right)+ O(L) \\
%  &= (\partial_{x}+2\omega \partial_{t})\mathcal{W}_{\psi}f(lL,\frac{m}{M},\frac{k}{Q}) + O(L) + O((MQ)^{-1}).
% \end{array}
%\]
In our algorithms it is essential that one works on the finite group with corresponding left-invariant vector fields. This is simply due to the fact that one computes finite Gabor-transforms (defined on the group $\gothic{h}_r$)
to avoid sampling errors on the grid. Standard finite difference approximations of the continuous left-invariant vector fields do not  appropriately deal with phase oscillations in the (discrete) Gabor transform.
%
%So, rather than considering the discrete left-invariant vector fields $\tilde{\mathcal{A}}_{i}^{D^{+}}$, $i=1,2,3$, acting on $G^{D}_{\PPP} \ul{f}$ as an approximation of the continuous left-invariant vector fields $\tilde{\mathcal{A}}_{i}$, $i=1,2,3$ acting on $\mathcal{G}_{\psi} f$, one should consider
%the continuous left-invariant vector fields $\tilde{\mathcal{A}}_{i}$ acting on $\mathcal{G}_{\psi} f$ as an approximation of the discrete left-invariant vector fields $\tilde{\mathcal{A}}_{i}^{D^{+}}$, $i=1,2,3$ acting on $G^{D}_{\PPP} \ul{f}$.
%\begin{definition}
%\end{definition}
\begin{remark}
In the PDE-schemes which we will present in the next sections, such as for example the diffusion scheme in Section \ref{ch:covdiff}, the solutions will leave the space of Gabor-transforms. In such cases one has to apply a left-invariant finite difference to a smooth function $U \in \mathbb{L}_{2}(H_r)$ defined on the Heisenberg-group $H_r$ or, equivalently, one has to apply a finite difference to a smooth function $\tilde{U} \in \mathbb{L}_{2}(\R^2)$ defined on phase space. If $U$ is not the Gabor transform of an image it is usually inappropriate to use the final results in (\ref{backwarddiff}) and (\ref{forwarddiff}) on the group $H_{r}$. Instead one should just use
\begin{equation}\label{right}
\begin{array}{ll}
(\mathcal{A}_{1}^{D^{+}} U)[l,m,k] &= ({\rm d}\mathcal{R}^{D^{+}}[1,0,0] U)[l,m,k]= \frac{U[[l,m,k][1,0,0]]-U[l,m,k]}{K^{-1}} \\
(\mathcal{A}_{1}^{D^{-}} U)[l,m,k] &= ({\rm d}\mathcal{R}^{D^{-}}[1,0,0] U)[l,m,k]= \frac{U[l,m,k]-U[[l,m,k][-1,0,0]]}{K^{-1}} \\
(\mathcal{A}_{2}^{D^{+}} U)[l,m,k] &= ({\rm d}\mathcal{R}^{D^{+}}[0,1,0] U)[l,m,k]= \frac{U[[l,m,k][0,1,0]]-U[l,m,k]}{N\, M^{-1}} \\
(\mathcal{A}_{2}^{D^{-}} U)[l,m,k] &= ({\rm d}\mathcal{R}^{D^{-}}[0,1,0] U)[l,m,k]= \frac{U[l,m,k]-U[[l,m,k][0,-1,0]]}{N\, M^{-1}} \\[7pt]
\end{array}
\end{equation}
which does not require any interpolation between the discrete data iff $\frac{Q}{2P} \in \mathbb{N}$.
The left-invariant operators on phase space (\ref{DiscretegenFS1}) and (\ref{DiscretegenFS2}) are naturally extendable to $\mathbb{L}_{2}(\R^2)$. For example,
$(\mathcal{A}_{1}^{D^{+}} U)[l,m]=[\mathcal{S}^{D}\circ \mathcal{A}_{1}^{D^{+}}\circ (\mathcal{S}^{D})^{-1}\tilde{U}][l,m]=
K(e^{-\frac{2\pi im}{P}}\tilde{U}[l+1,m]-\tilde{U}[l,m])$ for all $\tilde{U} \in \ell_{2}(\{0,\ldots,K-1\}\times\{0,\ldots,M-1\})$.
\end{remark}

\section{Algorithm for the PDE-approach to Differential Reassignment \label{ch:upwind}}

Here we provide an explicit algorithm on the discrete Gabor transform $G^{D}_{\PPP}\ul{f}$ of the discrete signal $\ul{f}$, that consistently corresponds to the theoretical PDE's on the continuous case as proposed in \cite{Daudet}, i.e. convection equation (\ref{conpde}) where we apply the first choice (\ref{conchoice}). Although that the PDE studied in \cite{Daudet} is not as simple as the second approach in (\ref{conchoice}) (which corresponds to a standard erosion step on the absolute value $|\mathcal{G}_{\psi}f|$ followed by a restoration of the phase afterwards) we do provide an explicit numerical scheme of this PDE, where we stay entirely in the \emph{discrete phase space}.

It should be stressed that taking straightforward central differences of the continuous differential operators of section \ref{ch:reass} does not work, due to the fast oscillations (of the phase) of Gabor transforms. We need the lef-invariant differences on discrete Heisenberg groups discussed in the previous subsection.
%For details and non-trivial motivation of left-invariant differences on discrete Heisenberg groups see %\cite{DuitsFuehrJanssen}.
\\
%QUESTION TO BART: ARE THE CAUCHY RIEMANN EQUATIONS STILL VIOLATED ?? ALSO WHEN I STAY AWAY FROM THE BOUNDARY ???
%\\
\\
{\small
\textbf{Explicit upwind scheme with left-invariant finite differences} in pseudo-code for $\mathcal{M}=1$ \\
\ \\
For $l=1,\ldots, K-1$, $m=1,\ldots M-1$ set $\tilde{W}[l,m,0]:= G^{D}_{\PPP}\ul{f}[l,m]$. \\
For $t=1,\ldots, T$ \\
 \  \ For $l=0,\ldots, K-1$, for $m=1, \ldots, M-1$ set \\
 \  \ $\tilde{v}^{1}[l,m]:= -\frac{aK}{2}(\log |\tilde{W}|[l+1,m,t=0]- \log |\tilde{W}|[l-1,m,t=0])$  \\
 \ \  $\tilde{v}^{2}[l,m]:= -\frac{aM}{2}(\log |\tilde{W}|[l,m+1,t=0]|-\log|\tilde{W}|[l,m-1,t=0])$ \\
 \  \  $\tilde{W}[l,m,t]:=\tilde{W}[l,m,t-1]+ K \, \Delta t  \left( z^{+}(\tilde{v}^{1})[l,m] \; [\tilde{\mathcal{A}}_{1}^{D^{-}}\tilde{W}][l,m,t] +
z^{-}(\tilde{v}^{1})[l,m] \; [\tilde{\mathcal{A}}_{1}^{D^{+}}\tilde{W}][l,m,t]\right) + $ \\
\ \ $M\Delta t \left( z^{+}(\tilde{v}^{2})[l,m] \; [\tilde{\mathcal{A}}_{2}^{D^{-}}\tilde{W}][l,m,t] +
z^{-}(\tilde{v}^{2})[l,m] \; [\tilde{\mathcal{A}}_{2}^{D^{+}}\tilde{W}][l,m,t]\right)$.
}
 \\ \\
Explanation of all involved variables:
{\small
\[
\boxed{
\begin{array}{ll}
l & \textrm{discrete position variable $l=0,\ldots, K-1$.}\\
m & \textrm{discrete frequency variable $m=1, \ldots, M-1$.} \\
t & \textrm{discrete time $t=1,\ldots T$, where $T$ is the stopping time.} \\
%n & \textrm{discrete dummy intragtion}
\PPP &\textrm{discrete kernel }\PPP=\PPP_{a}^{C}=\{\psi_{a}(n N^{-1})\}_{n=-(N-1)}^{N-1} \textrm{ or } \PPP=\{\PPP_{a}^{D}[n]\}_{n=-(N-1)}^{N-1} \textrm{see below.} \\
G^{D}_{\PPP}\ul{f}[l,m] & \textrm{discrete Gabor transform computed by diagonalization via Zak transform \cite{GuidoJanssen}.} \\
\tilde{W}[l,m,t]& \textrm{discrete evolving Gabor transform evaluated at position $l$, frequency $m$ and time $t$.} \\
\tilde{\mathcal{A}}_{i}^{D^{\pm}}& \textrm{forward (+), backward (-) left-invariant position ($i=1$) and
frequency ($i=2$) shifts. }
%\end{array}
%\]
 \\
z^{\pm}& z^{+}(\phi)[l,m,t]=\max\{\phi(l,m,t),0\}, z^{-}(\phi)[l,m,t]=\min\{\phi(l,m,t),0\} \textrm{ for upwind.}
\end{array}
}
\]
}
%Recall Eq.~(\ref{DiscretegenFS1}).
%{\small
%\begin{equation}\label{DiscretegenFS1}
%%\begin{array}{ll}
%\begin{array}{ll}
%(\tilde{\mathcal{A}}_{1}^{D^{+}} \tilde{\Phi})[l,m] = K(e^{-\frac{2\pi i L m}{M}}\tilde{\Phi}[l+1,m]-\tilde{\Phi}[l,m]), & \ \ \
%(\tilde{\mathcal{A}}_{1}^{D^{-}} \tilde{\Phi})[l,m] = K(\tilde{\Phi}[l,m]-e^{\frac{2\pi i Lm}{M}}\tilde{\Phi}[l,m]), \\
%(\tilde{\mathcal{A}}_{2}^{D^{+}} \tilde{\Phi})[l,m] = MN^{-1}(\tilde{\Phi}[l,m+1]-\tilde{\Phi}[l,m]), &\ \ \
%(\tilde{\mathcal{A}}_{2}^{D^{-}} \tilde{\Phi})[l,m] = MN^{-1}(\tilde{\Phi}[l,m]-\tilde{\Phi}[l,m-1]), %\\[7pt]
%\end{array}
%\end{equation}
%}
%The discrete Gabor transform equals
%$
%G^{D}_{\PPP}\ul{f}[l,m]=\frac{1}{N} \sum \limits_{n=0}^{N-1} \overline{\PPP[n-lL]}f[n]e^{-\frac{2\pi i n(n-lL)}{M}}\ ,
%$
%where $M/L$ denotes the (integer) oversampling factor and $N=KL$.
The discrete Cauchy Riemann kernel $\PPP^{D}_{a}$ is derived in \cite{DuitsFuehrJanssen} and satisfies the system
{\small
\begin{equation}\label{discreteCR}
\begin{array}{l}
\forall_{l=0,\ldots,K\!-\!1}\forall_{m=0,\ldots,M\!-\!1} \forall_{\ul{f} \in \ell_{2}(I)}\; : \;
\frac{1}{a}(\tilde{\mathcal{A}}_{2}^{D^{+}}\! \!+\!\tilde{\mathcal{A}}_{2}^{D^{-}}) + i \, a (\tilde{\mathcal{A}}_{1}^{D^{+}}\! \!+\!\tilde{\mathcal{A}}_{1}^{D^{-}})(G_{\PPP^{D}_{a}}^{D}\ul{f})[l,m]=0 \ ,
\end{array}
\end{equation}
}
which has a unique solution in case of extreme oversampling $K=M=N$, $L=1$.
%\begin{figure}
\begin{figure}[b]
\centerline{
\hfill
\includegraphics[width=0.29\hsize]{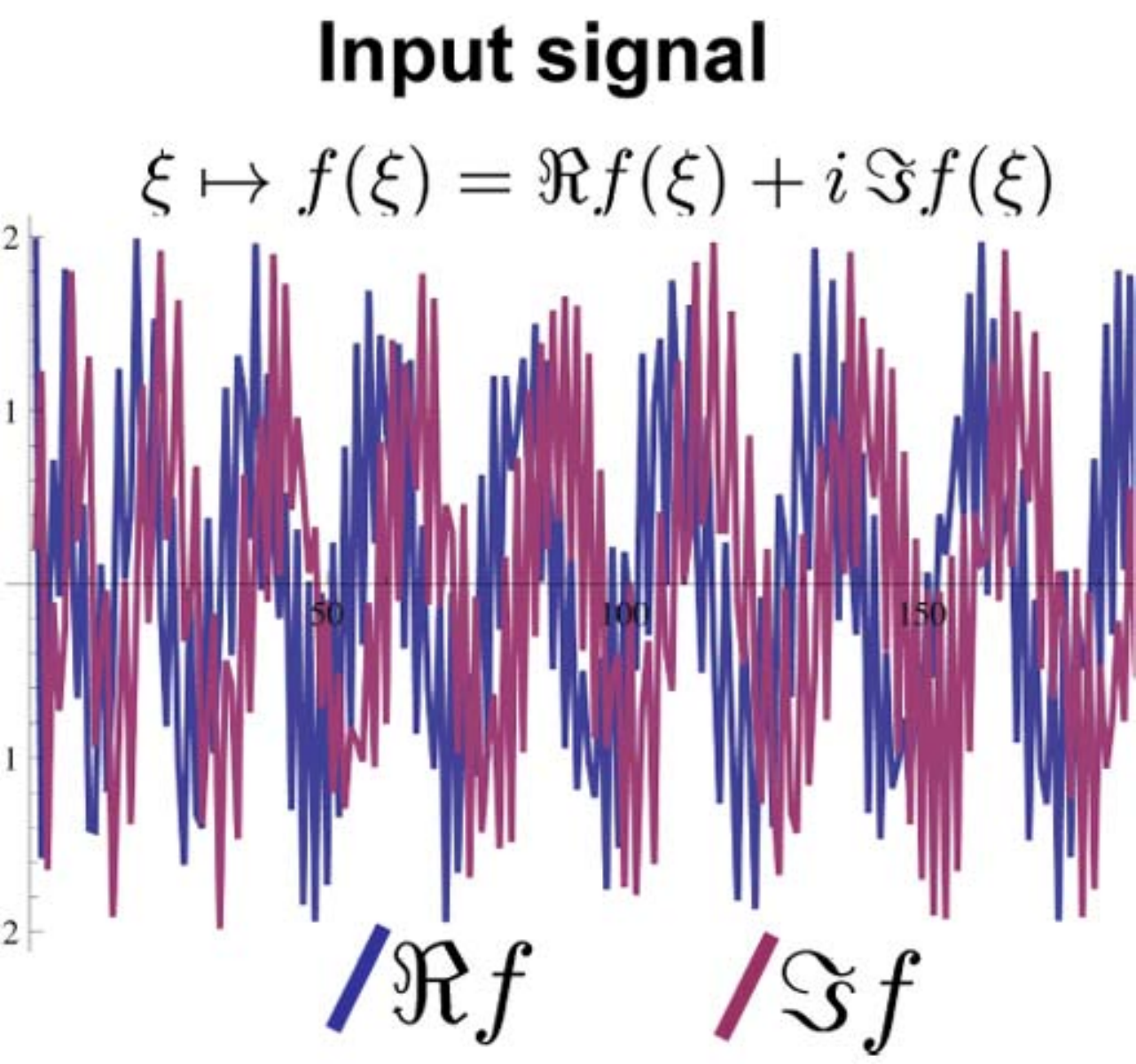}
\hfill
\includegraphics[width=0.54\hsize]{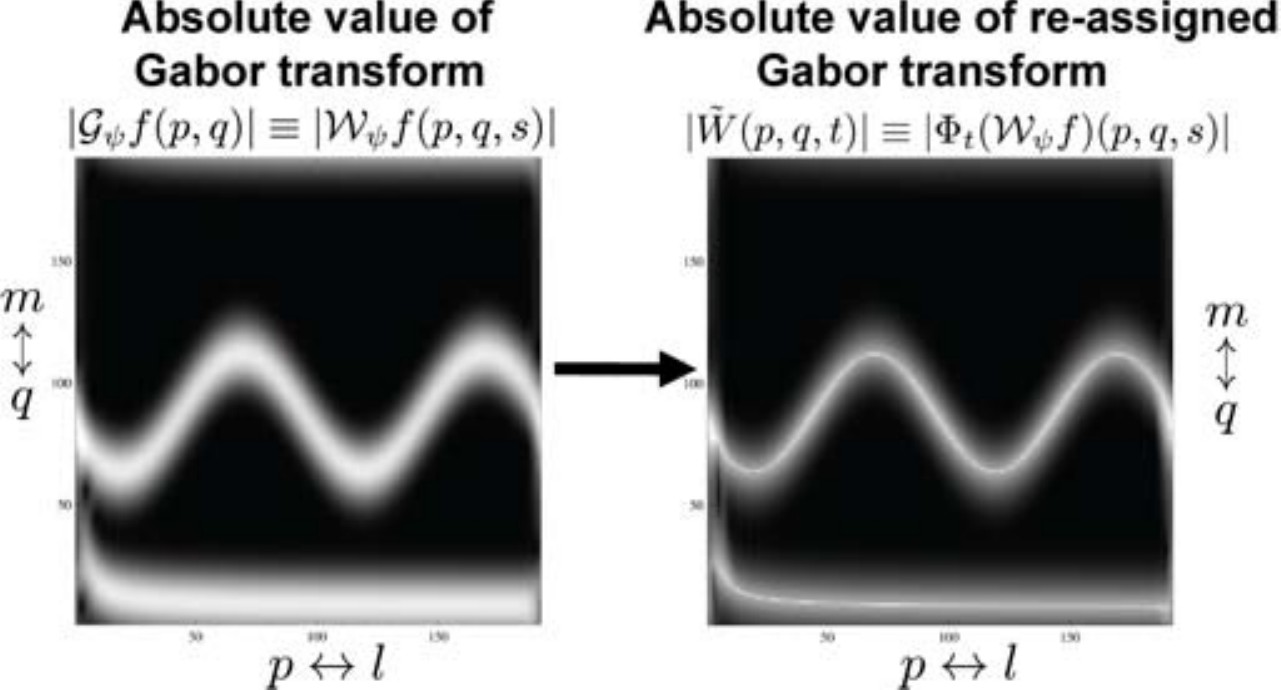}
\hfill
}
\caption{Illustration of reassignment by adaptive phase-invariant convection explained in Section \ref{ch:reass}, using the upwind scheme of subsection \ref{ch:upwind} applied on
 a Gabor transform. %Left: input signal, middle: absolute value of the corresponding Gabor transform, right: result of re-assignment in the %Gabor domain.
 }\label{fig:reass}
\end{figure}
%\vspace{-01cm} \mbox{}
\subsection{Evaluation of Reassignment \label{ch:evaluate}}

We distinguished between two approaches to apply left-invariant adaptive convection on discrete Gabor-transforms (that we diagonalize by discrete Zak transform \cite{GuidoJanssen}, recall subsection \ref{ch:zak}).
%\footnote{The induced frame operator can be efficiently diagonalized by Zak-transform, \cite{GuidoJanssen}, %boiling down to diagonalization of inverse Fourier transform on $H_{r}$, \cite[ch:2.3]{DuitsFuehrJanssen}. We %used this in our algorithms. }
Either we apply the numerical upwind PDE-scheme described in subsection \ref{ch:upwind} using the discrete left-invariant vector fields (\ref{DiscretegenFS1}), or we apply erosion (\ref{erosion}) on the modulus and restore the phase afterwards.
%These two approaches corresponds to respectively the first and second choice in (\ref{conchoice}). All implementations take place in the %discrete phase space.
Within each of the two approaches, we can use the discrete Cauchy-Riemann kernel $\PPP^{D}_{a}$ or the sampled continuous Cauchy-Riemann kernel $\PPP^{C}_{a}$. %This brings the total amount of methods to four.

\begin{figure}
\centerline{
\includegraphics[width=0.34\hsize]{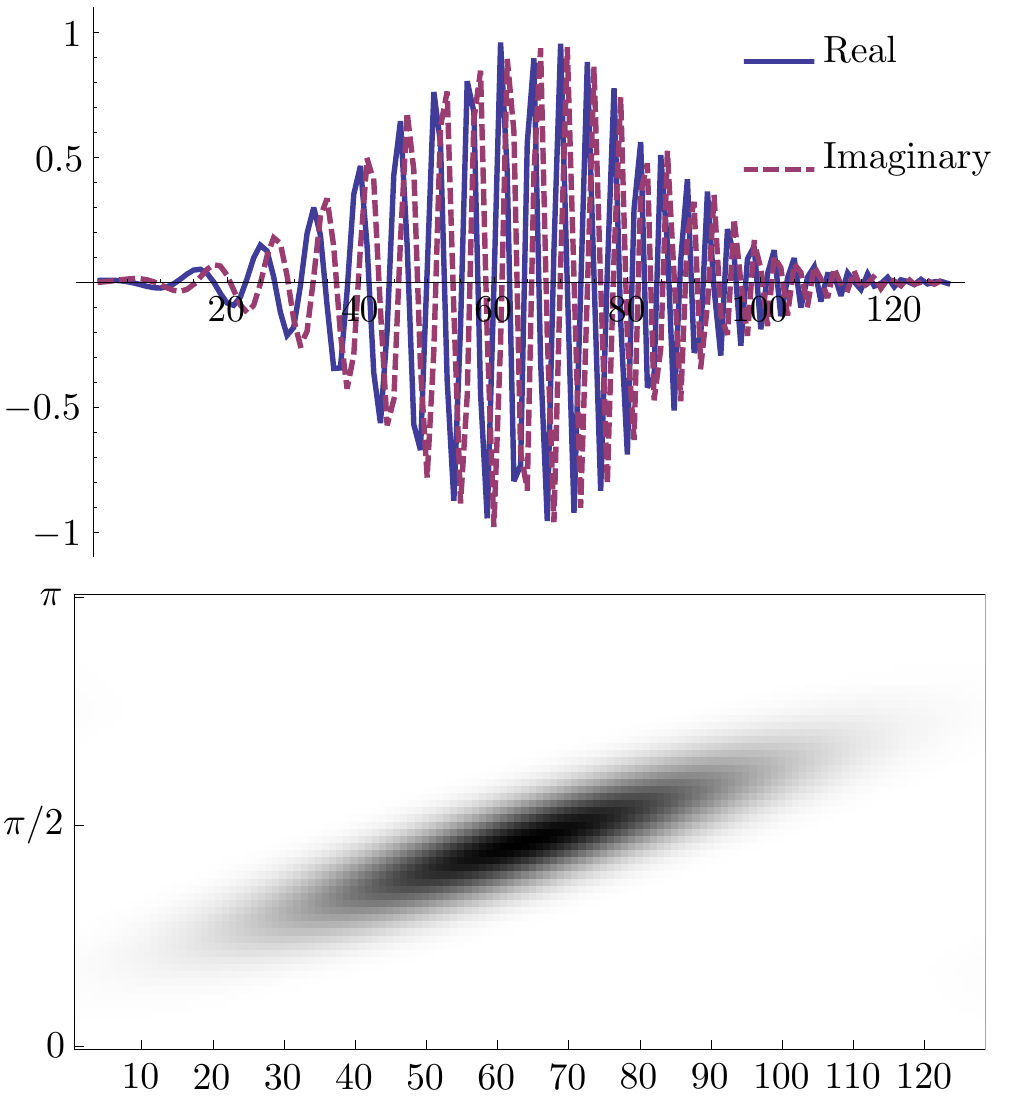}
\includegraphics[width=0.66\hsize]{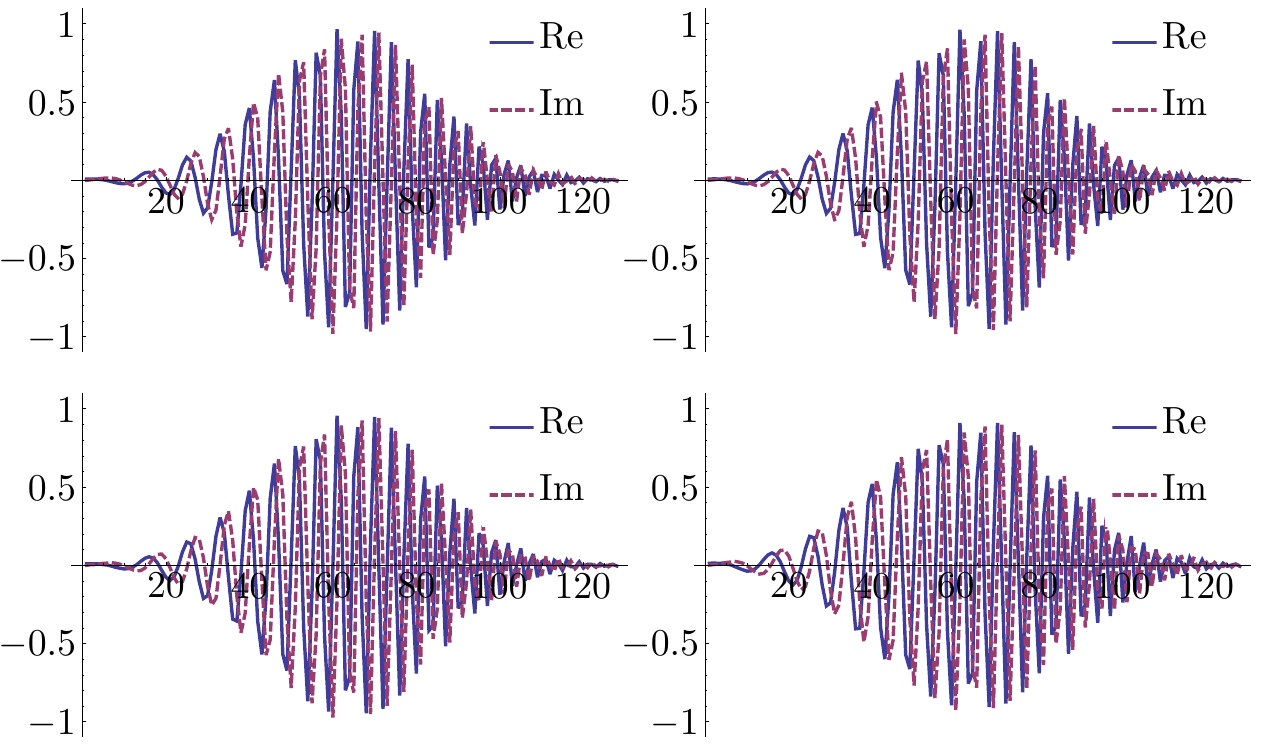}
}
\caption{
Reconstructions of the reassigned Gabor transforms of the original signal
that is depicted on the top left whose absolute value of the Gabor transform is depicted on bottom left. In the right: 1st row
corresponds to reassignment by the upwind scheme ($\mathcal{M}=1$) of subsection \ref{ch:upwind}, where again left we used $\PPP^{C}_{a}$ and right we used $\PPP^{D}_{a}$.
Parameters involved are grid constants
$K = M = N = 128$, window scale $a = 1/6$, time step $\Delta t = 10^{-3}$ and time $t = 0.1$.
2nd row to
reassignment by morphological erosion where in the left we used kernel $\PPP^{C}_{a}$
and in the right we used $\PPP^{D}_{a}$.
The goal of reassignment is achieved; all reconstructed signals are close to the original signal, whereas their corresponding Gabor transforms depicted in Fig. \ref{fig:5} are much sharper than the absolute value of the Gabor transform of the original depicted on the bottom left of this figure. }\label{fig:4}
\end{figure}

\begin{figure}
\centerline{
\includegraphics[width=0.4 \hsize]{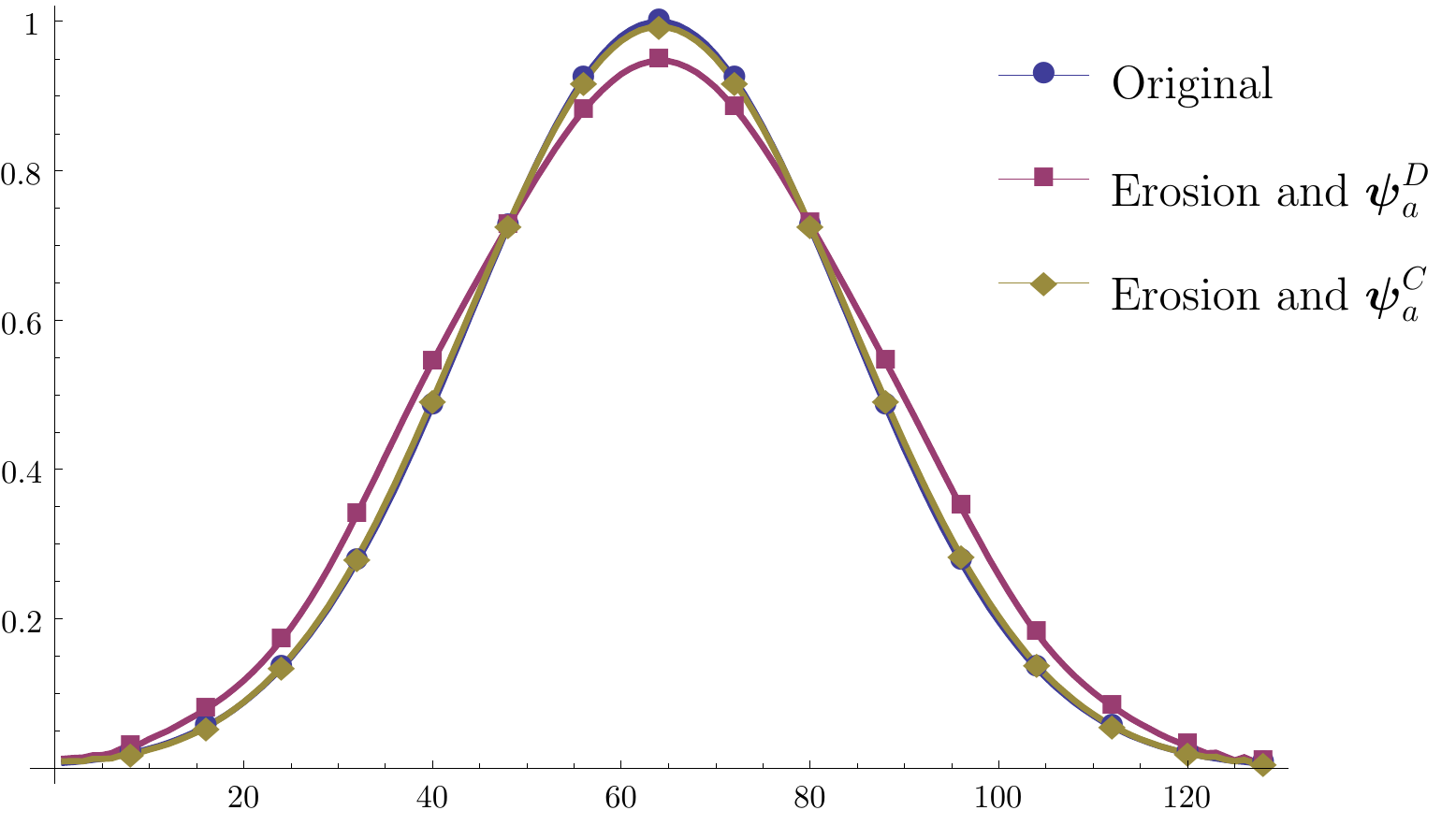}
}
\caption{The modulus of the signals in the bottom row of Fig. \ref{fig:4}. For erosion (\ref{erosion}) $\PPP^{C}_{a}$ performs better than erosion applied on a Gabor transform constructed by $\PPP^{D}_{a}$.}\label{fig:8}
\end{figure}

\begin{figure}
%\centerline{
%\includegraphics[width=0.5\hsize]{GT}
%}
\centerline{
\includegraphics[width=0.3\hsize]{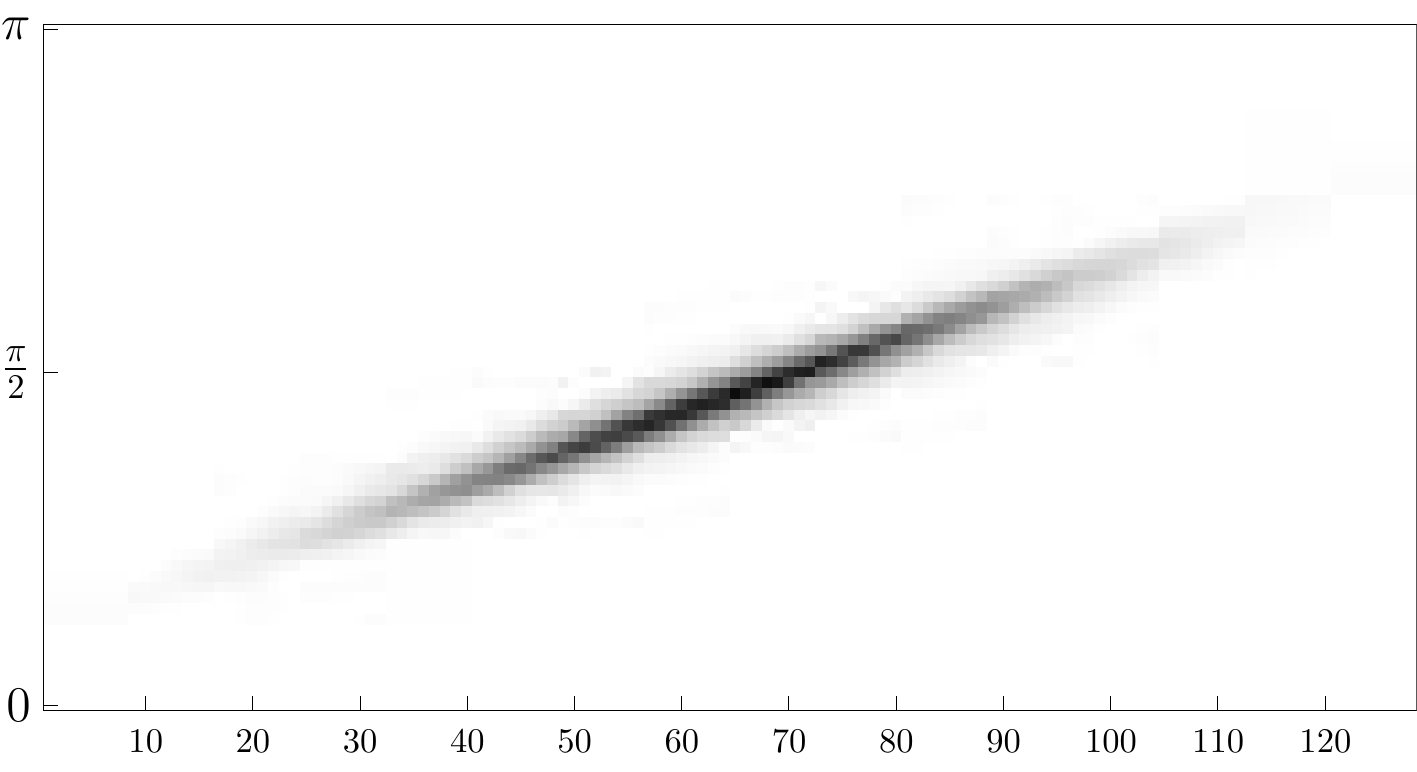}
\includegraphics[width=0.3\hsize]{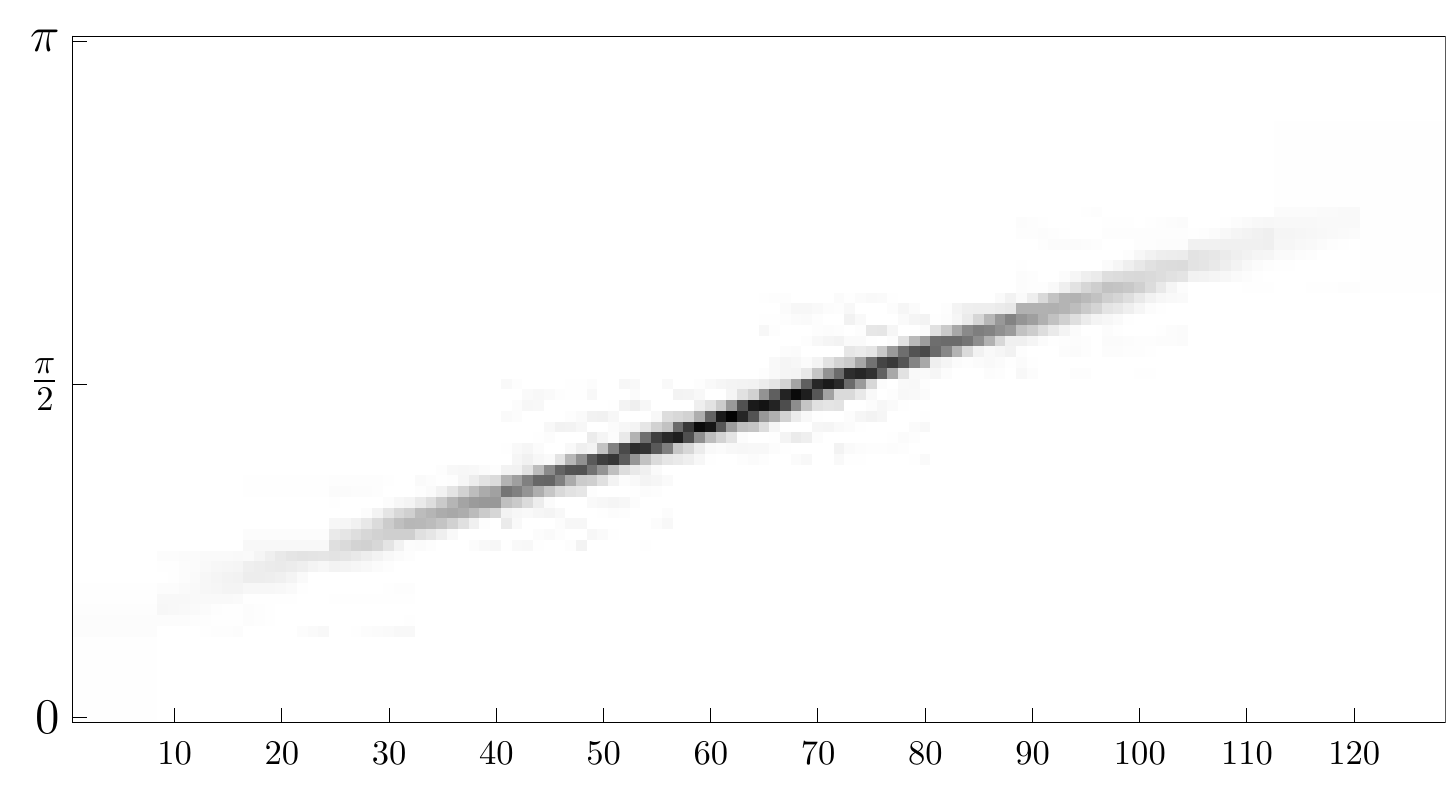}
}
\centerline{
\includegraphics[width=0.3\hsize]{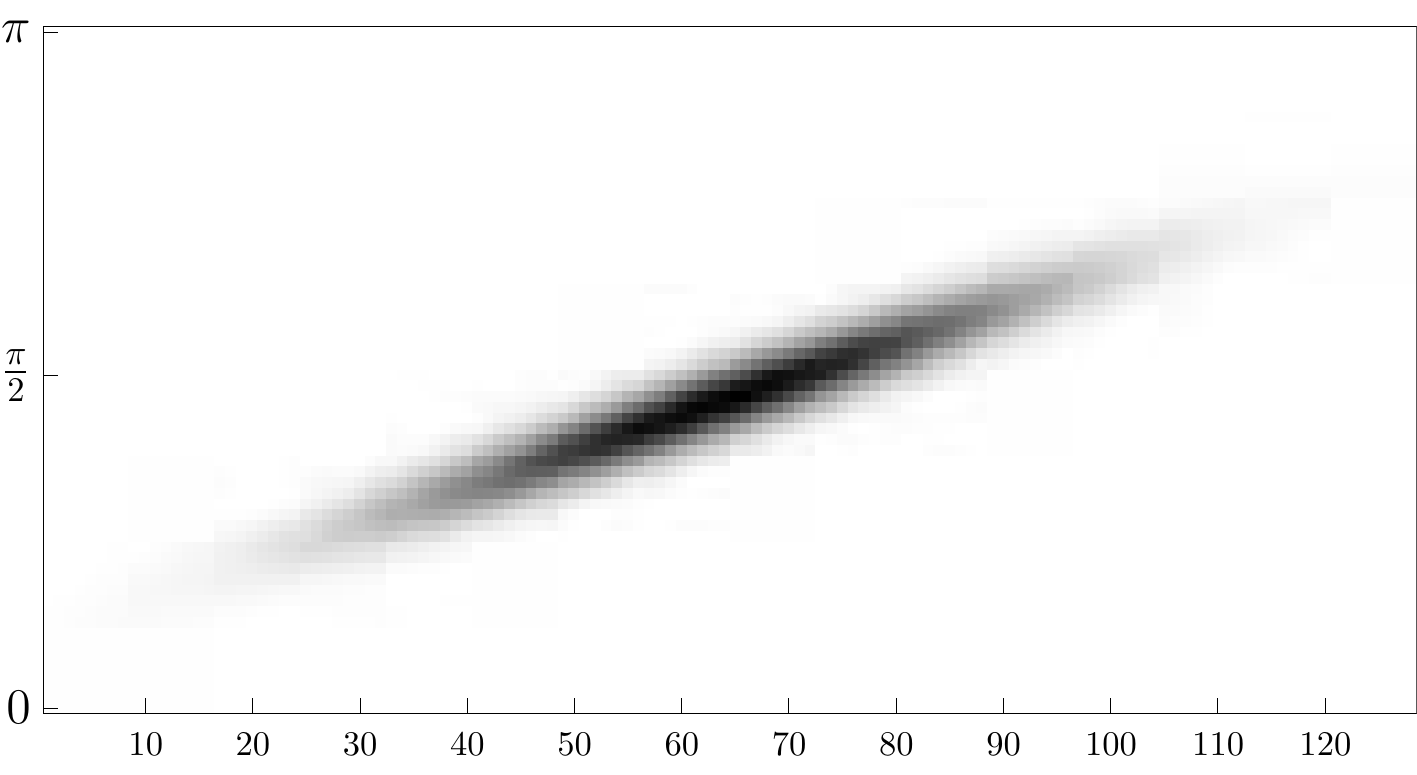}
\includegraphics[width=0.3\hsize]{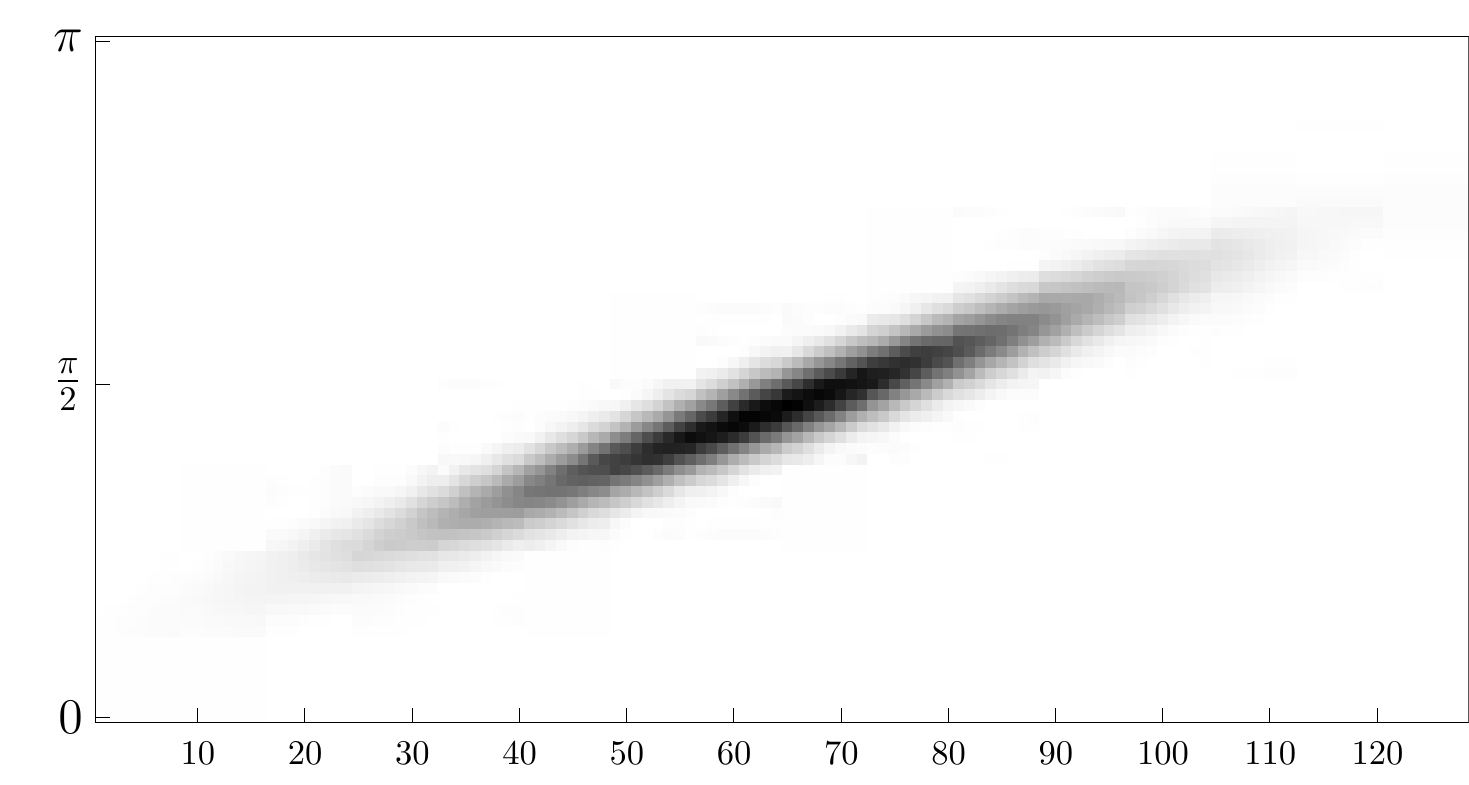}
}
%\caption{The reassigned Gabor transforms of the signals depicted in Fig. } \label{fig:7}
%\end{figure}
\caption{Absolute value of reassigned Gabor transforms of the signals depicted in the right of Fig. \ref{fig:4}.
}\label{fig:5}
\end{figure}

To evaluate these 4 methods we apply the reassignment scheme to the
reassignment of a linear chirp that is multiplied by a modulated Gaussian and
is sampled using $N =128$ samples. %The input signal is an analytic signal so it
%suffices to show its Gabor transform from 0 to $\pi$.
A visualization of this complex
valued signal can be found Fig. \ref{fig:4} (top). The other signals
in this figure are the reconstructions from the reassigned Gabor transforms that
are given in Fig. \ref{fig:5}. Here the topmost image shows the Gabor transform of
the original signal. One can also find the reconstructions and reassigned Gabor transforms
respectively using the four methods of reassignment. The
parameters involved in generating these figures are $N = 128$, $K = 128$, $M = 128$, $L=1$.
Furthermore $a = 1/6$ and the time step for
the PDE based method is set to $\Delta t =10^{-3}$. All images show a snapshot of the
reassignment method stopped at $t = 0.1$. The signals are scaled such
that their energy equals the energy of the input signal. This is needed to correct
for the numerical diffusion the discretization scheme suffers from.
Clearly the reassigned signals resemble the input signal quite well. The PDE scheme
that uses the sampled continuous window shows some defects. In contrast, the
PDE scheme that uses $\PPP^{D}_{a}$ resembles the modulus of the original signal the
most. Table \ref{table7} shows the relative
$\ell_{2}$-errors for all 4 experiments.
%, i.e.
%\[
%\begin{array}{ll}
%\epsilon_{1}= \frac{\|\ul{f}-\tilde{\ul{f}}\|_{\ell_{2}(I)}}{\|\ul{f}\|_{\ell_{2}(I)}}\ ,
% &
%\epsilon_{2}= \frac{\|\, |\ul{f}|-|\tilde{\ul{f}}|\, \|_{\ell_{2}(I)}}{\|\ul{f}\|_{\ell_{2}(I)}}\ .
%\end{array}
%\]
%where $\ul{f} \in \ell_{2}(I)$ is the original discrete signal and $\tilde{\ul{f}}$ the reconstructed signal after normalized re-assignment %in the Gabor domain.
%$\epsilon_1$ and the relative
%error of its modulus
%$\epsilon_2$. %, where the time step $\Delta t$ for the PDE
%scheme using $\PPP_{a}^{D}$ is set to $10^{-4}$ instead of $10^{-3}$.
Advantages of the erosion scheme (\ref{erosion}) over the PDE-scheme of Section \ref{ch:upwind} are:
\begin{enumerate}
\item
The erosion scheme does not
produce numerical approximation-errors in the phase, which is evident since the phase is not used in the computations.
\item
The erosion scheme does not involve numerical diffusion as it does not suffer from finite step-sizes.
\item
The separable erosion scheme is much faster. % from a computational point of view.
%\footnote{For example the erosions (\ref{erosion}) are separable, since we use a quadratic structure element (\ref{erosionkernel}) which is %a sum of a square in $p$ and a square in $q$, reducing the complexity of the erosion by order $O(K+M)$ in stead of $O(KM)$.}
\end{enumerate}
The convection time in the erosion scheme is different than the convection time in the upwind-scheme, due to violation of the Cauchy-Riemann equations.
%Typically, to get similar visual sharpening of the re-assigned Gabor transforms, the %convection time of the PDE-scheme should be taken larger than the convection time of the %erosion scheme (due to numerical blur in the PDE-scheme).
Using a discrete window (that satisfies the discrete Cauchy-Riemann relations) in the PDE-scheme is more accurate. Thereby one can obtain more visual sharpening in the Gabor domain while obtaining similar relative errors in the signal domain. For example $t=0.16$ for the PDE-scheme roughly corresponds to erosion schemes with $t=0.1$ in the sense that the $\ell_{2}$-errors nearly coincide, see Table~1.
%, but even here the erosion scheme seems to sharpen more than the PDE-scheme.
%thereforeF
%$\epsilon_{1}$ is smaller then one would expect based on inspection of
%Figure \ref{fig:6}.
\begin{table}[t]\label{table7}
\begin{center}
\begin{tabular}{l|c|c|c|}
  & $\epsilon_1$ & $\epsilon_2$ & $t$ \\ \hline Erosion continuous window & $2.41 \cdot 10^{-2}$ & $8.38 \cdot 10^{-3}$ & $0.1$ \\ Erosion discrete window & $8.25 \cdot 10^{-2}$ & $7.89 \cdot 10^{-2}$ & $0.1$ \\ PDE continuous window & $2.16 \cdot 10^{-2}$ & $2.21 \cdot 10^{-3}$ & $0.1$ \\ PDE discrete window  & $1.47 \cdot 10^{-2}$ & $3.32 \cdot 10^{-4}$ & $0.1$ \\ \ & $2.43 \cdot 10^{-2}$ & $6.43 \cdot 10^{-3}$ & $0.16$ \end{tabular}\caption{The first column shows {\small $\epsilon_{1}= (\|\, \ul{f}-\tilde{\ul{f}}\, \|_{\ell_{2}(I)})\|\ul{f}\|^{-1}_{\ell_{2}}$}, the relative error of the complex valued reconstructed signal compared to the input signal. In the second column {\small$\epsilon_{2}= (\|\, |\ul{f}|-|\tilde{\ul{f}}|\, \|_{\ell_{2}(I)})\|\ul{f}\|^{-1}_{\ell_{2}}$} can be found which represents the relative error of the modulus of the signals. Parameters involved are $K=M=N=128$, window scale $a=\frac{1}{8}$ and convection time $t=0.1$, with times step $\Delta t=10^{-3}$ if applicable. PDE stand for the upwind scheme presented in subsection~\ref{ch:upwind} and erosion means the morphological erosion method given by eq. (\ref{erosion}). \vspace{-0.3cm}\mbox{}} \end{center} \label{tab.reassignment.RelativeErrors}
\end{table}
The method that uses a
sampled version of the continuous window shows large errors and indeed
in Fig.~\!\ref{fig:5} the defects are visible. This shows the importance of the
window selection, i.e. in the PDE-schemes it is better to use window $\PPP_{a}^{D}$ rather than window $\PPP_{a}^{C}$. However, Fig.~\!\ref{fig:8} and Table~\!1 clearly indicate that in the erosion schemes it is better to choose window $\PPP_{a}^{C}$ than $\PPP_{a}^{D}$. %For further details see \cite[ch:5.7, 5.8]{Janssen2009b}.
%\begin{figure}
%\centerline{
%\includegraphics[width=0.75 \hsize]{recsignals}
%}
%\caption{
%The modulus of the signals corresponding to the Gabor transforms in Figure \ref{fig:5}. See the
%caption of that figure for further details. Note that ``$PDE$-sampled window'' stands for the PDE-approach of subsection \ref{ch:upwind} (with $M=1$) applied to discrete Gabor transforms constructed by kernels $\PPP^{C}_{a}$, recall Fig.\ref{fig:comparison}.
%The erosion scheme does not perform well in the signal domain if the discrete kernel $\PPP^{D}_{a}$ is used. The PDE-scheme performs better if $\PPP^{D}_{a}$ is used (rather than $\PPP^{C}_{a}$).
%}\label{fig:6}
%\end{figure}
% TABEL !!!!!!!
%
%\section*{Appendix}
%\addcontentsline{toc}{section}{Appendix}
%
%
%When placed at the end of a chapter or contribution (as opposed to at the end of the book), the numbering of tables, figures, and equations %in the appendix section continues on from that in the main text. Hence please \textit{do not} use the \verb|appendix| command when writing %an appendix at the end of your chapter or contribution. If there is only one the appendix is designated ``Appendix'', or ``Appendix 1'', or %``Appendix 2'', etc. if there is more than one.
\section{The Exact Analytic Solutions of Reassigned Gabor Transforms of 1D-chirp Signals \label{ch:chirp}}

In the previous section we have introduced several numerical algorithms for differential reassignment. We compared them experimentally on the special case where the initial 1D-signal (i.e. $d=1$) is a chirp. In this section we will derive the analytic solution of this reassignment in the Gabor domain. Furthermore, we will show the geometrical meaning of
the Cauchy-Riemann relation (\ref{CRMetric}) in this example.

\begin{lemma}\label{lemma:erode}
Let $t>0$, $\eta \in [\frac{1}{2},\infty)$.
Let $\Phi_{t}^{\eta}: \mathcal{H}_{n} \to C(H(3))$ be the operator that maps  a Gabor transform $\mathcal{W}_{\psi_{a}}f$ to its reassigned Gabor transform
\[
Q_{t,a,\eta}(|\mathcal{W}_{\psi_{a}}f|)\; e^{i \textrm{arg}\{\mathcal{W}_{\psi_{a}}f\}}
\]
where $Q_{t,a,\eta}$ maps the absolute value $|\mathcal{W}_{\psi_{a}}f|$ of the Gabor transform to
the unique viscosity solution $Q_{t,\eta}(|\mathcal{W}_{\psi_{a}}f|):=\tilde{W}(\cdot,t)$ of the Hamilton-Jacobi equation
\begin{equation}\label{PSHam}
\left\{
\begin{array}{l}
\frac{d \tilde{W}}{dt}(p,q,t)= -\frac{1}{2\eta}\left(a^{-2}\frac{d^2 \tilde{W}}{dp^2}(p,q,t) + a^2\frac{d^2 \tilde{W}}{dq^2}(p,q,t)\right)^{2\eta}\ , \qquad (p,q) \in \R^{2}, t>0, \\
\tilde{W}(p,q,0)=|\mathcal{W}_{\psi_{a}}f(p,q,s)|=|\mathcal{G}_{\psi}(p,q)|, \qquad (p,q) \in \R^{2},
\end{array}
\right.
\end{equation}
at fixed time $t>0$. Set $\Xi^{a,\alpha}_{t}:= \Phi_{t}^{\eta} \circ \mathcal{W}_{\psi_{a}}$. Let
$\mathcal{D}_{a}$ be the unitary scaling operator given by (\ref{dil}).
%\[
%(\mathcal{D}_{a}f)(x)= a^{-d/2}f(x/a)
%\]
Then
\begin{equation}\label{scaling}
(\Xi^{a,\eta}_{t} f)(p,q,s)=\sqrt{a}\, (\Xi^{1,\eta}_{t} \mathcal{D}_{a^{-1}}f)(a^{-1}p,aq,s)
\end{equation}
\end{lemma}
\textbf{Proof }First of all one has by substitution in integration that
\begin{equation} \label{eqt}
(\mathcal{W}_{\psi_{a}}f)(p,q,s)= \sqrt{a} (\mathcal{W}_{\psi_{a=1}} \mathcal{D}_{a^{-1}}f)(\frac{p}{a}, aq,s)\ .
\end{equation}
for all $f \in \mathbb{L}_{2}(\R^d)$, $(p,q,s) \in H(3)$.
So if we introduce the scaling operator $\mathbf{D}_{a}: \mathbb{L}_{2}(H(3)) \to \mathbb{L}_{2}(H(s))$ by
\[
(\mathbf{D}_{a}U)(p,q,s)=U(a^{-1}p, qa ,s)\ , a>0,
\]
then (\ref{eqt}) can be written as $\mathcal{W}_{\psi_{a}}= \sqrt{a} \mathbf{D}_{a}\mathcal{W}_{\psi_{1}}$ and
by the chain-law for differentiation one has $Q_{t,a,\eta}=\mathbf{D}_{a}Q_{t,1,\eta}\mathbf{D}_{a^{-1}}$ (where we note that the viscosity condition, cf. \cite{Evans} is scaling covariant). Consequently, we obtain
\[
\begin{array}{ll}
(\Xi^{a,\eta}_{t} f)(p,q,s) & =\sqrt{a}\, e^{i \textrm{arg}\{\mathcal{W}_{\psi_{a}}f(p,q,s)\}} \,(\mathbf{D}_{a} Q_{t,1,\eta} \mathbf{D}_{a^{-1}} \mathbf{D}_{a} \mathcal{W}_{\psi_{1}}\mathcal{D}_{a^{-1}}f)(p,q,s)\\ & = \sqrt{a} \, e^{i \textrm{arg}\{\mathcal{W}_{\psi_{1}}f(a^{-1}p,aq,s)\}}\, (Q_{t,1,\eta}(\mathcal{D}_{a^{-1}}f))(a^{-1}p,aq,s) \\ &=
\sqrt{a}\,(\Xi^{1,\eta}_{t} f)(a^{-1}p,aq,s)
 \hfill \Box
\end{array}
 \]
\begin{corollary}
By means of the scaling relation Eq. (\ref{scaling}) we may as well restrict ourselves to the case $a=1$ for analytic solutions.
\end{corollary}
\begin{lemma}\label{lemma:comp}
Let $g(\xi)=e^{w \xi^2 + u \xi +t}$, with $\textrm{Re}(w)<0$ and $\textrm{Im}(w)>0$. Then
\[
\int_{\R} g(\xi)\, {\rm d}\xi=\sqrt{\frac{\pi}{-w}} e^{\frac{4tw-u^2}{4w}}
\]
where the complex square root is taken using the usual branch-cut along the positive real axis.
\end{lemma}
\textbf{Proof }
By contour integration where the contour encloses a the two-sided section in the complex plane given by the intersection of the ball with radius $R$ with $0<\arg(z)<-\frac{1}{2}\arg(-w)$ (positively) and $-\pi <\arg(z)<-\pi + -\frac{1}{2}\arg(-w)$ (negatively) then one has by Cauchy's formula of integration and letting $R \to \infty$ that
\[
\int_{\R} g(\xi)\, {\rm d}\xi= \int_{\arg(z)= -\frac{1}{2}\arg(-w)} e^{|w|e^{i \arg(w)}z^2-\frac{u^2}{4w}}\, {\rm d}z=
\sqrt{\frac{|w|}{-w}} \sqrt{\frac{\pi}{|w|}} e^{-\frac{u^2}{4w}}
\]
from which the result follows. $\hfill \Box$
\begin{theorem} \label{th:exact1}
Let $r,b>0$. Let $f$ be (a chirp signal) given by
\begin{equation}\label{chirp}
f(\xi)= e^{-\frac{\xi^2}{2 b^2}}e^{i\pi r \xi^2}\ ,
\end{equation}
Then its Gabor transform (with $\psi(\xi)=e^{-\pi \xi^2}$) equals
\begin{equation} \label{GTC}
\mathcal{W}_{\psi}f(p,q,s)=\sqrt{\frac{1}{2\pi b^2-ir+1}}\, e^{-2\pi i (s+\frac{pq}{2})} e^{(p,q) B_{r,b} (p,q)^{T}}
\end{equation}
with $B_{rb}=\textrm{Re}(B_{r,b})+i \, \textrm{Im}(B_{r,b})$, $\textrm{Re}(B_{r,b})<0$, $(\textrm{Im}(B_{r,b}))^{T}=\textrm{Im}(B_{r,b})$, $(\textrm{Re}(B_{r,b}))^{T}=\textrm{Re}(B_{r,b})$  given by
\[
\begin{array}{ll}
B_{rb} &= \frac{1}{r^2 b^4 +(\frac{1}{2\pi}+b^2)^2}
\left(
\left(
\begin{array}{cc}
-\frac{1}{2}(b^2 +\frac{1}{2\pi})-\pi r^2 b^4 & \pi r b^4 \\
\pi r b^4 & -\pi b^2(b^2 +\frac{1}{2\pi})
\end{array}
\right)  \right. \\ & \left.
+ \, i \,
\left(
\begin{array}{cc}
\pi r b^4 & \frac{1}{4\pi}+\frac{b^2}{2}+\pi r^2 b^4 \\
\frac{1}{4\pi} +\frac{b^2}{2}+\pi r^2 b^4 & -\pi r b^4
\end{array}
\right)
\right)
\end{array}
\]
and thereby we have for $\eta \in (\frac{1}{2},1]$:
\begin{equation}\label{result}
\begin{array}{l}
(\Xi^{a,\eta}_{t} f)(p,q,s)=\sqrt{a}\, (\Xi^{1,\eta}_{t} \mathcal{D}_{a^{-1}}f)(a^{-1}p,aq,s)\ , a>0, \textrm{ with } \\
(\Xi^{1,\eta}_{t}f)(p,q,s)= e^{2\pi i(s-\frac{pq}{2})} e^{i \, \textrm{Im}(B_{rb})} (k_{t}^{\eta} \ominus
e^{(\cdot,\cdot) \textrm{Re}(B_{rb}) (\cdot,\cdot)^{T}})(p,q)\ , \\
\end{array}
\end{equation}
with (positive) erosion kernel at time $t>0$ given by
\begin{equation} \label{kerneleta}
k_{t}^{\eta}(p,q)= \frac{2\eta-1}{2\eta} t^{-\frac{1}{2\eta-1}} (p^2+q^2)^{\frac{\eta}{2\eta-1}}
\end{equation}
for $\eta \in (\frac{1}{2},1]$
and
\begin{equation} \label{kernel}
k_{t}^{\eta=\frac{1}{2}}(p,q)=
\left\{
\begin{array}{ll}
\infty & \textrm{if }p^{2}+q^2 \geq t^2 \\
0 & \textrm{if } p^{2}+q^2 < t^2
\end{array}
\right.
\end{equation}
\end{theorem}
\textbf{Proof }
Consider Eq. (\ref{GT_matrix_coeff}) where we set $n=1$, $d=1$, Eq.~(\ref{chirp}) and $\psi(\xi)=e^{-\pi \xi^2}$. We apply
Lemma~\ref{lemma:comp} with $w=i\pi r -\frac{1}{2b^2}-\pi$, $u=-2\pi p -2\pi i q$, $t=2 \pi i pq-\pi p^2$ and Eq.~(\ref{GTC}) follows. Then we note that when transporting along equiphase planes, phase-covariance is the same as phase invariance and indeed by Lemma \ref{lemma:erode} the main result (\ref{result}) follows. Finally, we note that
the viscosity solutions of the erosion PDE is given by morphological convolution with the erosion kernel which by the Hopf-Lax formula equals $k_{t}^{\eta}(p,q)= t \mathcal{L}_{\eta}(t^{-1}(p,q)^{T})$ where the Lagrangian $\mathcal{L}_{\eta}(p,q)= (\gothic{F}H_{\eta})(p,q)$ is obtained, cf. \cite[ch:3.2.2]{Evans}, by the Fenchel-transform of the hamiltonian $H_{\eta}(a^{-1}u,a v)=\frac{1}{2\eta}(a^{-2}(u)^2+a^2(v)^2)^{2\eta}$ that appears in the righthand side of the Hamilton-Jacobi equation Eq. (\ref{PSHam}), cf.~\cite[ch:2,p.24]{Rund}, Eq.~(\ref{erosion}), from which the result follows. $\hfill \Box$
\begin{remark}
In theorem we have set centered the chirp signal in position, frequency and phase. The general case follows by:
%\[
$(\mathcal{W}_{\psi} \mathcal{U}_{(p_0,q_0,s_{0})}f)(p,q,s)=(\mathcal{W}_{\psi}f)((p_{0},q_{0},s_{0})^{-1}(p,q,s))$.
%\]
\end{remark}
\begin{remark}
Note that $\textrm{trace}\{\textrm{Re}(B_{rb})\}<0$ and $\textrm{Det}\{\textrm{Re}(B_{rb})\}>0$, so both eigenvalues
of the symmetric matrix $\textrm{Re}(B_{rb})$ are negative. Consequently, the equi-contours of the spectrogram $(p,q) \mapsto |G_{\psi}f|(p,q)$ are ellipses. In a chirp signal (without window, i.e. $b\to \infty$) frequency increases linear with $\xi$ via rate $r$ and thereby one expects the least amount of decay in the spectrogram along $q=r\, p$, i.e. one expects $(1,r)^{T}$ to be the eigenvector with smallest absolute eigenvalue of $\textrm{Re}(B_{r,b})$. This is indeed only the case if $b\to \infty$ as
\[
\textrm{Re}(B_{r,b})
\left(
\begin{array}{c}
1 \\
r
\end{array}
\right)= \frac{(b^2/2)}{r^{2}b^4+ ((2\pi)^{-1}+b^2)^2}
\left( \left(
\begin{array}{c}
1 \\
r
\end{array}
\right)+
\left(
\begin{array}{c}
\frac{1}{2\pi b^2} \\
0
\end{array}
\right)
\right).
\]
%If we set $(p,q)$
%For suitable choice of parameters $(a,b)$ this indeed holds
%since $|r-v^{2}|=\left| \frac{1+2\beta}{\frac{r}{2}-r^{-1}\beta+\sqrt{1+((r/2)+r^{-1}\beta)^2}}\right|$, where
%$\beta= \frac{1}{8\pi^2 r^{2}b^4} +\frac{1}{8\pi^2 a^2 b^2} - \frac{a^2}{2b^2}-\frac{1}{2}$ and $\beta=-\frac{1}{2}$
%if $a=\frac{\sqrt{1+\sqrt{1+16 \pi^2 b^4}}}{2 \sqrt{2} \pi b}$.
\end{remark}
\begin{remark}
Note that $\textrm{Det}\{\textrm{Re}(B_{rb})\}<0$, so the eigenvalues
of the symmetric matrix $\textrm{Re}(B_{rb})$ have different sign and consequently the equiphase lines in phase space (where we have set $s=-pq/2$) are hyperbolic.
\end{remark}
In case $\eta \to \infty$ the Lagrangian is homogeneous, the Hamilton-Jacobi equation (that describes the evolution of geodesically equidistant surfaces $\{(p,q) \in \R^2\; |\; \tilde{W}(p,q,t)=c\}$) becomes time-independent, cf.~\cite[ch:4,p.170]{Rund}
\[
1= a^{-2}\frac{d^{2}}{dp^2}\tilde{W}(p,q,t) + a^{2}\frac{d^{2}}{dp^2}\tilde{W}(p,q,t).
\]
In case $\eta=\frac{1}{2}$ the Hamiltonian is homogeneous and the erosion kernels are flat accordingly. The non-flatness of the erosion kernels can be controlled via parameter $\eta \in (\frac{1}{2},\infty)$.
Furthermore, the Hamilton-Jacobi
equation in (\ref{PSHam}) is invariant under monotonic transformations iff $\eta=\frac{1}{2}$ and this allows us to compute $(\Xi^{1,\frac{1}{2}}_{t}f)(p,q,s)$.
\begin{lemma}\label{lemma:3}
Let $\lambda_{1}, \lambda_{2}$ denote the eigenvalues (with $|\lambda_{1}|<|\lambda_{2}|$) of $\textrm{Re}(B_{r,b})$ with respective normalized eigenvectors $\ul{k}_{1}$ and $\ul{k}_{2}$.
Then each vector in $\R^{2}$ can be written as
$(p,q)=\alpha_{1}(p,q)\ul{k}_{1}+\alpha_{2}(p,q)\ul{k}_{2}$. One has
{\small
\[
\begin{array}{l}
\lambda_{1}= \frac{1}{d_{br}} \left( -1 -4 b^2 \pi(1+ \pi(b^2+r^2)) + \sqrt{(\pi r b^4)^2+\frac{1}{4}\left(\frac{1}{4\pi} +b^2 \pi(r^2-b^2)\right)^2} \right) \ , \\
\lambda_{2}= \frac{1}{d_{br}} \left( -1 -4 b^2 \pi(1+ \pi(b^2+r^2)) - \sqrt{(\pi r b^4)^2+\frac{1}{4}\left(\frac{1}{4\pi} +b^2 \pi(r^2-b^2)
\right)^2}\right)\ , \\
\alpha_{1}(p,q)=N_{br}^{1}\; d_{br} \; \left( \frac{4b^4\pi^2r}{c_{br}}p+  \left(\frac{1}{2}-\left(\frac{-1+4b^2 \pi(b^2-r^2)}{c_{br}}\right)\right)q\right)\ ,\\
\alpha_{2}(p,q)=N_{br}^{2}\; d_{br}\left( -\frac{4b^4\pi^2r}{c_{br}}p+  \left(\frac{1}{2}+\left(\frac{-1+4b^2 \pi(b^2-r^2)}{c_{br}}\right)\right)q\right)\ ,  \\
\end{array}
\]
}
with $d_{br}=r^2b^4+(b^2+\frac{1}{2\pi})^2$ and
\[
\begin{array}{l}
c_{br}=\sqrt{1+ 8b^2\pi^2(r^2+b^2(-1+2\pi^2(b^4+2b^2(2b^2-1)r^2+r^4)))}\ , \\
N_{br}^{k}=\sqrt{1+\frac{|1-4b^4\pi^2+4b^2\pi^2r^2(-1)^{k}c_{br}|^2}{64 \pi^4 b^8r^2}} \textrm{ for }k=1,2.
\end{array}
\]
\end{lemma}
We omit the proof as the result follows by direct computation.
\begin{lemma}
Consider the ellipsoidal isolines given by
\[
\{(p',q')\in \R^{2}\;|\;
(p',q') B_{r,b} (p',q')^{T}=\lambda_{1}|\alpha_{1}(p',q')|^2+\lambda_{2}|\alpha_{2}(p',q')|^2=c\},
 \]
 with $c>0$ arbitrary. Consider the line
spanned by the principal direction with smallest absolute eigenvalue. Applying the settings of Lemma \ref{lemma:3} this line is given by $\alpha_{2}(p',q')=0$. Consider another point $(p,q)$ on this line, i.e. $\alpha_{2}(p,q)=0 \desda p=\tau q$, with $p>0$. Then there exists a unique $p^{*}=p^{*}(p,p\tau)=\frac{p}{1-\frac{\sigma}{1+\tau^2}}$
such that the circle $(p-p^*)^{2}+(\tau p-\tau p^{*})^{2}=t^2$
is tangent to one of the ellipses $(p',q') B_{r,b} (p',q')^{T}=\lambda_{1}|\alpha_{1}(p',q')|^2+\lambda_{1}|\alpha_{2}(p',q')|^2=c$ in $(p^{*},\tau p^{*})$.
The corresponding time equals
\[
t_{max}(p,p\tau)= \sqrt{1+\tau^2}\, \sqrt{(p-p^{*}(p,q))^2}= \frac{\sqrt{\sigma}|p|}{\left|1-\sqrt{\frac{\sigma}{1+\tau^2}}\right|}.
\]
where the constants $\sigma$ and $\tau$ are given by
{\small
\[
\begin{array}{l}
\sigma= b^4 \left(4\frac{N_{br}^{2}}{N_{br}^{1}}\right)^2 \pi^2 r \, \frac{c_{br} d_{br}(-1 + c_{br} + 4 b^2 \pi^2 (b - r) (b + r)) (1 + c_{br} +
   4 b^2 \pi (1 + \pi (b^2 + r^2)))}{((1 + c_{br} - 4 b^4 \pi ^2)^2 +
   8 b^2 \pi^2 (1 + c_{br} + 4 b^4 (-1 + 2 b^2) \pi^2) r^2 +
   16 b^4 \pi^4 r^4)^{\frac{3}{2}
 } (1 - c_{br} + 4 b^2 \pi (1 + \pi (b^2 + r^2)))},\\
\tau=\frac{4b^4\pi^2r}{\frac{c_{br}}{2}-1+4b^2 \pi(b^2-r^2)}.
\end{array}
\]
}
\end{lemma}
\textbf{Proof }
In order to have the circle touching the ellipse in $(p^{*},\tau p^{*})$, we have to solve the following system
\[
\begin{array}{l}
%\alpha^{1}(p',q')'=
\alpha_{2}(p',q')=0, \\
(p-p')^2+(q-q')^2 =t^2, \\
\kappa^{2}(p',q')= \frac{1}{t^2}.
\end{array}
\]
The curvature along isocontours of $u(p,q)=(p',q') B_{r,b} (p',q')^{T}$ is expressed as
{\small
\begin{equation} \label{draak}
\begin{array}{l}
\kappa(p',q')= -
\frac{1}{\left(\left(\frac{\partial u(p',q')}{\partial p'}\right)^2+\left(\frac{\partial u(p',q')}{\partial q'}\right)^2\right)^{\frac{3}{2}}} \cdot \\\left(\frac{\partial u(p',q')}{\partial q'} \frac{\partial^{2} u(p',q')}{\partial p'^2}-2 \frac{\partial u(p',q')}{\partial q'}\frac{\partial u(p',q')}{\partial p'} \frac{\partial^{2} u(p',q')}{\partial p' \partial q'}+ \left(\frac{\partial u(p',q')}{\partial p'}\right)^2 \frac{\partial^{2} u(p',q')}{\partial q'^2}\right)= \sigma^{-1/2} (p')^{-1}
\end{array}
\end{equation}
}
then along $\ul{k}_{1}$ we have $\alpha_{2}(p',q')=0$ i.e. $q'=\tau p'$ and substitution in (\ref{draak}) yields
\[
\kappa^{-2}(p',\tau p')=  (p')^2 \, \sigma=t^2=(p-p')^2+(q-q')^2=(p-p')^{2}(1+ \tau^2)
\]
this yields $p'(p,\tau p)=\frac{p}{1\pm\sqrt{\frac{\sigma}{1+\tau^2}}}$, so since $p^{*} \geq p$ we find $p^{*}(p,p\tau)=\frac{p}{1-\sqrt{\frac{\sigma}{1+\tau^2}}}$.$\hfill \Box$ \\
%\end{proof}
The next theorem provides the exact solution $\Xi^{1,\frac{1}{2}}_{t}f : H_{r} \to \mathbb{C}$ of an eroded Gabor transform of a
chirp signal at time $t>0$. In the subsequent corollary we will show that the isocontours of the spectogram are non-ellipsoidal Jordan curves that shrink towards the principal eigenvector of $\textrm{Re}(B_{r,b})$ where they collaps as $t$ increases.
This collapsing behavior can also be observed in the bottom rows of Figure \ref{Fig:Exact} and \ref{Fig:Exact2}.
\begin{theorem}\label{th:exact}
Let $\eta=\frac{1}{2}$ and let $t>0$. Let $f$ be a chirp signal given by (\ref{chirp}) then the reassigned Gabor transform of $f$ (with window scale $a=1$) is given by
\begin{equation} \label{Xisol}
\begin{array}{l}
\!
(\Xi^{1,\frac{1}{2}}_{t}f)(p,q,s) =
\sqrt{\frac{1}{2\pi b^2-ir+1}}\, e^{-2\pi i (s+\frac{pq}{2})} e^{i \, (p,q)\textrm{Im}(B_{rb})(p,q)^{T}}
 (k_{t}^{\eta=\frac{1}{2}} \ominus e^{(\cdot,\cdot) \textrm{Re}(B_{rb})(\cdot,\cdot)^{T} })(p,q) \\[8pt]
 =\sqrt{\frac{1}{2\pi b^2-ir+1}}\, e^{-2\pi i (s+\frac{pq}{2})} \cdot \\
   \left\{
 \begin{array}{ll}
 e^{(\lambda_{t}(p,q))^2 \left(
 \frac{|\alpha_{1}(p,q)|^2\lambda_{1}}{(\lambda_{1}-\lambda_{t}(p,q))^2}+ \frac{|\alpha_{2}(p,q)|^2\lambda_{2}}{(\lambda_{2}-\lambda_{t}(p,q))^2}
 \right)} & \textrm{ if }\alpha_{1}(p,q) \alpha_{2}(p,q) \neq 0, \\
e^{(\alpha_{2}(p,q)+t)^{2} \lambda_{2}}          & \textrm{ if }\alpha_{1}(p,q) = 0, \\
e^{(\alpha_{1}(p,q)+t)^{2} \lambda_{1}}          & \textrm{ if }\alpha_{2}(p,q) = 0 \textrm{ and }t \leq t_{max}(p,\tau p). \\
 \end{array}
 \right.
\end{array}
\end{equation}
where $\lambda_{1}, \lambda_{2}$ denote the negative eigenvalues (with $|\lambda_{1}|<|\lambda_{2}|$) of $\textrm{Re}(B_{r,b})$ and where
$(p,q)=\alpha_{1}(p,q)\ul{k}_{1}+\alpha_{2}(p,q)\ul{k}_{2}$ with $\{\ul{k}_{1},\ul{k}_{2}\}$ the normalized eigenvectors of $\textrm{Re}(B_{r,b})$ and where the Euler-Lagrange multiplier $\lambda_{t}:=\lambda_{t}(p,q)$ corresponds to the unique zero $P_{t}(\lambda_t)=0$ of the $4$-th order polynomial
\begin{equation} \label{4}
P_{t}(\lambda):=t^2\prod \limits_{i=1}^{2}(\lambda-\lambda_{i})^2  - \sum \limits_{i=1}^{2}(\lambda_{i} \alpha_{i}(p,q))^{2}(\lambda-\lambda_{-i+3})^2,
\end{equation}
such that, for $\alpha_{1}(p,q)\alpha_{2}(p,q) \neq 0$, the corresponding unique solution $(p',q')$ of
\begin{equation} \label{pqaccent}
(\textrm{Re}(B_{r,b})-\lambda_{t}(p,q))\left(
\begin{array}{c}
p' \\
q'
\end{array}
\right)=
-\lambda_{t}(p,q)\left(
\begin{array}{c}
p \\
q
\end{array}
\right)
\end{equation}
has maximum $\sum \limits_{k=1}^{2}\lambda_{k}|\alpha_{k}(p',q')|^{2}$. The polynomial has at least 2 real-valued zeros. For each $(p,q) \in \R^{2}$ there exists a $t=t_{\textrm{max}}(p,q)\geq 0$ such that $P_{t}$ has 3 real-valued zeros and $P_{s}$ has 4 real-valued zeros for all $s>t$. Finally, the Lagrange multiplier satisfies the following scaling property
\begin{equation}\label{scaling}
\lambda_{t}(p,q)=\lambda_{\zeta t}(\zeta p,\zeta q) \textrm{ for all }\zeta>0.
\end{equation}
\end{theorem}
\textbf{Proof }
In case $\eta=\frac{1}{2}$ we have
\begin{equation}\label{optim}
\begin{array}{ll}
(k_{t}^{\frac{1}{2}} \ominus f)(p,q) &= \min \limits_{(p-p')^2+(q-q')^2 < t^2} e^{(p', q')\textrm{Re}(B_{r,b}) (p',q')^T} %\\ &
= e^{\left(\min \limits_{(p-p')^2+(q-q')^2 < t^2} (p', q')\textrm{Re}(B_{r,b}) (p',q')^T \right)}\ .
\end{array}
\end{equation}
Now $\textrm{Re}(B_{r,b})<0$, hence the minimum of the associated continuous function
can be found on the boundary of the convex domain. Application of Euler-Lagrange yields the system
\begin{equation} \label{ELa}
\begin{array}{l}
\textrm{Re}(B_{r,b})\left(
\begin{array}{c}
p' \\
q'
\end{array}
\right)=
\lambda_{t}(p,q)\left(
\begin{array}{c}
p'-p \\
q'-q
\end{array}
\right) \ ,\\
(p-p')^2+(q-q')^2=t^2.
\end{array}
\end{equation}
First we must find the Euler-Lagrange multipliers $\lambda_{t}(p,q)$. If $\lambda_{t}:=\lambda_{t}(p,q)$ is not an eigenvalue of $\textrm{Re}(B_{r,b})$ the resolvent $(\textrm{Re}(B_{r,b})-\lambda_{t} I)^{-1}$ exists and one finds
\[
t^{2}=\|\textrm{Re}(B_{r,b})(\textrm{Re}(B_{r,b})-\lambda_{t} I)^{-1} (p,q)^{T}\|^2,
\]
which yields the 4th-order polynomial equation (\ref{4}), where we note that
$(p,q)=\sum_{k=1}^{2}\alpha_{k}(p,q)\ul{k}_{k}$ with $(\ul{k}_{i},\ul{k}_{j})=\delta_{ij}$ so that
\begin{equation} \label{helper}
t^2=\|\textrm{Re}(B_{r,b})(\textrm{Re}(B_{r,b})-\lambda_{t} I)^{-1} (p,q)^{T}\|^2= \sum \limits_{k=1}^{2}\left|\frac{\lambda_{k}\alpha_{k}(p,q)}{\lambda_{k}-\lambda_{t}(p,q)}\right|^2.
\end{equation}
Now that the (four, three or two) Lagrange multipliers are known, we select the one minimal associated to the
global minimum and
the general solution is given by
\begin{equation} \label{form}
e^{(\lambda_{t}(p,q))^2 (p,q) (\textrm{Re}(B_{r,b})-\lambda_{t} I)^{-1}\textrm{Re}(B_{r,b})(\textrm{Re}(B_{r,b})-\lambda_{t} I)^{-1}},
\end{equation}
and substitution of $(p,q)=\sum_{i=1}^{2}\alpha_{i}(p,q)\ul{k}_{i}$ straightforwardly yields the first case in (\ref{Xisol}). The scaling relation (\ref{scaling}) now directly follows by (\ref{form}) and $\alpha_{i}(\zeta p,\zeta q)=\zeta \alpha_{i}(p,q)$, for all $i=1,2$, $(p,q)\in \R^{2}$ and $\zeta>0$.

\begin{figure}
\centerline{\includegraphics[width=0.15\hsize]{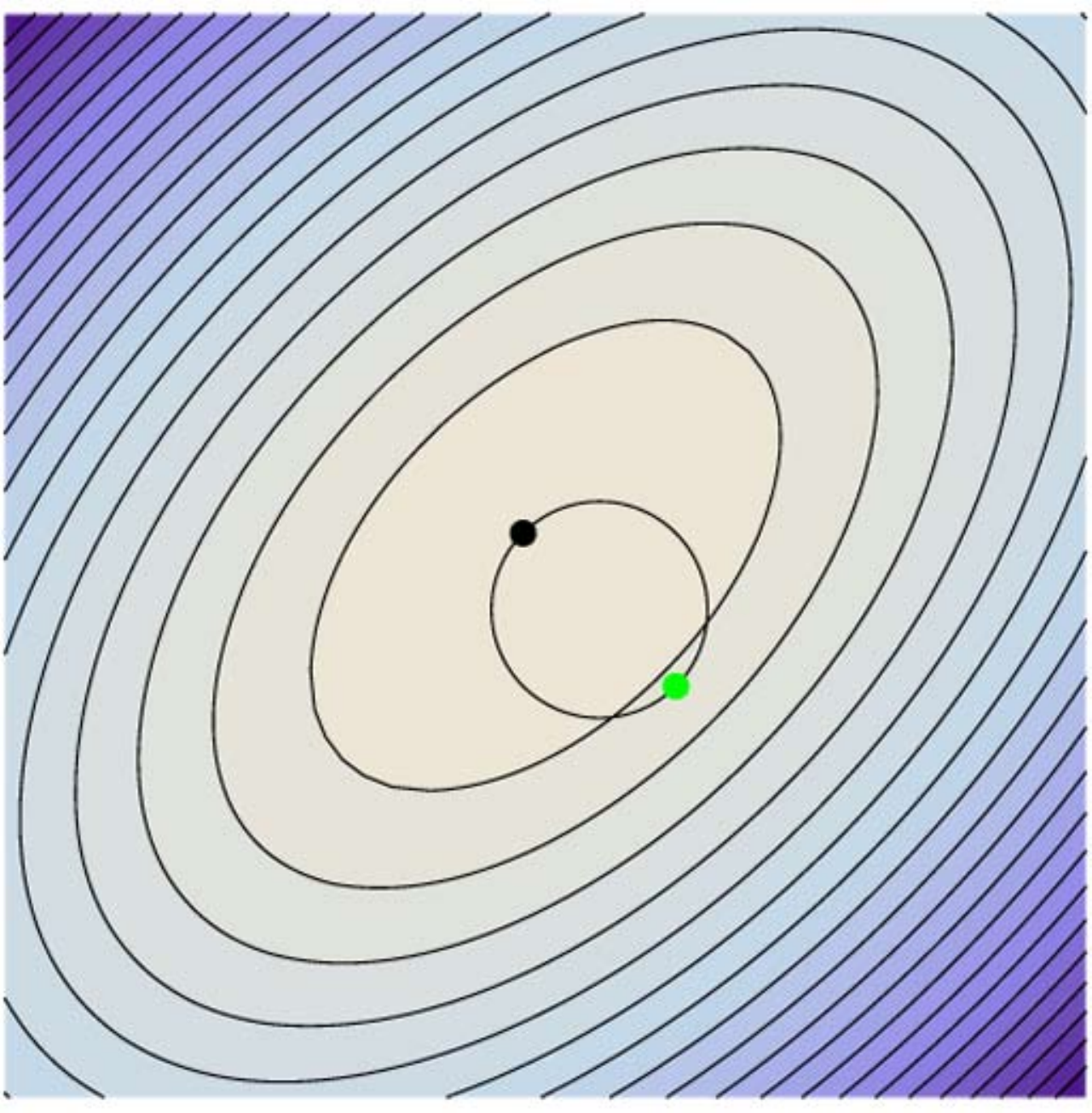}
\includegraphics[width=0.15\hsize]{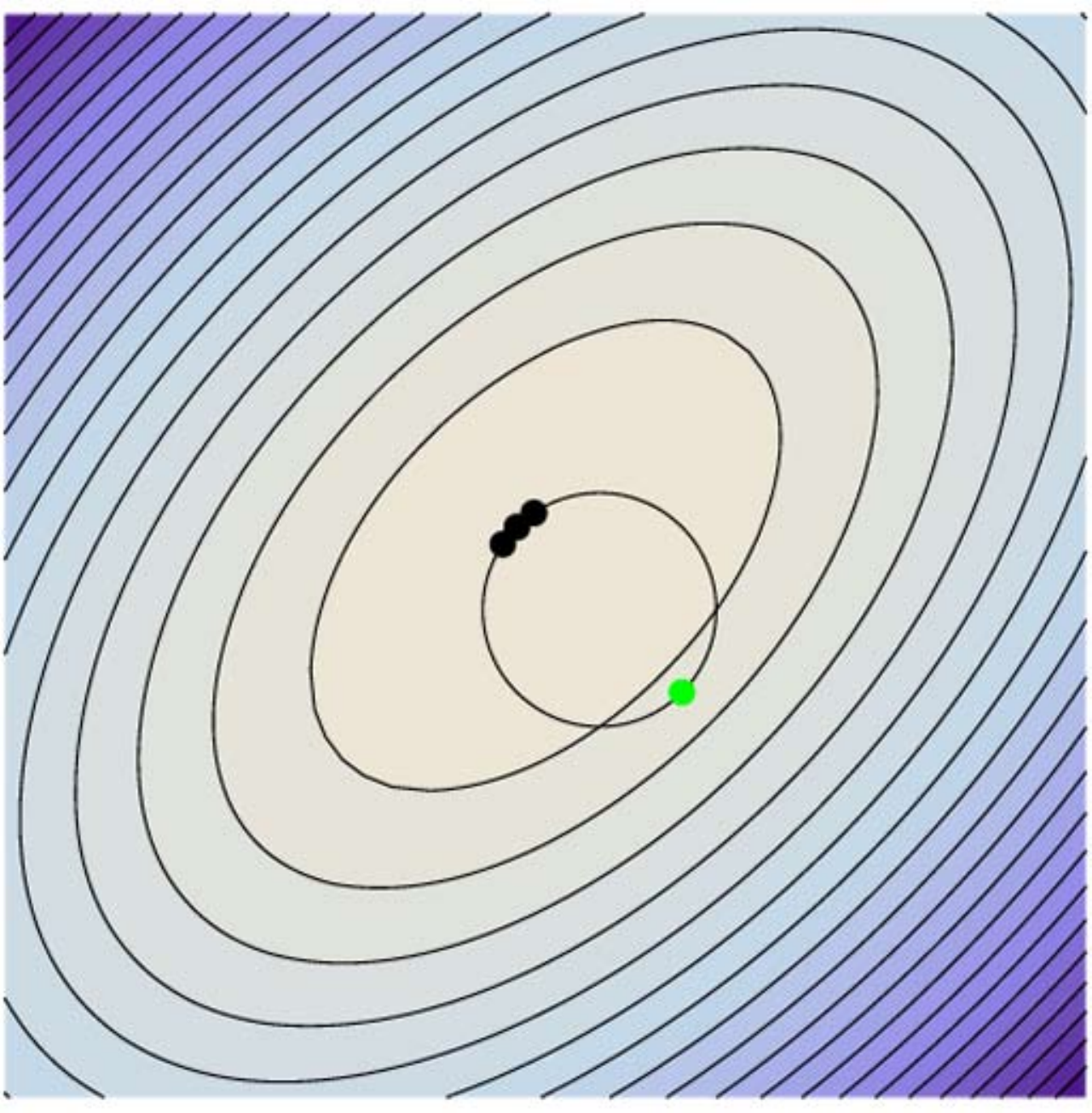}
\includegraphics[width=0.15\hsize]{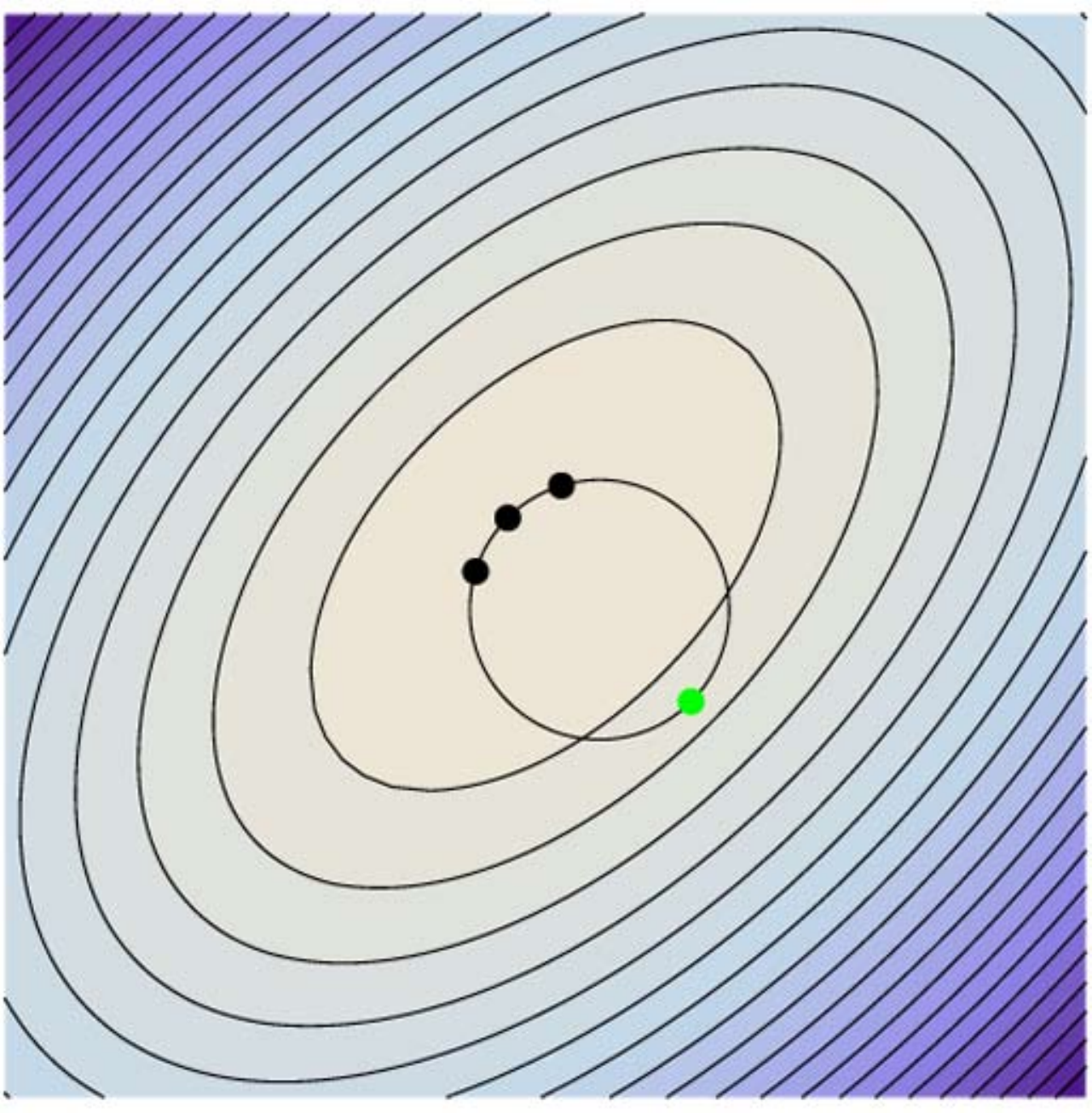}
\includegraphics[width=0.15\hsize]{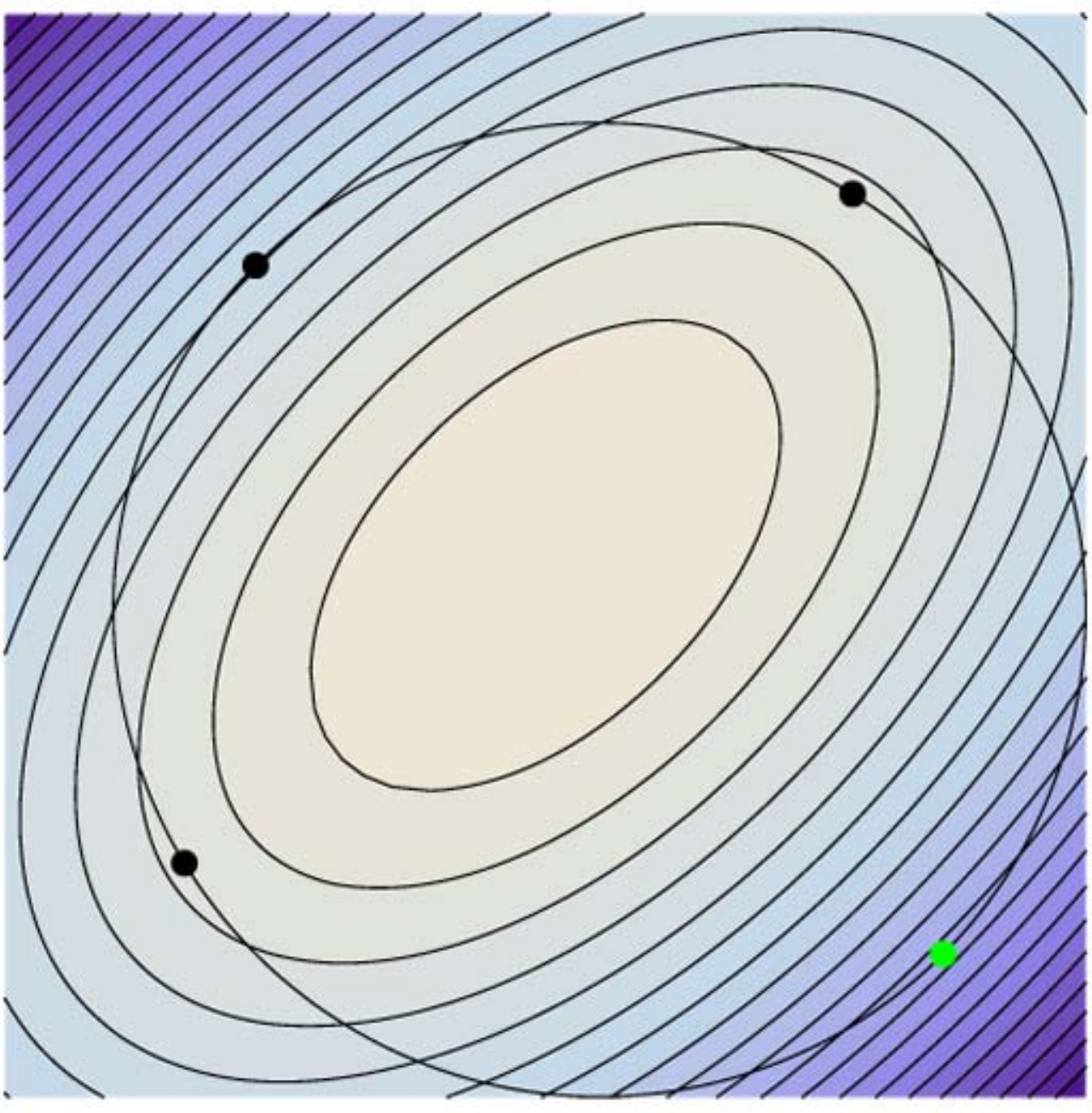}}
\centerline{\includegraphics[width=0.15\hsize]{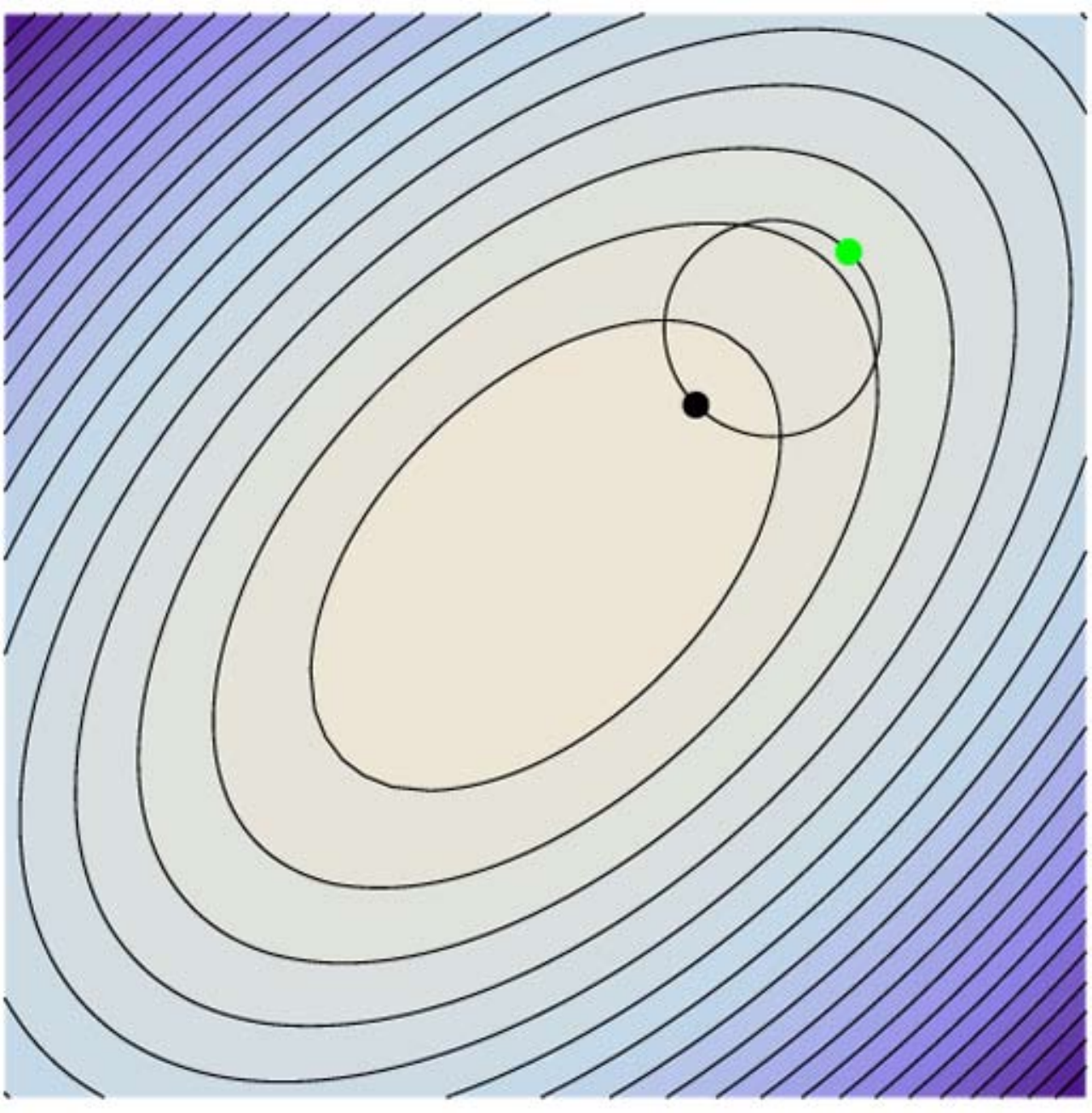}
\includegraphics[width=0.15\hsize]{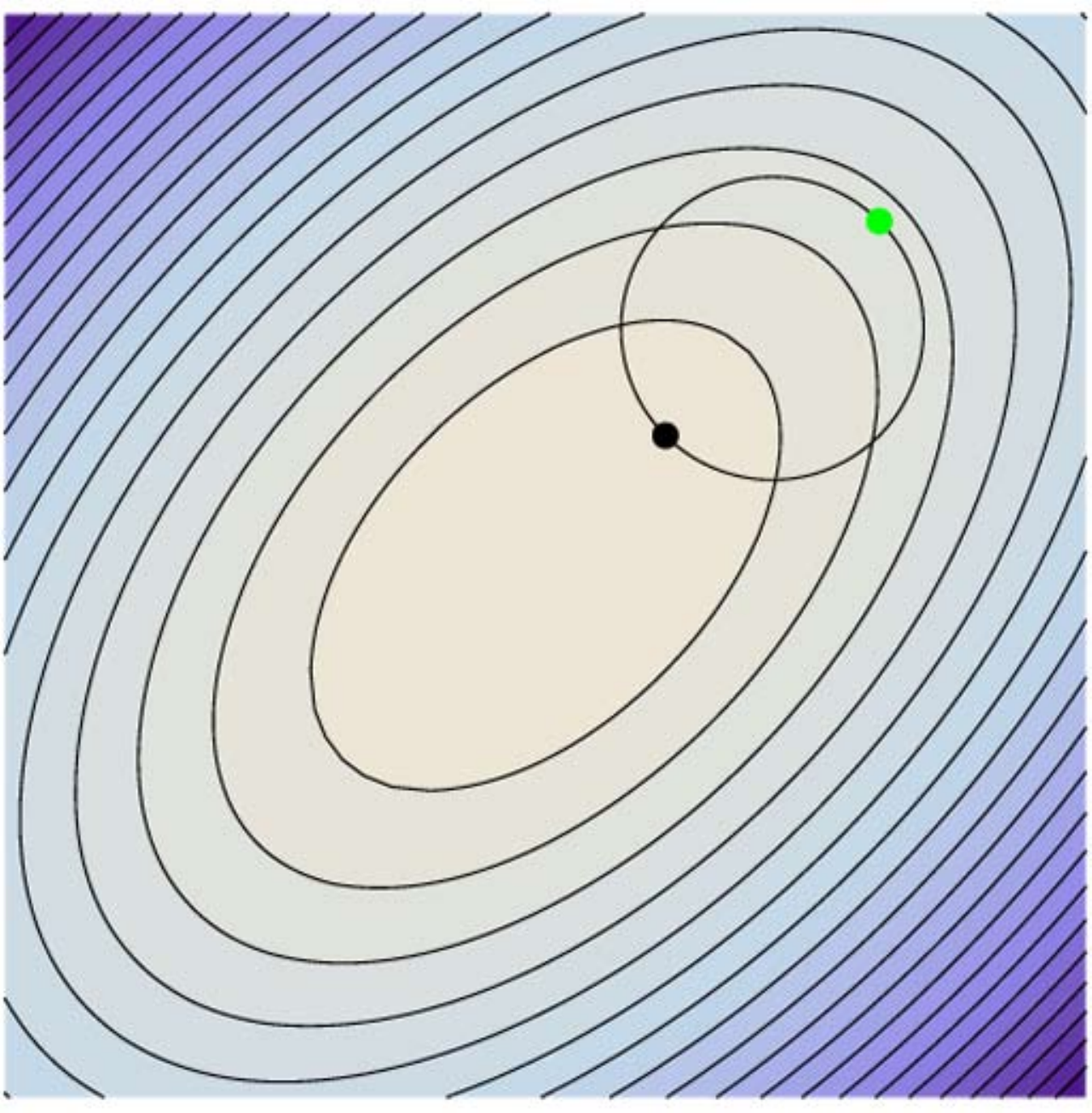}
\includegraphics[width=0.15\hsize]{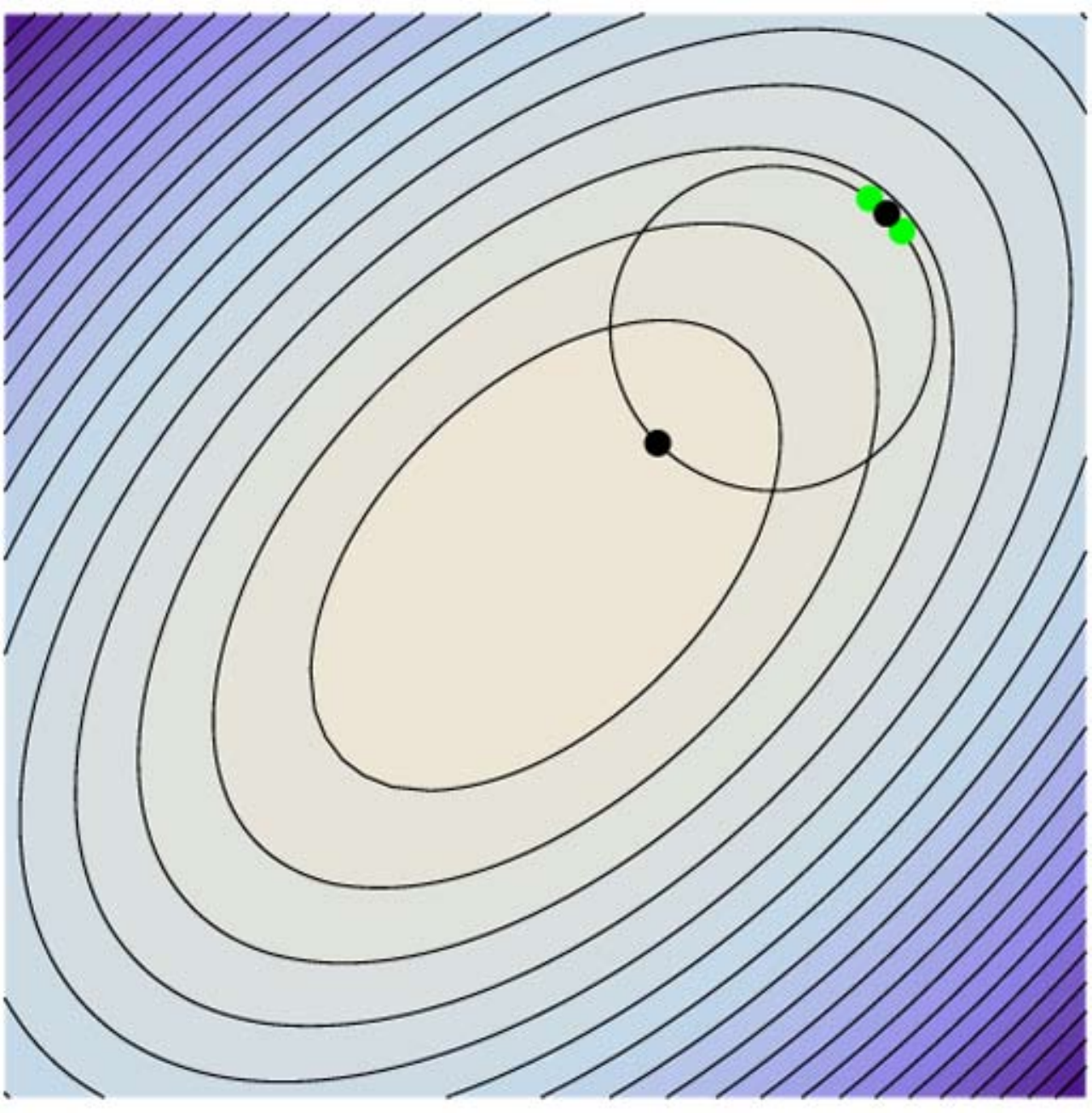}
\includegraphics[width=0.15\hsize]{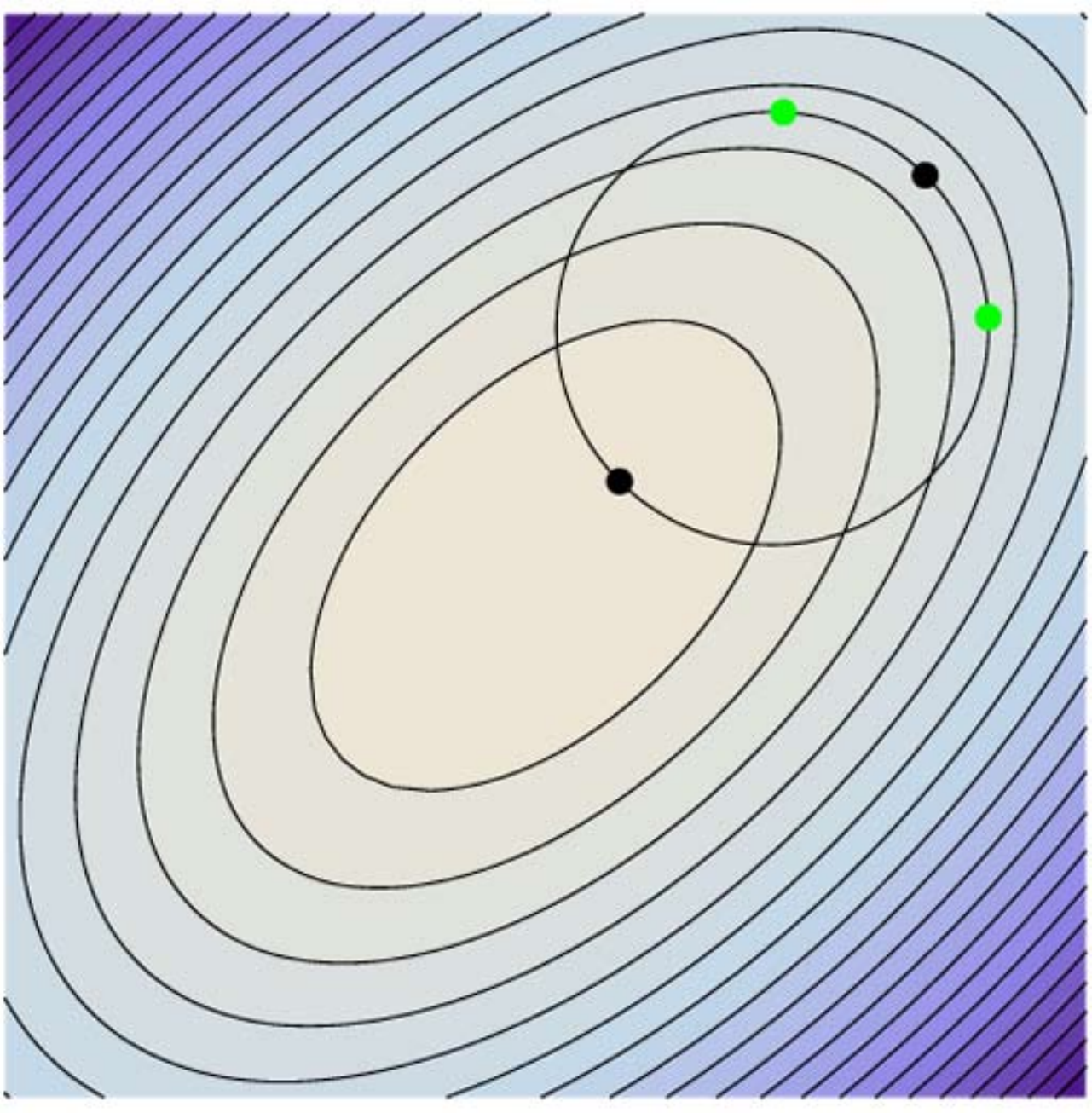}}
\centerline{\includegraphics[width=0.15\hsize]{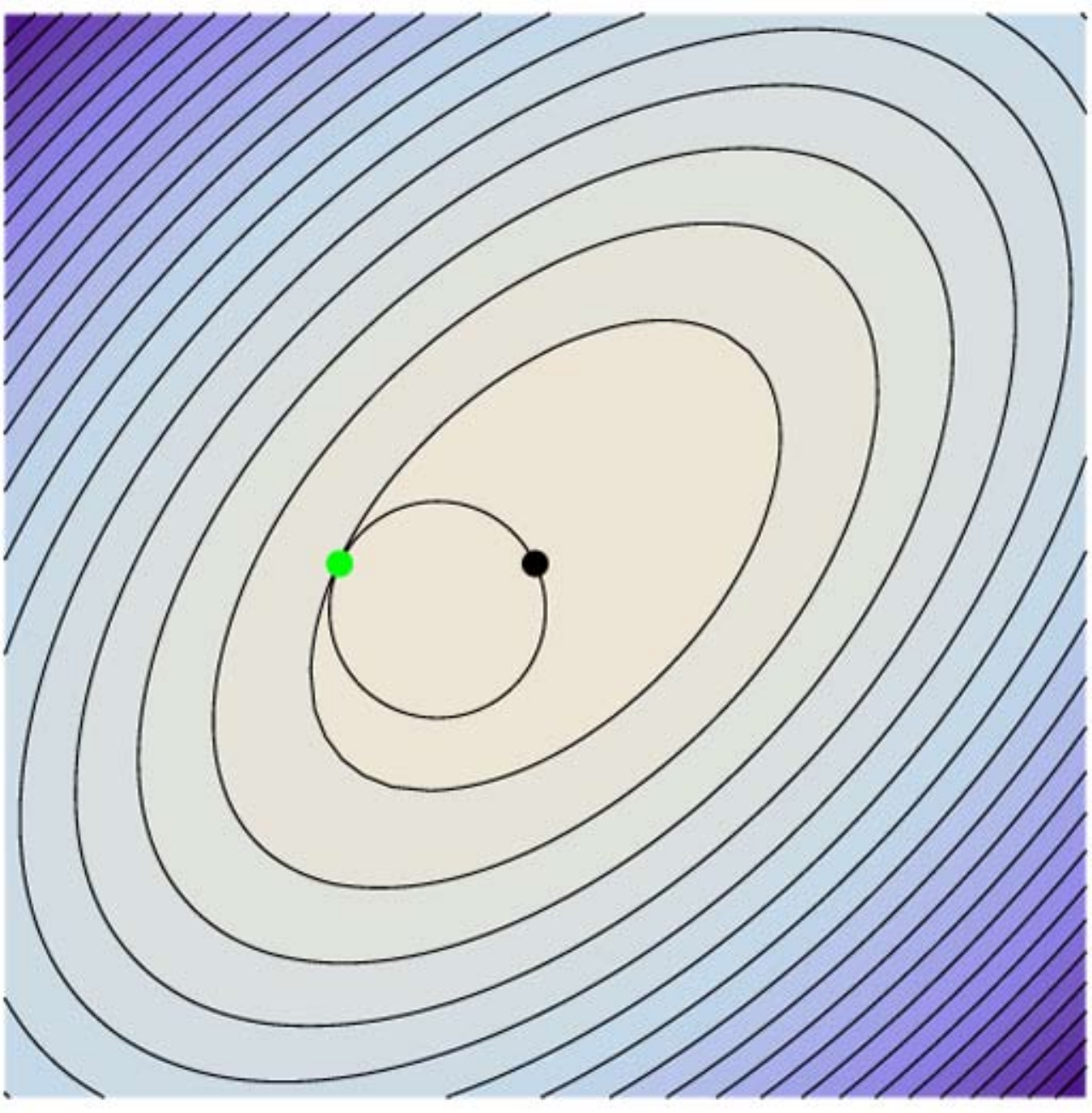}
\includegraphics[width=0.15\hsize]{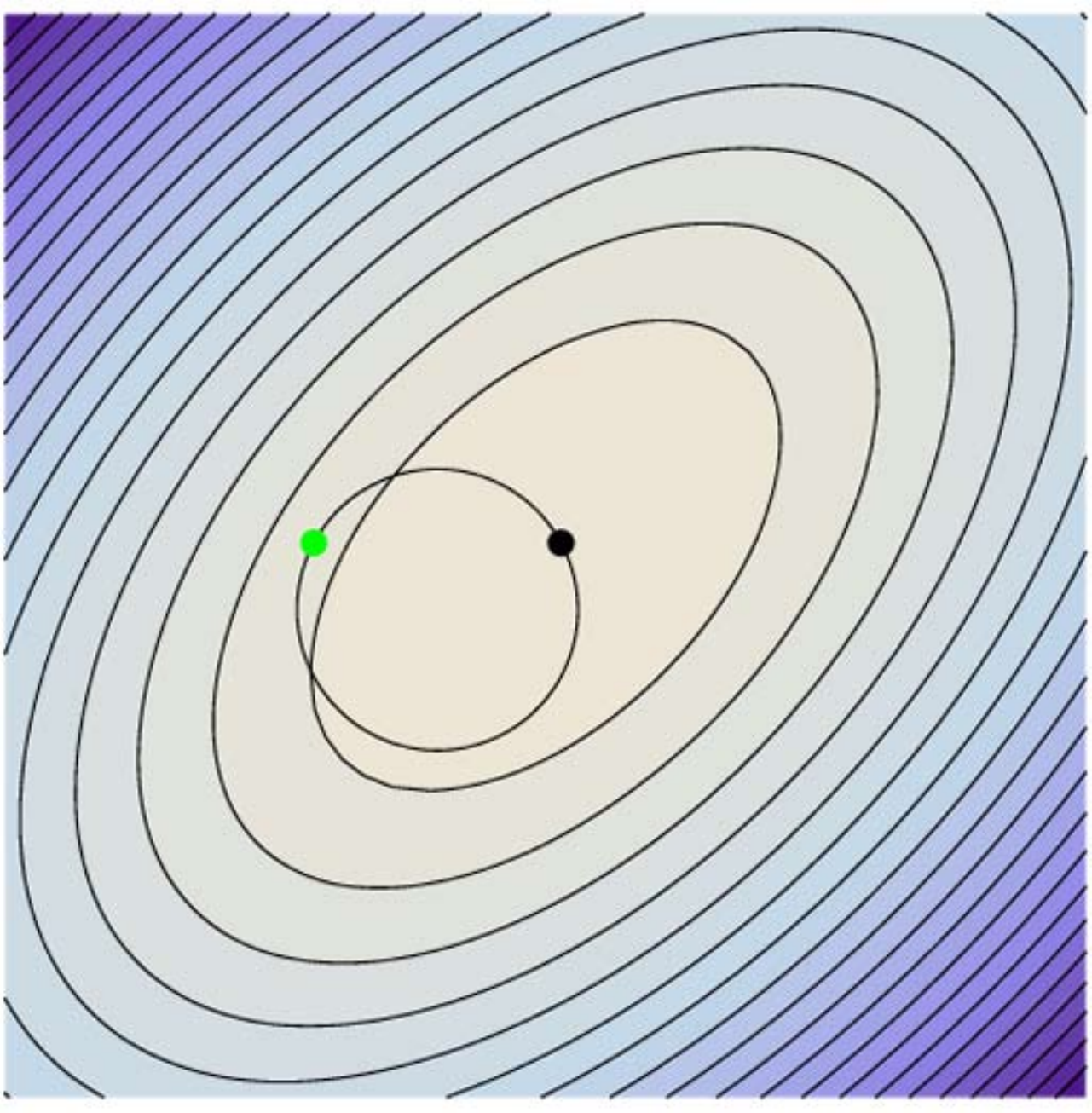}
\includegraphics[width=0.15\hsize]{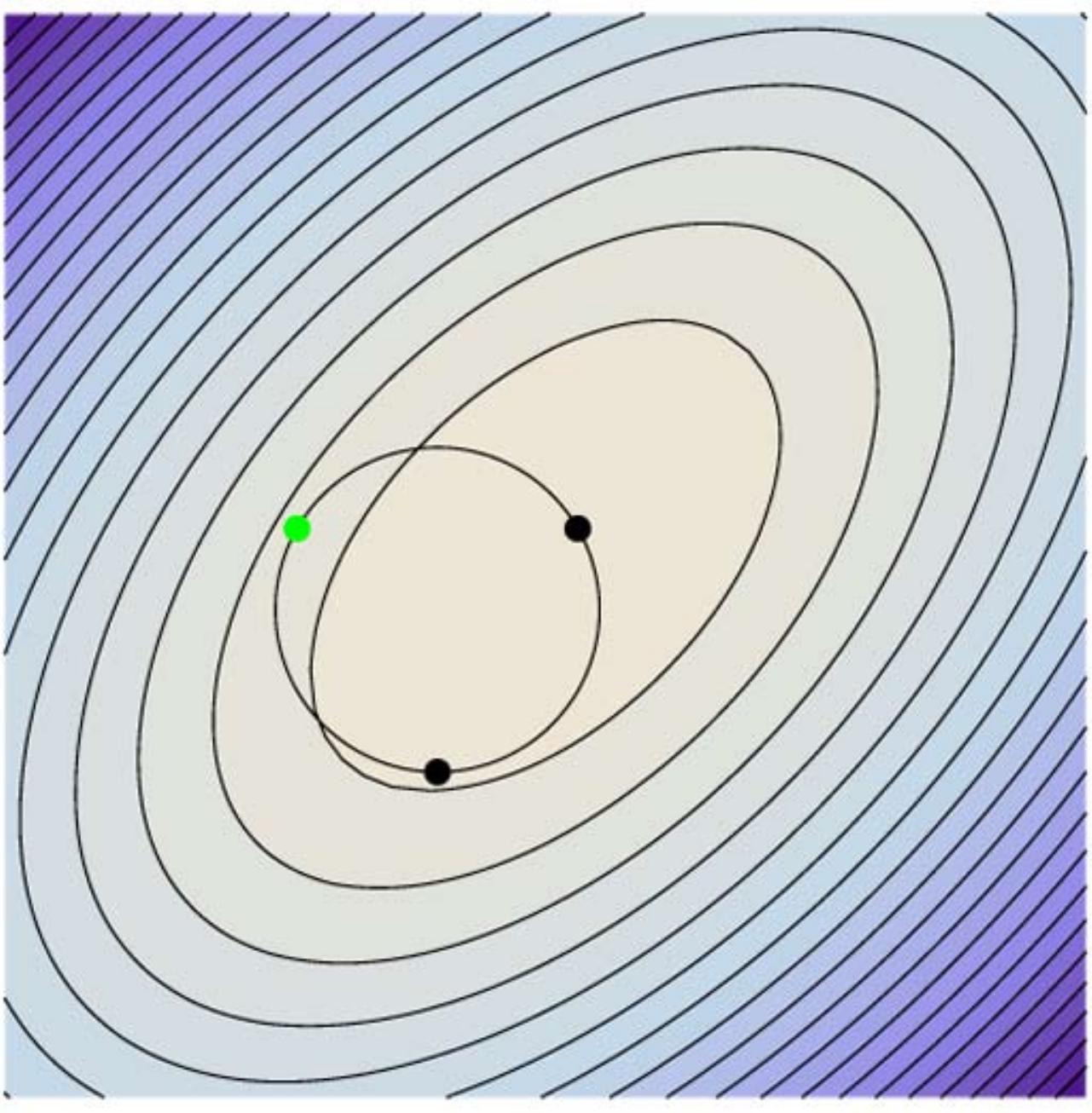}
\includegraphics[width=0.15\hsize]{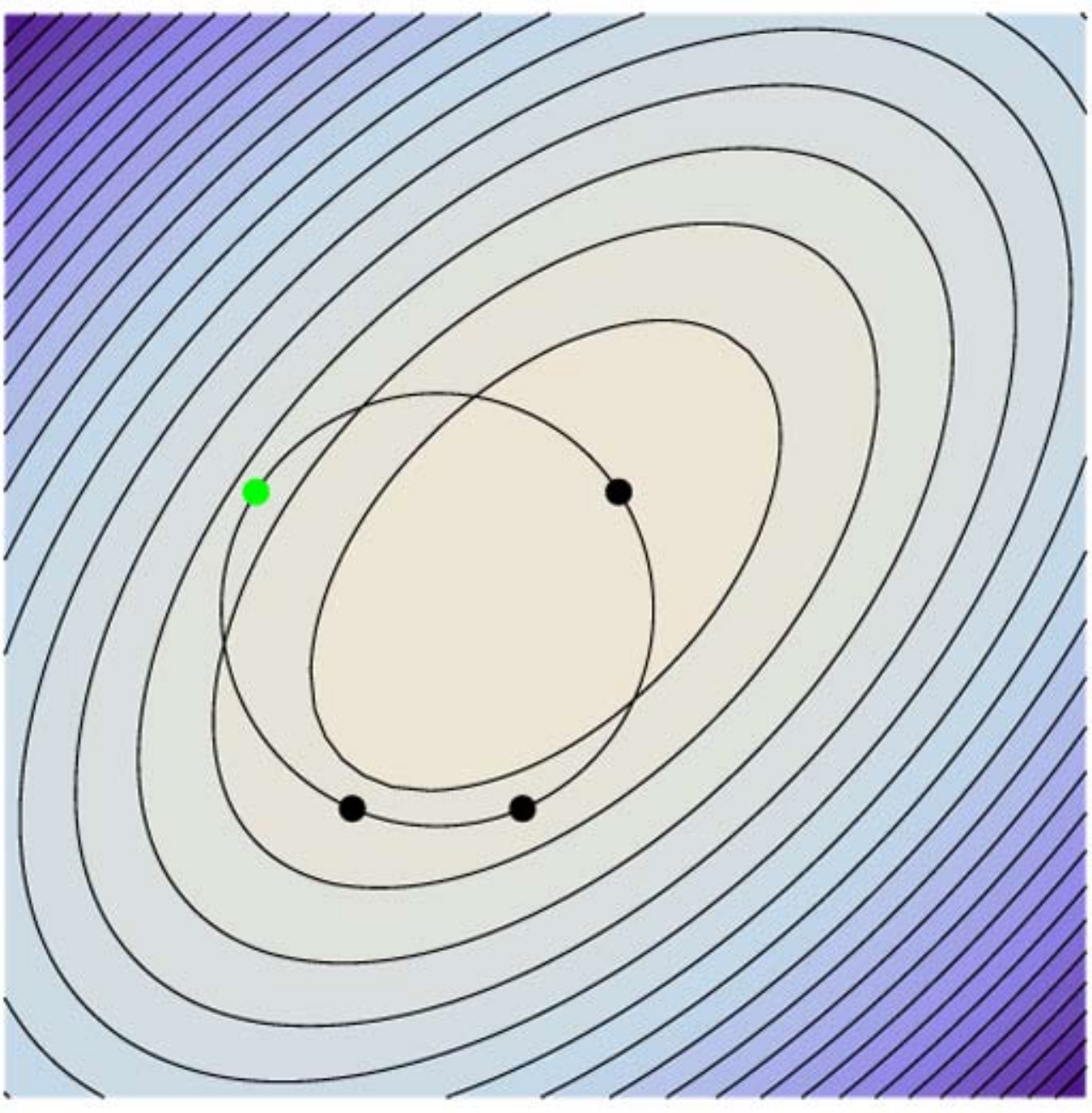}}
\caption{
Various plots of the solutions of optimization problem (\ref{optim}) arising in
$\Xi_{t}^{1,\frac{1}{2}}$. The horizontal axis, is the $p'$ axis, whereas the vertical axis is the $q'$-axis. Top row $(p,q)=(0.02,0.02)=0.02 \ul{k}_{1}$, middle row $(p,q)=(0.05,-0.05)=-0.05 \ul{k}_{2}$, bottom row $(p,q)=(-0.1,-0.05)$ for increasing time-frames $t>0$ from left to right. The minimizer(s) on the circle is(are) indicated in Green, whereas the other spurious solution(s) of Eq.~(\ref{pqaccent}) and Eq.~(\ref{4}) is(are) indicated in black.
% Bottom row: explanation color-coding of the phase in the corresponding Gabor transforms in the top row.
}
\label{fig:erosioncircles}
\end{figure}

In case $\lambda_{t}(p,q)$ is equal to an eigenvalue of $\textrm{Re}(B_{r,b})$, say $\lambda_{l}$ then (\ref{ELa}) yields $\lambda_{k}\alpha_{k}(p',q')=\lambda_{t}(p,q)(\alpha_{k}(p',q')-\alpha_{k}(p,q))$, $k=1,2$, from which we deduce that Eq.~(\ref{4}) is still valid if the resolvent does not exist, since then
$\lambda_{t}(p,q)=\lambda_{l} \Rightarrow \alpha_{l}(p,q)=0$. If $\alpha_{l}(p,q)=0$ then $(p,q)$ is aligned with the $(3-l)$-th principle axis of the ellipsoids $\{(p,q)\in \R^2 \; |\; (p,q) \textrm{Re}(B_{r,b}) (p,q)^{T}=c\}$, $c>0$ and the minimum $e^{(t+|\alpha_{3-l}|)^2\lambda_{3-l}}$ is obtained at $(p',q')=(p,q)+ t \textrm{sgn}(\alpha_{3-l}(p,q)) \ul{k}_{3-l}= (\alpha_{3-l}(p,q)+t \textrm{sgn}(\alpha_{3-l}))\, \ul{k}_{3-l}$, where only for $l=2$ we get the extra condition $t \leq t_{\textrm{max}}(p,p\tau)$. See the second row of Figure \ref{fig:erosioncircles}, where in the third column $t$ is chosen slightly larger than $t_{\textrm{max}}(p,p\tau)$. One can see that the minimum moves along the straight main principal direction until $(p^{*}(p,\tau p),\tau p^{*}(p,\tau p))$ is reached at time  $t_{\textrm{max}}(p,p\tau)$, where the single minimum is cut into two minima that evolve to the side.

For the $l=2$-case $\alpha_{2}(p,q)=0$ (and $t\leq t_{\textrm{max}}(p,p\tau)$) and the $l=1$-case $\alpha_{1}(p,q)=0$, see respectively first and second row of Figure \ref{fig:erosioncircles}. In these cases we find (by Eq.~(\ref{helper})) the Lagrange multiplier
\begin{equation} \label{exactl}
\alpha_{l}=0 \Rightarrow \lambda_{t}=(1+\frac{|\alpha_{3-l}|}{t})\lambda_{3-l} \end{equation}
which shows us that we indeed have a continuous transition of the solution across the principle directions in Eq.~(\ref{Xisol}). Finally, we note that along the ray
$\R^{+} \ni \mu \mapsto (1+\mu)(p^{*}(p,p\tau),\tau p^{*}(p,p\tau))$ there do not exist global minimizers of the optimization problem given by (\ref{optim}).
%Finally, the case $\alpha_{1}(p,q)=\alpha_{2}(p,q)=0$ occurs only if $(p,q)=(0,0)$ in which case the solution is %given by $e^{t \lambda_{1}}$, which is contained in (\ref{sol}).
%Finally we note that the center $(p',q')=(0,0)$ of the ellipsoidal quadratic form that is being minimized is contained within the open disk domain of our optimization iff $p^{2}+q^2<t^2$.
%
%iff $p^{2}+q^{2}<t^2$
%Now that general solution
%\end{proof}
$\hfill \Box$

\begin{figure}
\centerline{\includegraphics[width=0.75\hsize]{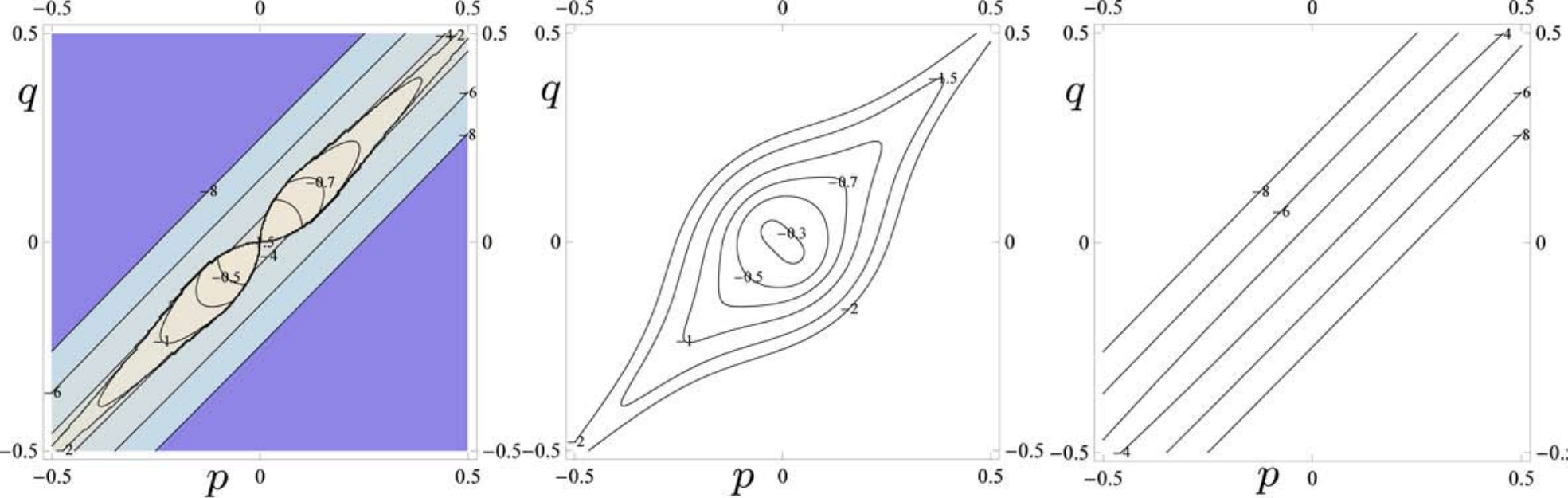}
}
\caption{The Lagrange multiplier $\lambda_{t}(p,q)$ for $(b,r,t)=(1,1,0.1)$ used in exact $\Xi_{t}^{1,1}f$, cf. Figure \ref{Fig:Exact2} and the corresponding approximations, in Eq. (\ref{approx}) for respectively $t > t_{\textrm{max}(p,q)}$ and $t \leq t_{\textrm{max}(p,q)}$ (where we have set $c=\frac{1}{2}$).
}\label{Fig:lambdat}
\end{figure}

\begin{remark}
We did not succeed in finding more tangible exact closed form expressions than the general tedious Gardano formula for the Lagrange-multipliers $\lambda_{t}(p,q)$ which are zeros of the fourth order polynomial given by Eq. (\ref{4}). For a plot of the graph of $(p,q) \mapsto \lambda_{t}(p,q)$, for $(b,r)=(\frac{1}{2},1)$, see Figure \ref{Fig:lambdat}.
These plots do suggest a close approximation of the type
{\small
\begin{equation} \label{approx}
\lambda_{t}(p,q) \approx %\overset{?}{=}
\left\{
\begin{array}{ll}
\lambda_{1}
%\left(1+ \frac{\sqrt{|\alpha_{1}(p,q)|^2+|\alpha_{2}(p,q)|^2\left( %\frac{|\alpha_{1}(p,q)|^2+ |\alpha_{2}(p,q)|^2}{t^2} %\right)^{\frac{3}{2}}}}{t} \right)
\left(1+ \frac{\sqrt{|\alpha_{1}(p,q)|^2+|\alpha_{2}(p,q)|^2\left( \frac{|\alpha_{2}(p,q)|^2+ |\alpha_{1}(p,q)|^2}{t^2} \right)^{\frac{3}{2}}}}{t} \right)
& \textrm{for }t_{\textrm{max}}(p,q) < t \\
\lambda_{2}\left(1+ \frac{\sqrt{|\alpha_{2}(p,q)|^2+c^2 \frac{|\lambda_{1}|^2}{|\lambda_{2}|^2}|\alpha_{1}(p,q)|^2}}{t} \right)
&\textrm{for }t_{\textrm{max}}(p,q) \geq t\ ,
\end{array}
\right.
\end{equation}
}
\noindent
with $c \in (0,1]$, where we recall that $t_{\textrm{max}}(p,q)$ is defined in
Theorem \ref{th:exact}. This approximation obeys the scaling property (\ref{scaling}) and is exact along the principal axes (i.e. exact if $\alpha_{1}=0$ or $\alpha_{2}=0$).
 % and is exact along $\ul{k}_{2}$ (i.e. $\alpha_{1}=0$). % and is exact along the principal axes.
%With respect to the second estimate, we note that $\lambda_{t} \leq \lambda_{2} \leq \lambda_{1}<0$ implies %that
%{\small $\left(\frac{\lambda_{t}-\lambda_{2}}{\lambda_{t}-\lambda_{1}} \right)^{2} \in (0,1)$}. This second %estimate is obtained by fixing this fraction to $c$ in Eq. (\ref{4}).
\end{remark}
\begin{corollary} \label{Cor:int}
Let $\lambda_{1},\lambda_{2}<0$ be the eigenvalues of $\textrm{Re}(B_{r,b})$ with $|\lambda_{1}|<|\lambda_{2}|$ and with corresponding eigenvectors $\ul{k}_{1}$, $\ul{k}_{2}$.
The function $\Xi_{t}^{1,\frac{1}{2}}f:H_{r}\to \mathbb{C}$ converges pointwise towards $0$ as $t \to \infty$.
%and $\lambda_{t}$ converges pointwise
%towards $\lambda_{2}$.
The isocontours of $|\Xi_{t}^{1,\frac{1}{2}}f|$ with $t>0$ are
non-ellipsoidal Jordan curves that retain the reflectional symmetry in the principal axes of $\textrm{Re}(B_{r,b})$.
The anisotropy of these Jordan curves (i.e. the aspect ratio of the intersection with the principal axes) equals
\begin{equation} \label{ani}
%\frac{\lambda_{2}}{\lambda_{1}} %\frac{|\lambda_{1}-\lambda_{t}(p,q)|^2}{|\lambda_{2}-\lambda_{t}(p,q)|^2},
\frac{\sqrt{\frac{c}{\lambda_{1}}}-t}{\sqrt{\frac{c}{\lambda_{2}}}-t}
\end{equation}
which tends to $\infty$ as $t \uparrow t_{\textrm{fin}}:=\sqrt{\frac{c}{\lambda_{2}}}$, which is the finite final time where the isocontour {\small $\{(p,q)\;|\; |\Xi_{t}^{1,\frac{1}{2}}f(p,q,s)|=c\}$} collapses to the span of $\ul{k}_{1}$.
\end{corollary}
\textbf{Proof }
Set $\eta=\frac{1}{2}$.
By the semi-group property of the erosion with kernel Eq.~(\ref{kernel}) and Eq.~(\ref{erosion})
(i.e. $k_{t}^{\eta} \ominus k_{\delta}^{\eta}=k_{t+\delta}^{\eta}$) and the fact
that the function {\small $|\Xi_{0}^{1,\frac{1}{2}}f|=|\mathcal{G}_{\psi}f|$}
vanishes at infinity and has a single critical point;
a maximum at $(p,q)=(0,0)$ it follows that
\[
\begin{array}{l}
|\Xi_{t+\delta}^{1,\eta}f|= k_{\delta}^{\eta} \ominus |\Xi_{t}^{1,\eta}f| < |\Xi_{t}^{1,\eta}f|, \
|\Xi_{t}^{1,\eta}f| \textrm{ attains a single extremum; that is a maximum at }(p,q)=(0,0), \\
|\Xi_{t}^{1,\eta}f| \textrm{ is continuous },
|\Xi_{t}^{1,\eta}f|(p,q) \to 0 \textrm{ as } \|(p,q)\| \to \infty,
\end{array}
\]
for all $t,\delta>0$. Consequently, the isocontour of $|\Xi_{t+\delta}^{1,\eta}f|$ with value $c>0$ is a
Jordan curve that is strictly contained in the interior of the isocontour of $|\Xi_{t}^{1,\eta}f|$ with the same value.
By Theorem \ref{th:exact1} the isocontours are ellipsoidal for $t=0$. For $t>0$ these are no longer
ellipsoidal since $(p,q)\mapsto \lambda_{t}(p,q)$
is not constant in Eq.~(\ref{Xisol}). For $t \to \infty$ the kernel $k_{t}^{\eta}$, cf.~Eq.~(\ref{kernel}),
vanishes and the erosion, cf.~Eq.~(\ref{erosion}), with the kernel
$k_{t}^{\eta}$ converges for each point $(p,q) \in \R^2$ to the global minimum of
$|\mathcal{G}_{\psi}f|$ which is zero.

Eq.~(\ref{4}) is invariant under $\alpha_{i} \mapsto -\alpha_{i} $ and thereby
$\lambda_{t}$ is invariant under reflections in the principal axes of $\textrm{Re}(B_{r,b})$.
As a result, by Eq.~(\ref{Xisol}), all isocontours of $|\Xi_{t}^{1,\eta}f|$ are invariant under these reflections.
Such an isocontour is given by
\begin{equation} \label{help3}
(\lambda_{t}(p,q))^2 \left(
 \frac{|\alpha_{1}(p,q)|^2\lambda_{1}}{(\lambda_{1}-\lambda_{t}(p,q))^2}+ \frac{|\alpha_{2}(p,q)|^2\lambda_{2}}{(\lambda_{2}-\lambda_{t}(p,q))^2}
 \right)=C,
\end{equation}
for some $C\leq 0$. Now set $\alpha_{1}=0$ in Eq.(\ref{help3}) and solve for $\alpha_{2}$, then set $\alpha_{2}=0$ in Eq.(\ref{help3}) and solve for $\alpha_{1}$. Here one should use the exact formula, Eq.~(\ref{exactl}) for $\lambda_{t}$ that holds along the principal axes. Finally, division of the 2 obtained results yields Eq.~(\ref{ani}) from which the result follows.
$\hfill \Box$\\ \\
See Figure \ref{fig:erosioncircles} for solutions of the optimization problem (\ref{optim}) for various settings of $(p,q,t) \in \R^{2} \times \R^{+}$.
\subsection{Plots of the Exact Re-assigned Gabor Transforms of Chirp Signals}
\begin{figure}
\centerline{\includegraphics[width=0.65\hsize]{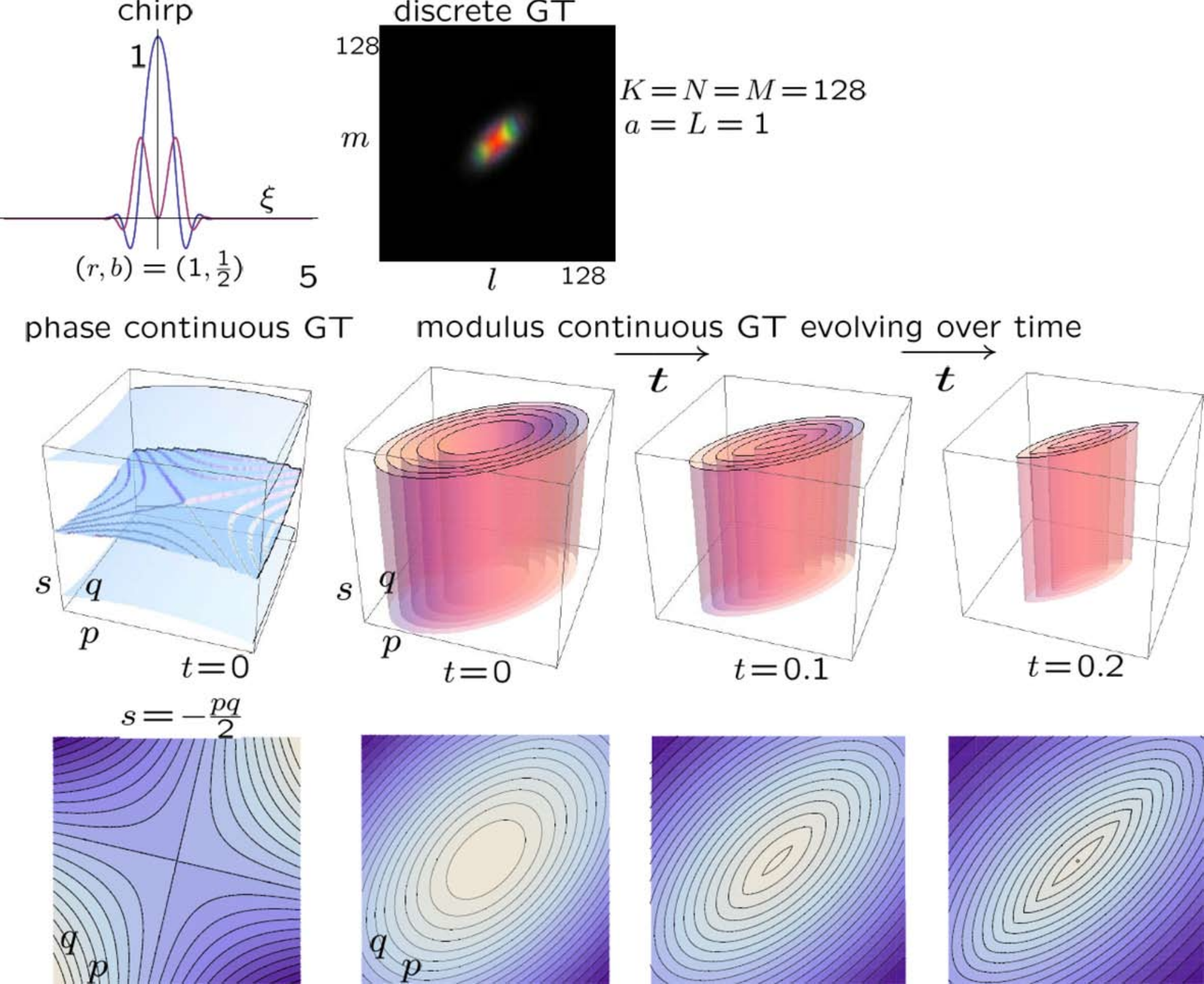}}
\caption{Visualization of
$(p,q,s) \mapsto (\Xi_{t}^{1,\frac{1}{2}}f)(p,q,s)$ where $f$ is the chirp signal (\ref{chirp}) with $(b,r)=(\frac{1}{2},1)$. Top row: left; the input chirp signal, right: We have sampled the chirp signal on $128$-grid and depicted the corresponding discrete Gabor transform ($K=M=N, L=a=1$) where color represents phase and where intensity represents the modulus. Middle row; equisurfaces of the phase and modulus of our exact solutions for {\tiny $\Xi_{t}^{1,\frac{1}{2}}f$} in Theorem \ref{th:exact}, phase is constant over time whereas the modulus is considered for $t=0$, $t=0.1$ and $t=0.2$. Iso-intensities have been fixed ($-0.05,-0.09,-0.14, -0.18, -0.23$) over time to show their collapsing behavior towards $\ul{k}_{1}$ in accordance with Corollary \ref{Cor:int}. Bottom row; Restrictions to the surface $s=-\frac{pq}{2}$ yields the solutions in phase space.
}\label{Fig:Exact}
\end{figure}

\begin{figure}
\centerline{\includegraphics[width=0.65\hsize]{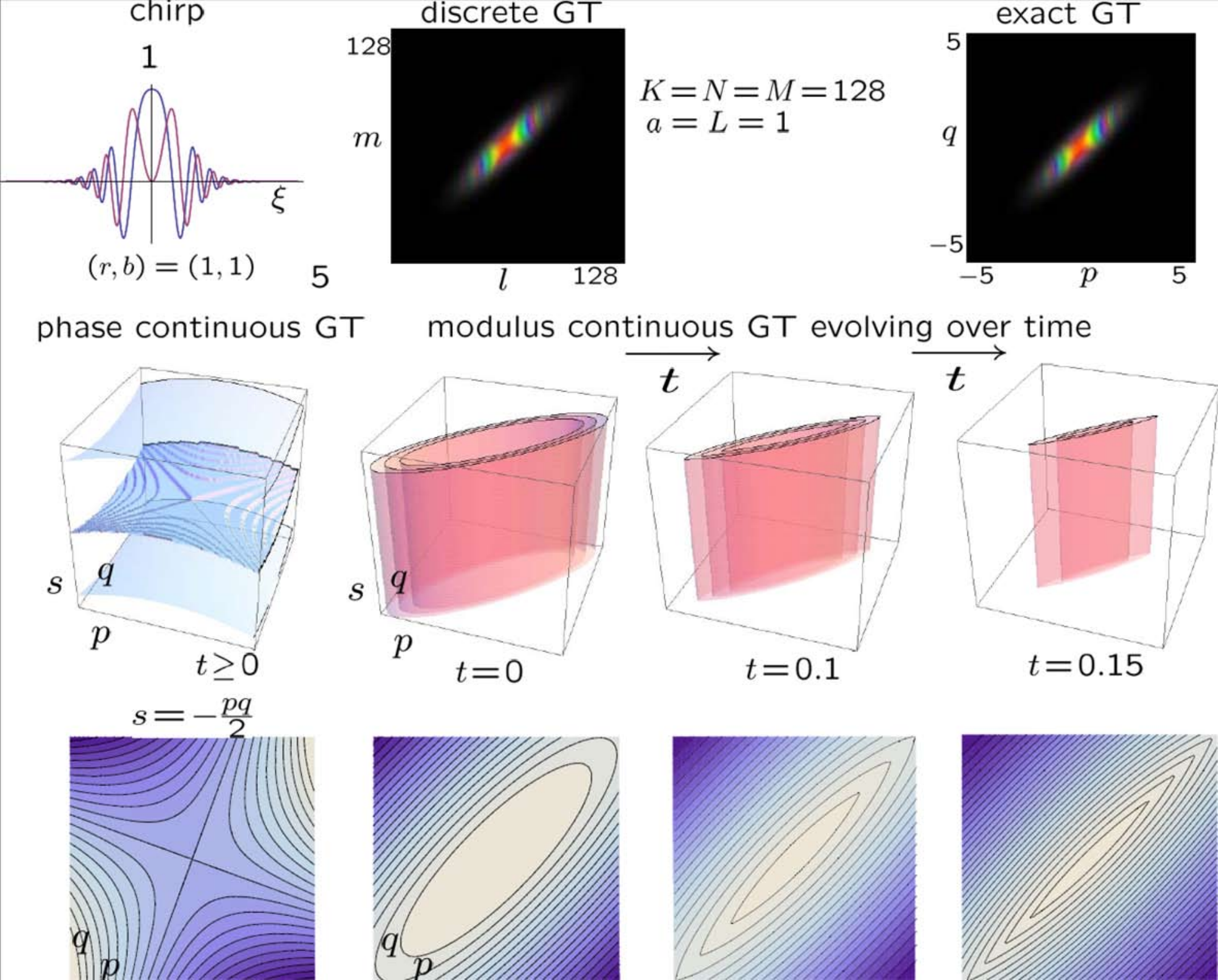}
}
\caption{Same settings as in Figure \ref{Fig:Exact} except for $(b,r)=(1,1)$ and iso-intensities ($-0.06,-0.08,-0.1$).
}\label{Fig:Exact2}
\end{figure}

For a plot of the exact $\Xi_{t}^{1,\frac{1}{2}}f$ in Theorem \ref{th:exact} using Gardano's formula for the zeros of the fourth order polynomial given by Eq. (\ref{4}) we refer to Figure \ref{Fig:Exact} (for the case $(b,r)=(\frac{1}{2},1)$) and Figure \ref{Fig:Exact2} (for the case $(b,r)=(1,1)$). Within these Figures one can clearly see that during the erosion process isocontours collapse towards the eigenspace $\langle \ul{k}_{1} \rangle $ with smallest eigenvalue. Moreover for $(b,r)=(1,1)$ one has $\ul{k}_{1} \approx \frac{1}{\sqrt{1+r^2}}(1,r)^{T}$ and finally, we note that the discrete Gabor transforms closely approximate the exact Gabor transforms as can be seen in the top row of Figure \ref{Fig:Exact2}. Although that, for $\eta=\frac{1}{2}$, the eroded ellipses are no longer ellipses the principal axis $\langle \ul{k}_{1}\rangle$ is preserved as a symmetry-axis for the iso-lines in phase space (where $s=-\frac{pq}{2}$).
\section{Left Invariant Diffusion on Gabor transforms \label{ch:covdiff}}

A common technique in imaging to enhance images $f:\R^{2}\to \R$ via non-linear adaptive diffusions are so-called coherence enhancing diffusion (CED) schemes, cf.~\!\cite{Weickert1999,Cottet,Nitzberg}, where one considers a diffusion equation, where the diffusivity/conductivity matrix in between the divergence and gradient is adapted to the Gaussian gradient (structure tensor) or Hessian tensor computed from the data.
Here the aim is to diffuse strongly along edges/lines and little orthogonal to them. In case of edge-adaptation the anisotropic diffusivity matrix is diagonalized along the eigenvectors of auxiliary
matrix
\begin{equation} \label{CED}
\ul{A}(u_f(\cdot,s))(x,y)= (G_{\epsilon}* (\nabla u_f(\cdot,s)  (\nabla u_f(\cdot,s))^T))(x,y)
\end{equation}
where $\nabla u_{f}$ denotes either the Gaussian gradient of image $f$ or the
gradient of the non-linearly evolved image $u_f(\cdot,s)$ and where the final convolution, with Gaussian kernel $G_{\epsilon}(x,y)=(4\pi \epsilon)^{-1}e^{-\frac{\|\ul{x}\|^{2}}{4\epsilon}}$ with $\epsilon >0$ small, is applied componentwise. In case of line-adaptation the auxiliary matrix is obtained by the Hessian:
\begin{equation}\label{CED2}
(\ul{A}(u_f(\cdot,s)))(\ul{x})=  (G_{\epsilon} * H u_f(\cdot,s))(\ul{x})= [(G_{\epsilon} * \partial_{x_i} \partial_{x_j} u_f(\cdot,s))(\ul{x})]_{i,j=1,2}= [( \partial_{x_i} \partial_{x_j} G_{\epsilon}* u_f(\cdot,s))(\ul{x})]_{i,j=1,2}\ .
\end{equation}
Given the eigen system of the auxiliary matrix, standard coherence enhancing diffusion equations on images cf.\!~\cite{Weickert1999} can be formulated as
{\scriptsize
\begin{equation} \label{CEDimages}
\left\{
\begin{array}{ll}
\partial_{s} u_{f}(\ul{x},s) & = \left(\nabla_{\ul{x}} \cdot \;
\mathbf{S}
\cdot\;
\left(
\begin{array}{cc}
\ve & 0 \\ 0& (1-\ve)e^{-\frac{c}{\lambda_{1}-\lambda_{2}}}+\ve
\end{array}
\right)
\cdot \;  \mathbf{S}^{-1} \cdot \; (\nabla_{\ul{x}} u_{f}(\cdot,s))^{T}\right)(\ul{x}) \\
 & =
\left(
\begin{array}{cc}
 \partial_{u} & \partial_{v}
 \end{array}
\right)
\left(
\begin{array}{cc}
\ve & 0 \\ 0& (1-\ve)e^{-\frac{c}{(\lambda_{1}-\lambda_{2})^2}}+\ve
\end{array}
\right)
\left(
\begin{array}{c}
\partial_{u} \\
\partial_{v}
\end{array}
\right) u_{f}(\ul{x},s), \ \   \ul{x} \in \R^2, s>0,
 \\[7pt]
u_{f}(\ul{x},0) &=f(\ul{x}), \qquad \ul{x} \in \R^2.
\end{array}
\right.
\end{equation}
}
Here we expressed the diffusion equations in both the global standard basis {\small $\{\ul{e}_{x},\ul{e}_{y}\}:=\{(1,0),(0,1)\} \leftrightarrow \{\partial_{x},\partial_{y}\}$} and in the locally adapted basis of eigenvectors of auxiliary matrix $\ul{A}(u_f(\cdot,s))(\ul{x})$, recall Eq.~\!(\ref{CED}) and Eq.~\!(\ref{CED2}):
\[
\{\ul{e}_{1},\ul{e}_{2}\}:=\{\ul{e}_{1}\left(\ul{A}(u_f(\cdot,s))(\ul{x})\right), \ul{e}_{2}\left(\ul{A}(u_f(\cdot,s))(\ul{x})\right)\} \leftrightarrow \{\partial_{u},\partial_{v}\},
\]
with respective eigenvalues $\lambda_{k}:=\lambda_{k}(\ul{A}(u_f(\cdot,s))(\ul{x}))$, $k=1,2$.
%of the so-called structure tensor $S(u_f(\cdot,s))(\ul{x})$ at position $\ul{x}\in \R^2$ at time $s>0$.
The corresponding orthogonal basis transform which maps the standard basis vectors to the eigenvectors $\{\ul{e}_{1},\ul{e}_{2}\}$ is denoted by {\small $\mathbf{S}= (\ul{e}_{1}\; |\; \ul{e}_{2}) $} and we have $(\partial_{u}\; \;  \partial_{v})=(\partial_{x}\; \; \partial_{y}) \cdot \mathbf{S}$. Note that
\[
\ul{A}= \ul{S} \cdot\textrm{diag}\{\lambda_{1},\lambda_{2}\} \cdot \ul{S}^{-1}. \]
At isotropic areas $\lambda_{1} \rightarrow \lambda_{2}$ and thereby the conductivity matrix becomes a multiple of the identity yielding isotropic diffusion only at isotropic areas, which is desirable for noise-removal.

\begin{figure}[h]
%\vspace{-0.3cm}\mbox{}
\centering
\begin{tabular}{ccc}
%\vspace{-0.3cm}\hbox
\scriptsize original image &\scriptsize CED & \scriptsize CEDOS \\
\includegraphics[width=0.195\linewidth]{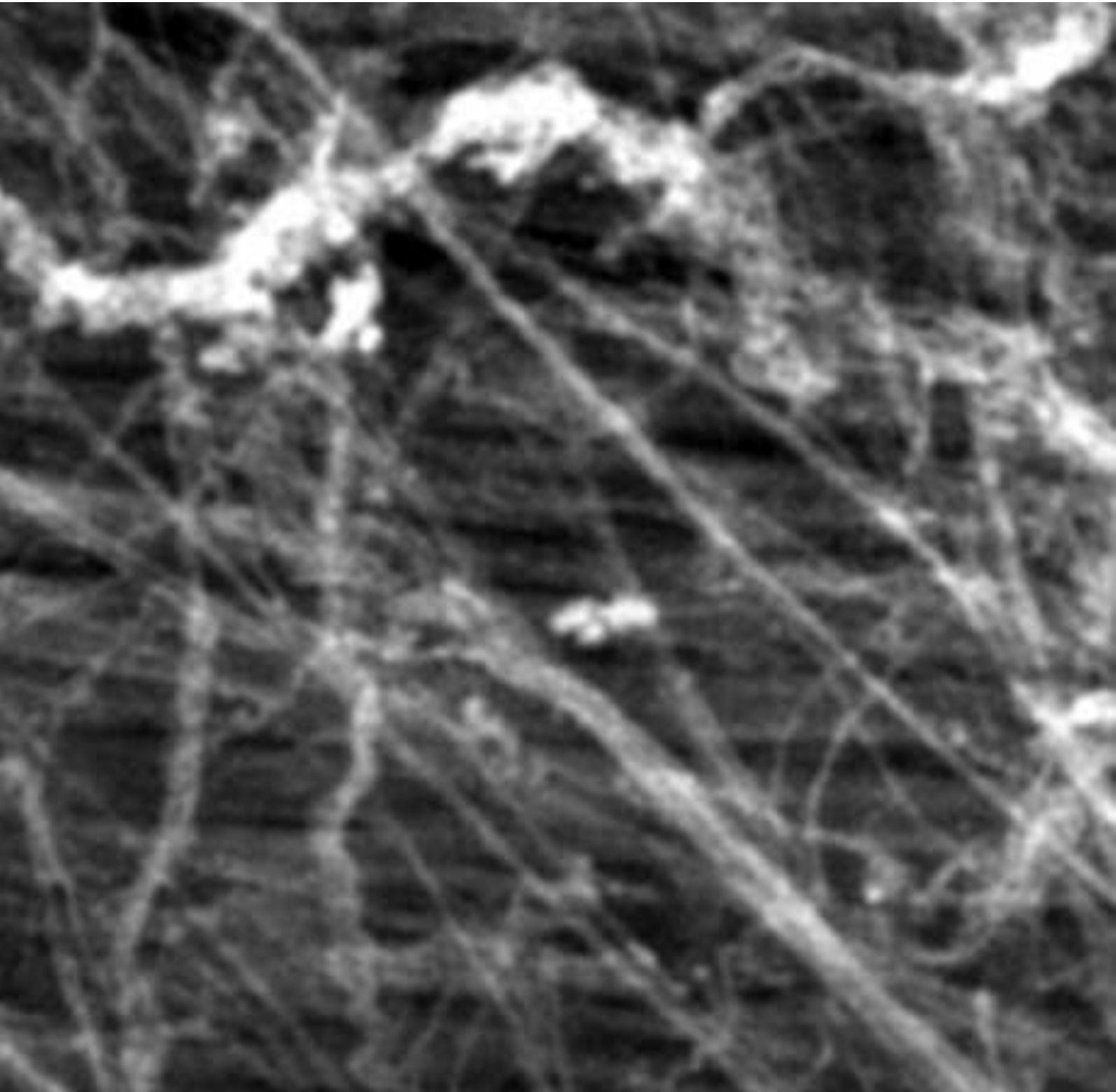}&
\includegraphics[width=0.195\linewidth]{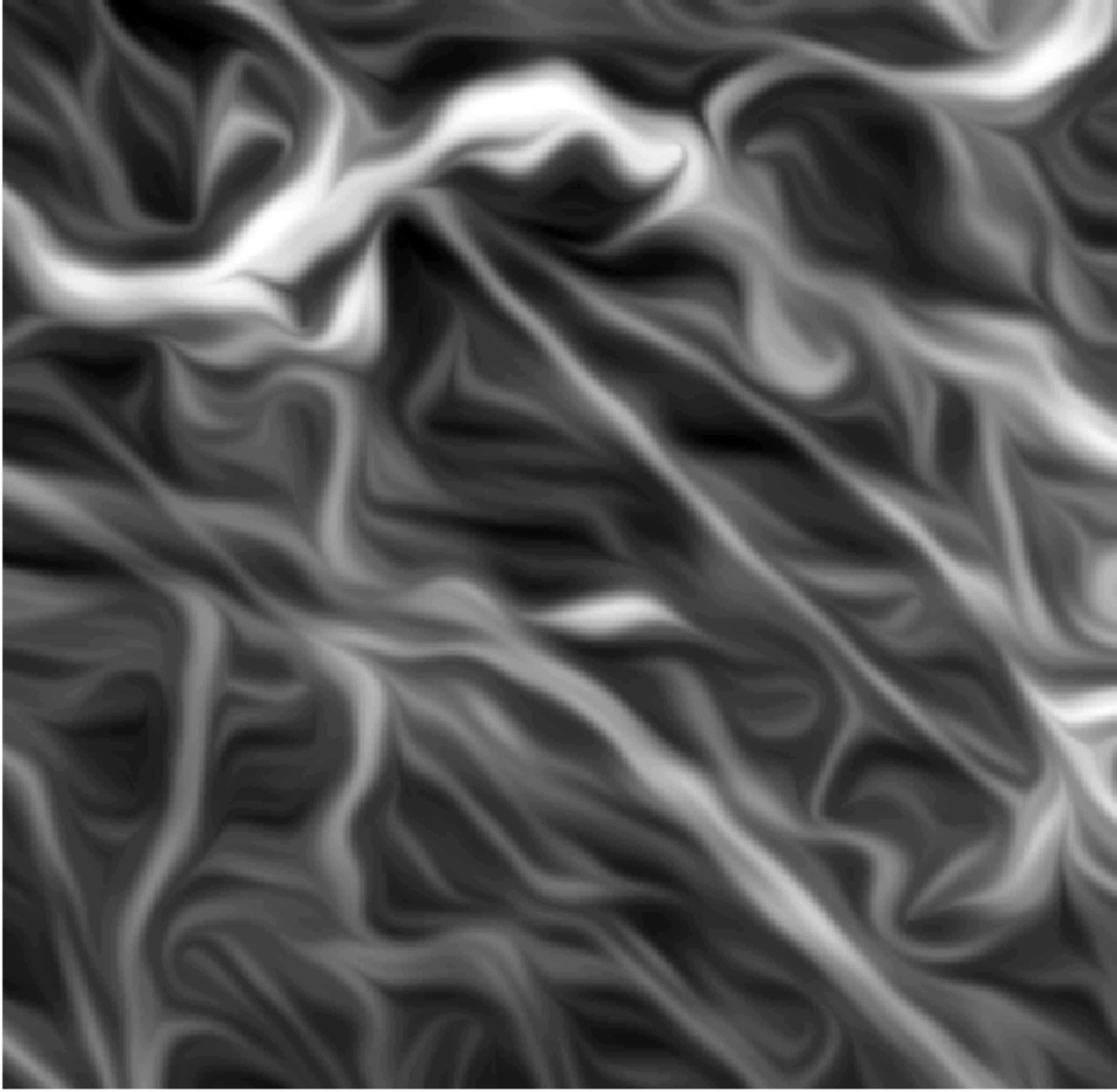}&
\includegraphics[width=0.195\linewidth]{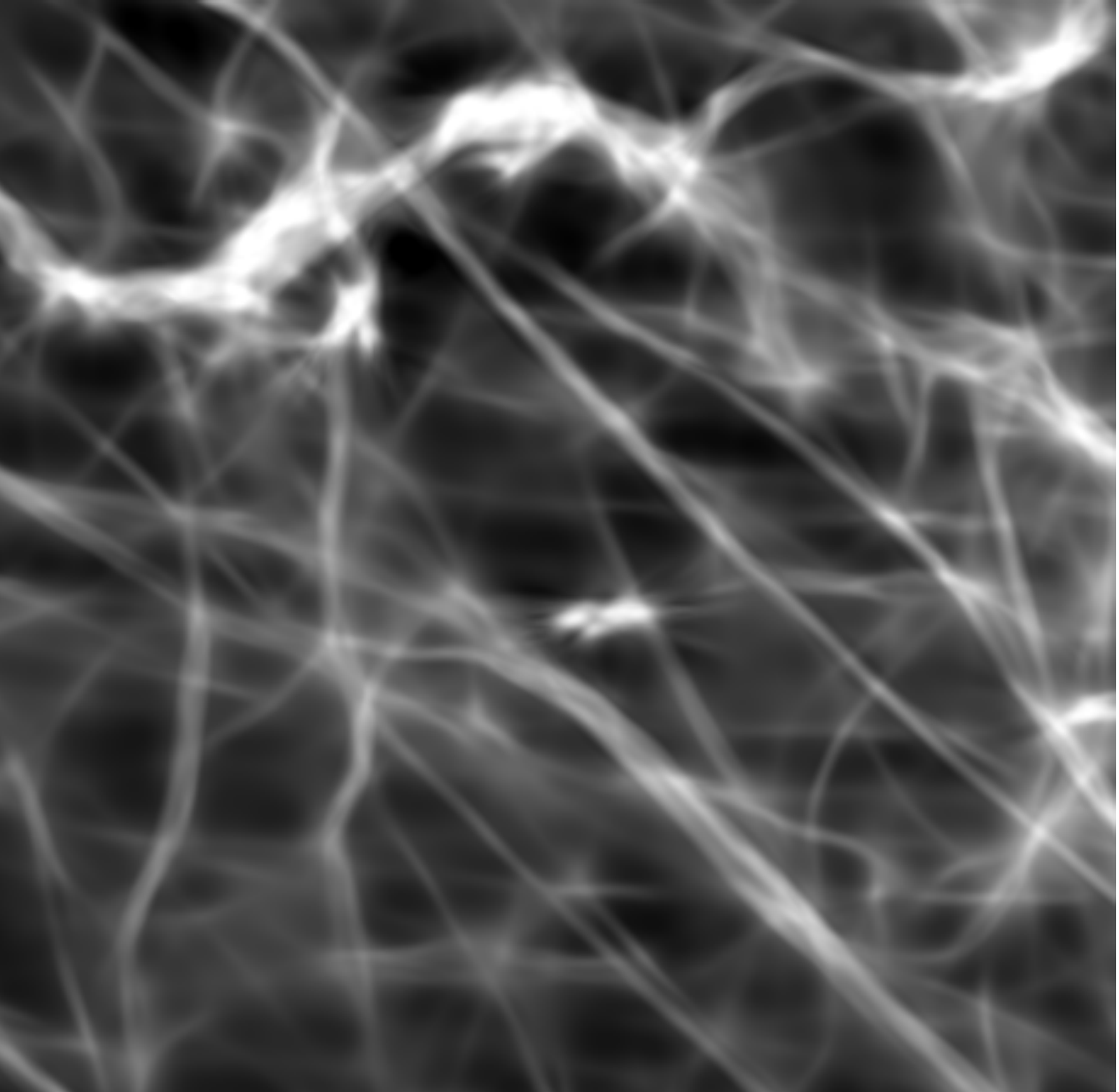} \\
\includegraphics[width=0.195\linewidth]{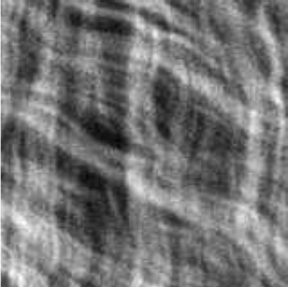}&
\includegraphics[width=0.195\linewidth]{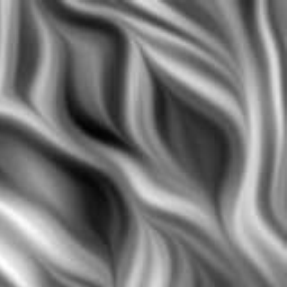}&
\includegraphics[width=0.195\linewidth]{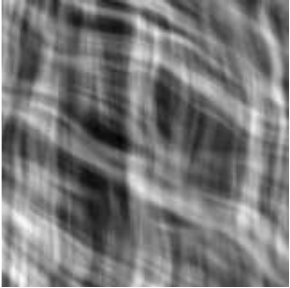}
\end{tabular}
\begin{tabular}{ccc}
\includegraphics[width=0.195\linewidth]{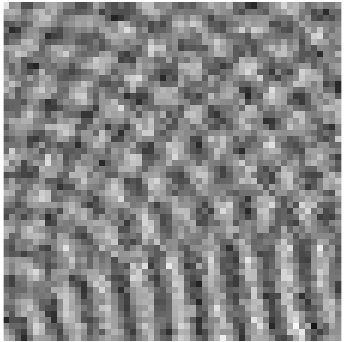}&
\includegraphics[width=0.195\linewidth]{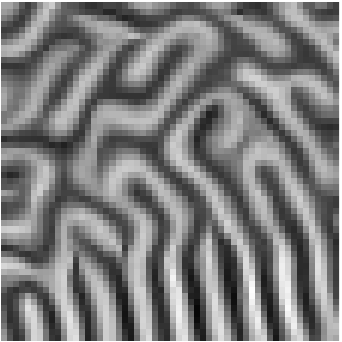}&
\includegraphics[width=0.195\linewidth]{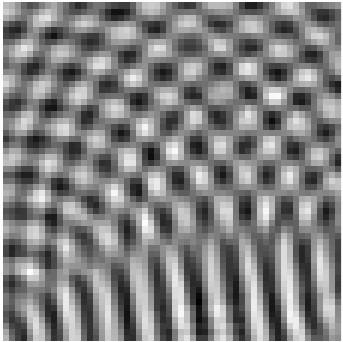}
\end{tabular}
\caption{Left-invariant processing via invertible orientation scores is the right approach to generically deal with crossings and bifurcations.
Left column: original images.
Middle column: result of coherence enhancing diffusion on images (CED), cf.~\!\cite{Weic99b}.
Right column: coherence enhancing diffusion via orientation score, cf.~\!\cite{DuitsAMS2,Fran2009}.
Top row: 2-photon microscopy image of bone tissue.
%2nd row: 2-photon microscopy image of muscle cell.
2nd row : collagen fibers of the heart.
3rd row: artificial pattern.
%All
%Typically, these applications clearly show that
coherence enhancing diffusion on orientation
scores (CEDOS) is capable of handling crossings and bifurcations, whereas (CED) produces spurious artifacts.
}\label{fig:OS}
\vspace{-0.4cm}\mbox{}
\end{figure}

A typical drawback of these coherence enhancing diffusion directly applied to images, is that %the
%local frame of reference is either based on
the direction of the image gradient (or Hessian eigenvectors) is ill-posed in the vicinity of crossings and junctions. Therefore, in their recent works on orientation scores, cf.~\cite{Fran2009,DuitsAMS2,FrankenPhDThesis} Franken and Duits developed
coherence enhancing diffusion via invertible orientation scores (CEDOS) where crossings are generically disentangled allowing crossing preserving diffusion, see Figure \ref{fig:OS}.
The key idea here is to extend the image domain $\R^{2}$ to a larger Lie group $SE(2)=\R^{2}\rtimes SO(2)$ where the left-invariant vector fields provide per spatial position a whole family of oriented reference frames cf.\!~\cite{DuitsR2006SS2,DuitsRTHESIS}. This allows the inclusion of coherent alignment (``context'') of disentangled local line fragments visible in the orientation score.

In this section we would like to extend this general idea to Gabor transforms defined on the
the Heisenberg group, where the left-invariant vector fields provide per spatial position a whole family of reference frames w.r.t. frequency and phase. Again we would like to achieve coherent alignment of all Gabor atoms in the Gabor domain via left-invariant diffusion.
%into separate oriented lines and where the inverse orientation score transform merges the enhanced lines to an %enhanced crossing.
\\
\\
In order to generalize the CED (coherence enhancing diffusion) schemes to Gabor transforms we must replace the left-invariant vector fields $\{\partial_{x},\partial_{y}\}$ on the additive group $(\R^{2},+)$, by the left-invariant vector fields on $H_{r}$. Furthermore, we replace the 2D-image by the Gabor transform of a $1D$-signal, the group $\R^{2}$ by $H_{r}$,
the standard inner product on $\R^2$ by the first fundamental form (\ref{metric}) parameterized by $\beta$.
%\footnote{The parameter $\beta^{-1}$ %(which equals %$a$ if $\beta=a^2$)
%has physical dimension length, so that the first fundamental form is consistent with physical dimensions.} %$\beta$.
Finally, we replace
the basis of normalized left-invariant vector fields $\{\ul{e}_{x},\ul{e}_{y}\}:=\{(1,0),(0,1)\} \leftrightarrow \{\partial_{x},\partial_{y}\}$ on $\R^2$ by the normalized left-invariant vector fields $\{\ul{e}_{x},\ul{e}_{y}\}:=\{(1,0),(0,1)\} \leftrightarrow \{\beta^{-2}\mathcal{A}_{1}, \mathcal{A}_{2}\}$ on $H_{r}$.

These steps produce
the following system for adaptive left-invariant diffusion on Gabor transforms:
%\[
%
%\]
{\scriptsize
\begin{equation}\label{CEDHr}
\hspace{-0.5cm}
\left\{
\begin{array}{ll}
\partial_{t} W(p,q,s,t) & = \left(
\left( \begin{array}{cc}
\beta^{-2}\,\mathcal{A}_{1} &  \mathcal{A}_{2}
\end{array} \right)
 \cdot
%\underset{{\tiny
%\epsilon \leftarrow \alpha}}{ S  }
\mathbf{S}
\cdot
{\small
\left(
\begin{array}{cc}
\ve & 0 \\ 0& (1-\ve)e^{-\frac{c}{\lambda_{1}-\lambda_{2}}}+\ve
\end{array}
\right)
}
  \cdot
  %\underset{{\tiny
 %\alpha \leftarrow \epsilon}}{ S  }
\mathbf{S}^{-1}
\left(
\begin{array}{c}
\beta^{-2}\,\mathcal{A}_{1} W \\
\mathcal{A}_{2} W
\end{array}
\right)
\right)(p,q,s,t) \\
 & =
\left(
\begin{array}{cc}
 \partial_{u} & \partial_{v}
 \end{array}
\right)
\left(
\begin{array}{cc}
\ve & 0 \\ 0& (1-\ve)e^{-\frac{c}{(\lambda_{1}-\lambda_{2})^2}}+\ve
\end{array}
\right)
\left(
\begin{array}{c}
\partial_{u} \\
\partial_{v}
\end{array}
\right) W(p,q,s,t), \ \  \\ &\textrm{ for all } (p,q,s) \in H_r , t>0 \textrm{ and fixed } c>0, \ve>0,
 \\[7pt]
W(p,q,s,0) &=\mathcal{W}_{\psi}f(p,q,s),  \textrm{ for all } (p,q,s) \in H_r.
\end{array}
\right.
\end{equation}
}
with {\small $\mathbf{S}:= (\ul{e}_{1}\; |\; \ul{e}_{2}) $} and $(\partial_{u}W\; \;  \partial_{v}W)=(\beta^{-2}\mathcal{A}_{1}W\; \; \mathcal{A}_{2}W) \mathbf{S}$ and where
\begin{equation}\label{evev}
\begin{array}{l}
\lambda_{k}:=\lambda_{k}(\ul{A}(|\mathcal{W}_{\psi}f|)(p,q,s)), k=1,2,  |\lambda_{1}|\leq |\lambda_{2}|\\
\ul{e}_{k}:=
\ul{e}_{k}(\ul{A}(|\mathcal{W}_{\psi}f|)(p,q,s)), k=1,2,
\end{array}
\end{equation}
denote the eigenvalues and $\{\ul{e}_{1},\ul{e}_{2}\} \leftrightarrow \{\partial_{u},\partial_{v}\}$ the corresponding normalized
eigenvectors of a local auxiliary matrix $\ul{A}(|\mathcal{W}_{\psi}f|)(p,q,s)$:
\[
\ul{A}(|\mathcal{W}_{\psi}f|)(p,q,s)\; \ul{e}_{k}= \lambda_{k}\, \ul{e}_{k}. % \textrm{ yielding } \ul{S}=(\ul{e}_{1}\;|\; \ul{e}_{2}).
\]
The eigenvectors are normalized w.r.t. first fundamental form (\ref{metric}) parameterized by $\beta$, i.e.
\[
\begin{array}{l}
\gothic{G}_{\beta}(\partial_{u},\partial_{u})=\gothic{G}_{\beta}(\partial_{v},\partial_{v})
=\gothic{G}_{\beta}(\beta^{-2}\mathcal{A}_{1},\beta^{-2}\mathcal{A}_{1})
=\gothic{G}_{\beta}(\mathcal{A}_{2},\mathcal{A}_{2})=1 \; \; \;
\leftrightarrow \\
\ul{e}_{1}^{T}\ul{e}_{1}=\ul{e}_{2}^{T}\ul{e}_{2}=\ul{e}_{x}^{T}\ul{e}_{x}=\ul{e}_{y}^{T}\ul{e}_{y}=1 \Rightarrow \mathbf{S}^{T}=\mathbf{S}^{-1}.
%\left(\underset{{\tiny
%\alpha \leftarrow \epsilon}}{ S  }\right)^T=\left(\underset{{\tiny
%\alpha \leftarrow \epsilon}}{ S  }\right)^{-1}=\underset{{\tiny
%\epsilon \leftarrow \alpha}}{ S  }.
\end{array}
\]
This auxiliary matrix $\ul{A}(|\mathcal{W}_{\psi}f|)(p,q,s)$ at each position $(p,q,s) \in H_{r}$ is chosen such that it depends only on the absolute value $|\mathcal{W}_{\psi}f|$ (so it is independent of phase parameter $s$) of the \emph{initial condition}. The same holds for the corresponding conductivity matrix-valued function appearing in Eq.~(\ref{CEDHr}):
\[
\ul{C}=
%\underset{{\tiny
%\epsilon \leftarrow \alpha}}{ S  }
\mathbf{S}
\cdot
{\small
\left(
\begin{array}{cc}
\ve & 0 \\ 0& (1-\ve)e^{-\frac{c}{\lambda_{1}-\lambda_{2}}}+\ve
\end{array}
\right)
}
  \cdot
  %\underset{{\tiny
%\alpha \leftarrow \epsilon}}{ S  }\
\mathbf{S}^{-1},
\]
where we note that $\ul{A}=\ul{S}\cdot \textrm{diag}\{\lambda_{1},\lambda_{2}\}\cdot \ul{S}$.

We propose the following specific choices of auxiliary matrices
\[
\begin{array}{l}
\ul{A}(|\mathcal{W}_{\psi}f|)(p,q,s)=  [( \partial_{p_i} \partial_{p_j} G_{\sigma}*|\mathcal{W}_{\psi}f|)(p,q)]_{i,j=1,2}\  , \\[7pt]
\ul{A}(|\mathcal{W}_{\psi}f|)(p,q,s)= [( \partial_{p_i} G_{\sigma}*|\mathcal{W}_{\psi}f|)(p,q) \;  (\partial_{p_j} G_{\sigma}*|\mathcal{W}_{\psi}f|)(p,q)]_{i,j=1,2}\
\end{array}
\]
with $p_{1}=\beta^{2}p \textrm{ and }p_{2}= q$
and separable Gaussian kernel
$
G_{\sigma}(p,q)= \frac{\beta^{2}}{2\pi \sigma^2} e^{-\frac{\beta^4p^2+  q^2}{2\sigma^2}}$.
Note that $\mathcal{A}_{1}|\mathcal{W}_{\psi}f|=\partial_{p}|\mathcal{W}_{\psi}f|$ and $\mathcal{A}_{2}|\mathcal{W}_{\psi}f|=\partial_{q}|\mathcal{W}_{\psi}f|$. Now since $|\mathcal{W}_{\psi}f|=|\mathcal{G}_{\psi}f|$ and $\tilde{\mathcal{A}}_{i}=\mathcal{S} \circ \mathcal{A}_{i} \circ \mathcal{S}^{-1}$ (recall (\ref{phasespacegen})) and $\mathcal{S}\circ \mathcal{W}_{\psi}f= \mathcal{G}_{\psi}f$ we get the following equivalent non-linear left-invariant diffusion equations on phase space:
{\scriptsize
\begin{equation} \label{CEDPS}
\hspace{-0.5cm}
\left\{
\begin{array}{ll}
\partial_{t} \tilde{W}(p,q,t) & = \left(
\left( \begin{array}{cc}
\beta^{-2}\tilde{\mathcal{A}}_{1} &  \tilde{\mathcal{A}}_{2}
\end{array} \right)
 \cdot
%\underset{{\tiny
%\epsilon \leftarrow \alpha}}{ S  }
\mathbf{S}
\cdot
{\small
\left(
\begin{array}{cc}
\ve & 0 \\ 0& (1-\ve)e^{-\frac{c}{\lambda_{1}-\lambda_{2}}}+\ve
\end{array}
\right)
}
  \cdot
  %\underset{{\tiny
\mathbf{S}^{-1}%\alpha \leftarrow \epsilon}}{ S  }
\left(
\begin{array}{c}
\beta^{-2}\tilde{\mathcal{A}}_{1} \tilde{W} \\
 \tilde{\mathcal{A}}_{2} \tilde{W}
\end{array}
\right)
\right)(p,q,t) \\
 & =
\left(
\begin{array}{cc}
 \partial_{\tilde{u}} & \partial_{\tilde{v}}
 \end{array}
\right)
\left(
\begin{array}{cc}
\ve & 0 \\ 0& (1-\ve)e^{-\frac{c}{(\lambda_{1}-\lambda_{2})^2}}+\ve
\end{array}
\right)
\left(
\begin{array}{c}
\partial_{\tilde{u}} \\
\partial_{\tilde{v}}
\end{array}
\right) \tilde{W}(p,q,t), \ \  \\ &\textrm{ for all } (p,q) \in \R^2 , t>0 \textrm{ and fixed } c>0, \ve>0,
 \\[7pt]
\tilde{W}(p,q,0) &=\mathcal{G}_{\psi}f(p,q),  \textrm{ for all } (p,q) \in \R^2.
\end{array}
\right.
\end{equation}
}
with again {\small $\mathbf{S}= (\ul{e}_{1}\; |\; \ul{e}_{2}) $} and $(\partial_{u}W\; \;  \partial_{v}W)=( \beta^{-2}\mathcal{A}_{1}W\; \;\mathcal{A}_{2}W)\, \mathbf{S}$ and where we recall from (\ref{evev}) that
%\[
%\begin{array}{l}
$\lambda_{k}$, %\\
$\ul{e}_{k}$
%\end{array}
%\]
denote the eigenvalues and normalized eigenvectors of the local auxiliary-matrix $\ul{A}(|\mathcal{W}_{\psi}f|)(\ul{x},s)=\ul{A}(|\mathcal{G}_{\psi}f|)(\ul{x},s)$.
In phase space, these eigenvectors (in $\R^2$) correspond to
\[
\{\ul{e}_{1},\ul{e}_{2}\}
\leftrightarrow \{\partial_{\tilde{u}},\partial_{\tilde{v}}\}=
\{\mathcal{S}\circ \partial_{u} \circ \mathcal{S}^{-1}, \mathcal{S} \circ \partial_{v} \circ \mathcal{S}^{-1}\}\ ,
\]
so that we indeed get the right correspondence between (\ref{CEDPS}) and (\ref{CEDHr}):
{\scriptsize
\[
\begin{array}{ll}
\forall_{(p,q,s) \in H_r , t>0}  \; :\;
& \partial_{t} W(p,q,s,t)=
\left(
\begin{array}{cc}
 \partial_{u} & \partial_{v}
 \end{array}
\right)
\left(
\begin{array}{cc}
\ve & 0 \\ 0& (1-\ve)e^{-\frac{c}{(\lambda_{1}-\lambda_{2})^2}}+\ve
\end{array}
\right)
\left(
\begin{array}{c}
\partial_{u} \\
\partial_{v}
\end{array}
\right) W(p,q,s,t)  \desda \\
\forall_{(p,q,s) \in H_r , t>0}  \; :\;
& \partial_{t} (\mathcal{S}^{-1}\tilde{W})(p,q,s,t)=
\left(
\begin{array}{cc}
 \partial_{u} & \partial_{v}
 \end{array}
\right)
\left(
\begin{array}{cc}
\ve & 0 \\ 0& (1-\ve)e^{-\frac{c}{(\lambda_{1}-\lambda_{2})^2}}+\ve
\end{array}
\right)
\left(
\begin{array}{c}
\partial_{u} \\
\partial_{v}
\end{array}
\right) (\mathcal{S}^{-1}\tilde{W})(p,q,s,t) \ \ \desda \\
\forall_{(p,q,s) \in H_{r}, t>0}  \; :\; &
 \partial_{t} (\mathcal{S}^{-1}\tilde{W})(p,q,s,t)=
\left(
\begin{array}{cc}
 \partial_{u} & \partial_{v}
 \end{array}
\right) \mathcal{S}^{-1}
\left(
\begin{array}{cc}
\ve & 0 \\ 0& (1-\ve)e^{-\frac{c}{(\lambda_{1}-\lambda_{2})^2}}+\ve
\end{array}
\right) \mathcal{S}
\left(
\begin{array}{c}
\partial_{u} \\
\partial_{v}
\end{array}
\right) (\mathcal{S}^{-1}\tilde{W})(p,q,s,t) \desda \\
\forall_{(p,q) \in \R^2, t>0}  \; :\;
&\partial_{t} \tilde{W}(p,q,t)=
\left(
\begin{array}{cc}
 \partial_{\tilde{u}} & \partial_{\tilde{v}}
 \end{array}
\right)
\left(
\begin{array}{cc}
\ve & 0 \\ 0& (1-\ve)e^{-\frac{c}{(\lambda_{1}-\lambda_{2})^2}}+\ve
\end{array}
\right)
\left(
\begin{array}{c}
\partial_{\tilde{u}} \\
\partial_{\tilde{v}}
\end{array}
\right) \tilde{W}(p,q,t)
%\forall_{p,q \in \R, t>0}  \, :\,
%\partial_{t} (\mathcal{S}^{-1}W)(p,q,t)=
%\mathcal{S}
%\left(
%\begin{array}{cc}
% \partial_{u} & \partial_{v}
% \end{array}
%\right)
%\mathcal{S}^{-1}
%\left(
%\begin{array}{cc}
%\ve & 0 \\ 0& (1-\ve)e^{-\frac{c}{(\lambda_{1}-\lambda_{2})^2}}+\ve
%\end{array}
%\right)\mathcal{S}
%\left(
%\begin{array}{c}
%\partial_{u} \\
%\partial_{v}
%\end{array}
%\right) (S^{-1} \tilde{W})(p,q,t)
\end{array}
\]
}
so that the uniqueness of the solutions $\tilde{W}$ of (\ref{CEDPS}) and $W$ of (\ref{CEDHr}) implies that
\[
\tilde{W}(\cdot,\cdot,0)=\mathcal{G}_{\psi}f=\mathcal{S} \circ \mathcal{W}_{\psi}f =\mathcal{S} \circ W(\cdot,\cdot,\cdot,0) \Rightarrow
\forall_{t \geq 0} : \tilde{W}(\cdot,\cdot,t)=\mathcal{S} \circ W(\cdot,\cdot,\cdot,t)\ .
\]
Clearly, Problem (\ref{CEDPS}) is preferable over problem (\ref{CEDHr}) if it comes to numerical schemes as it is a 2D-evolution.

\begin{figure}
\centerline{\includegraphics[width=0.8\hsize]{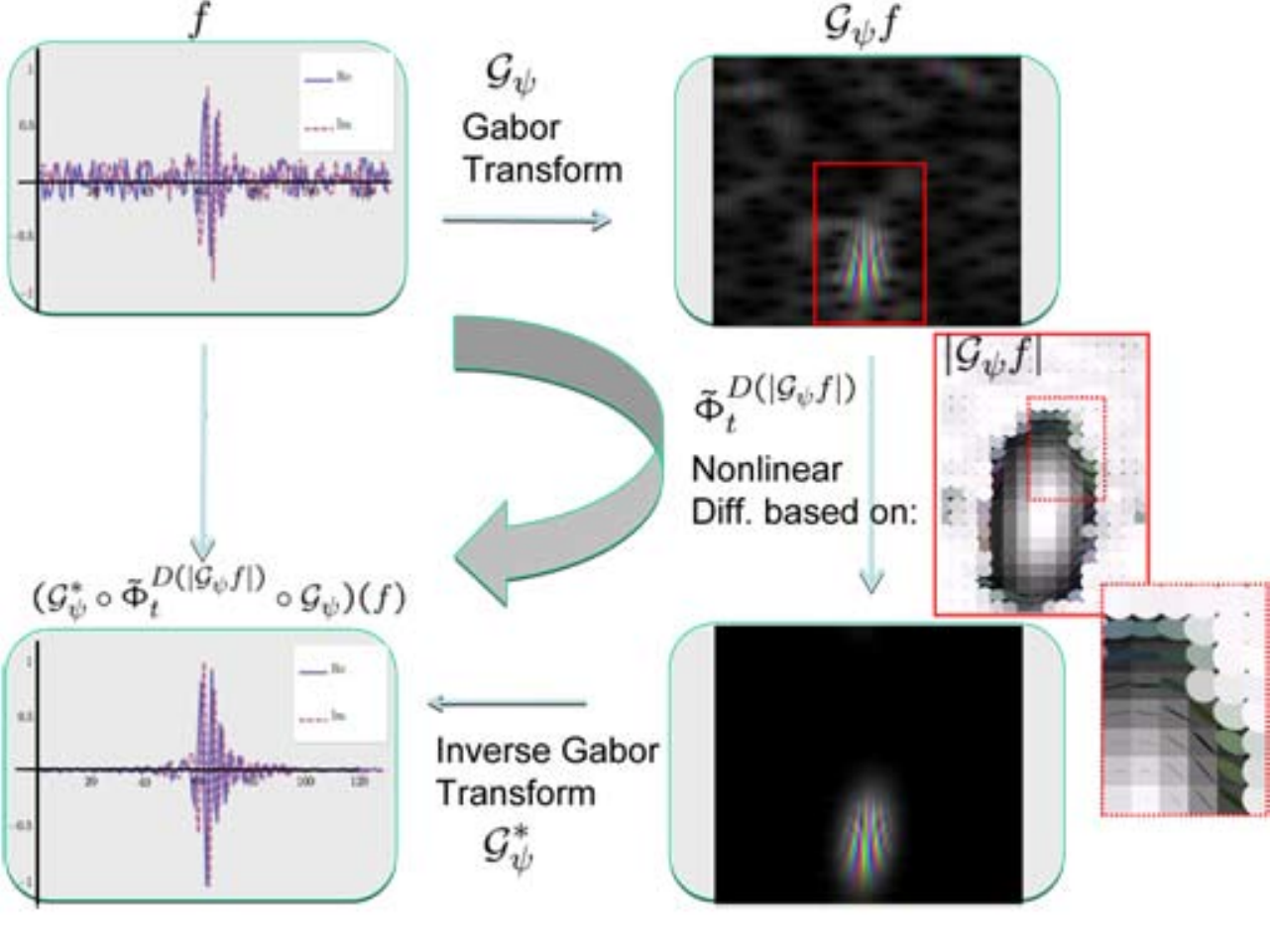}
}
\caption{Illustration of the left-invariant diffusions on Gabor transforms that are adapted to the Hessian of the
absolute value. The corresponding operator $f \mapsto \mathcal{G}_{\psi}^{*} \tilde{\Phi}_{t}^{D} \mathcal{G}_{\psi}$ in the signal domain smoothly enhances the signal without tresholding of Gabor coefficients.
In the most right plot we depicted ellipsoids representing the local eigenvectors of the Hessian of $|\mathcal{G}_{\psi}|$ along which the diffusion locally takes place. The directions of the ellipsoids coincide with the eigenvectors of the Hessian, whereas the anisotropy of the ellipsoids is determined by the fraction $|\lambda_{1}/\lambda_{2}|$ of the eigenvalues $\{\lambda_{1},\lambda_{2}\}$ with $|\lambda_{1}|>|\lambda_{2}|$, the colors indicate the directions of the largest eigenvector and the intensity reflects the relative damping factor $(p,q) \mapsto 1+ \epsilon^{-1} (1-\epsilon)e^{-c (\lambda_{1}(p,q)-\lambda_{2}(p,q))^{-2}} \lambda_{1}(p,q)$ of the second eigenvector.}\label{fig:ellips}
\end{figure}

See Figure \ref{fig:ellips} for an explicit example, where we applied adaptive nonlinear diffusion on the Gabor transform of a noisy chirp signal. Due to left-invariance of our diffusions the phase is treated appropriately and thereby the effective corresponding denoising operator on the signal is sensible. Moreover, the tails of the enhanced chirp are smoothly dampened.
%We implemented the diffusions (\ref{CEDPS}) by a standard explicit Euler-Forward scheme. Figure %\ref{fig:translation} shows a promising result of adaptive phase-covariant diffusions on an explicit example.
%For illustration of the local adaptivity of the phase-covariant diffusions in phase space (\ref{CEDPS}) within %this specific example, see Figure \ref{fig:ellips}.

\section{Reassignment, Texture Enhancement and Frequency Estimation in 2D-images ($d=2$) \label{ch:2D}}

\begin{figure}
\centerline{
\includegraphics[width=0.25\hsize]{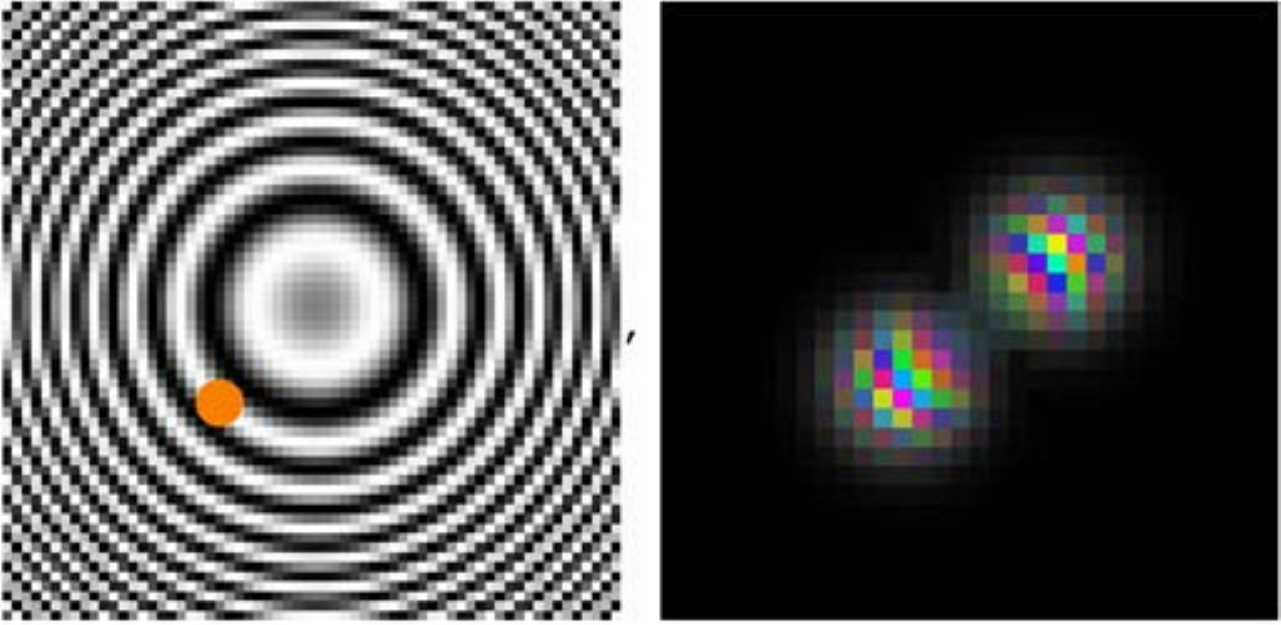}
\includegraphics[width=0.25\hsize]{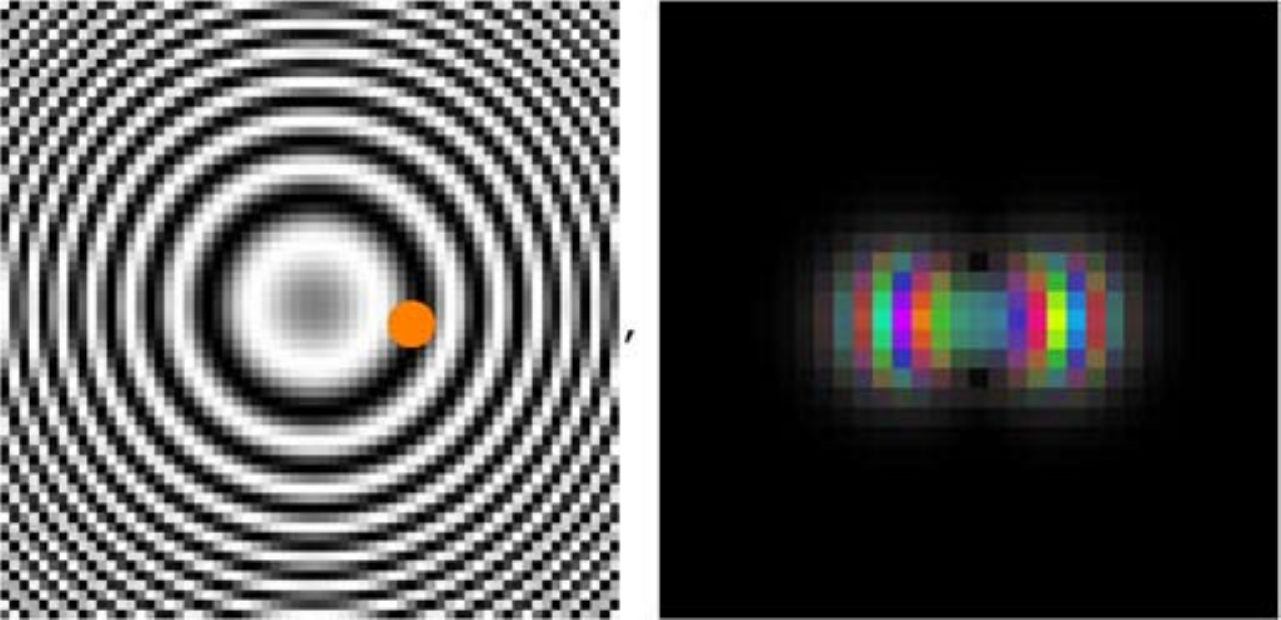}
}
\centerline{
\includegraphics[width=0.25\hsize]{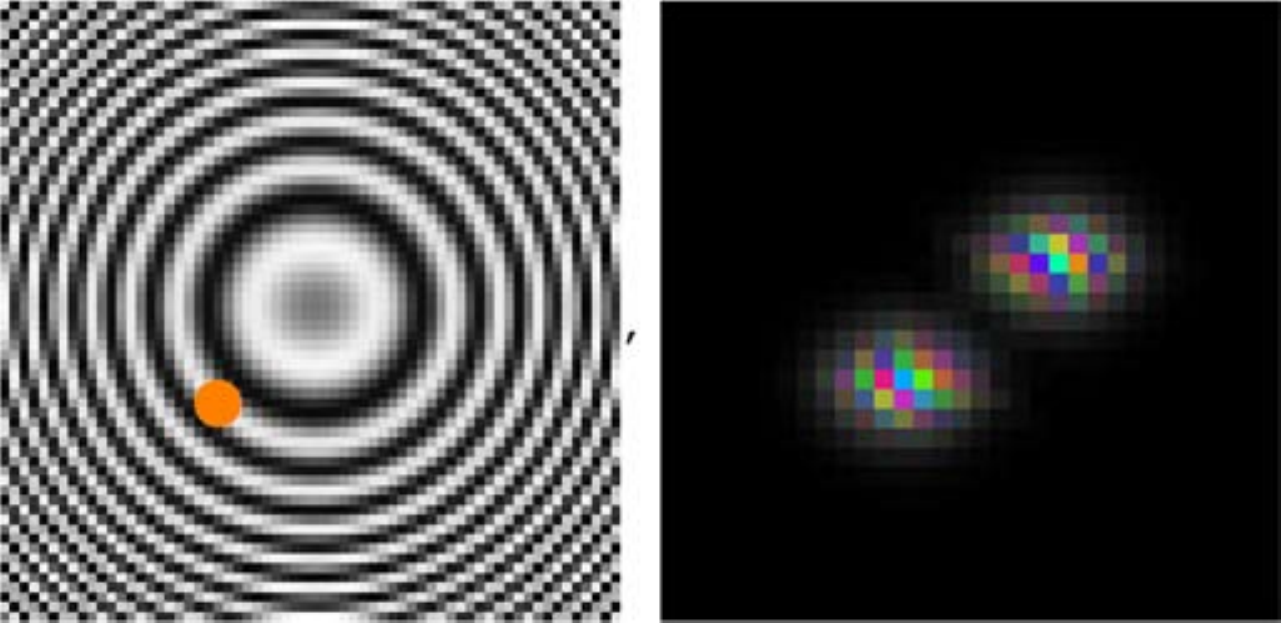}
\includegraphics[width=0.25\hsize]{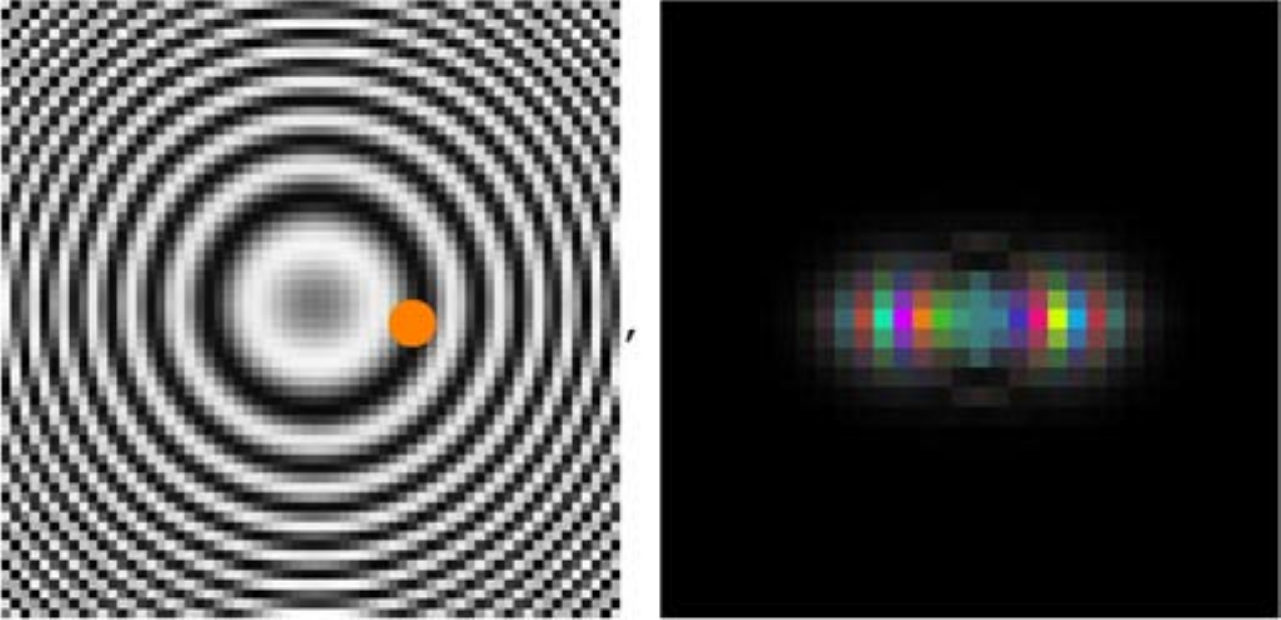}
}
\centerline{
\includegraphics[width=0.25\hsize]{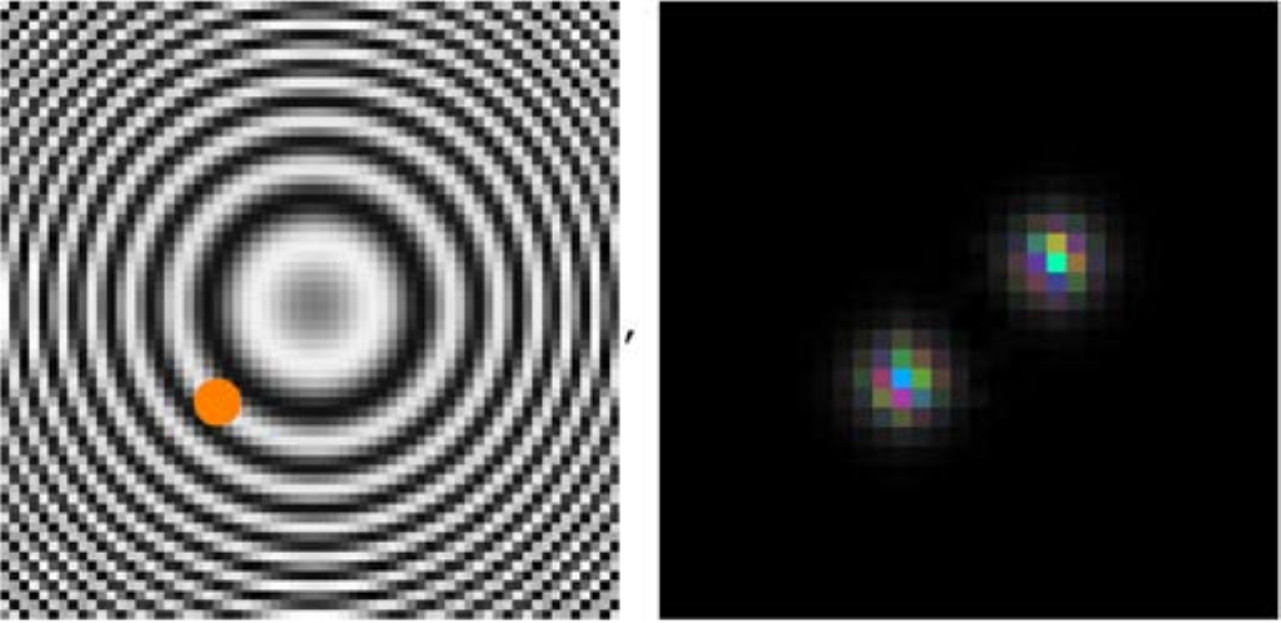}
\includegraphics[width=0.25\hsize]{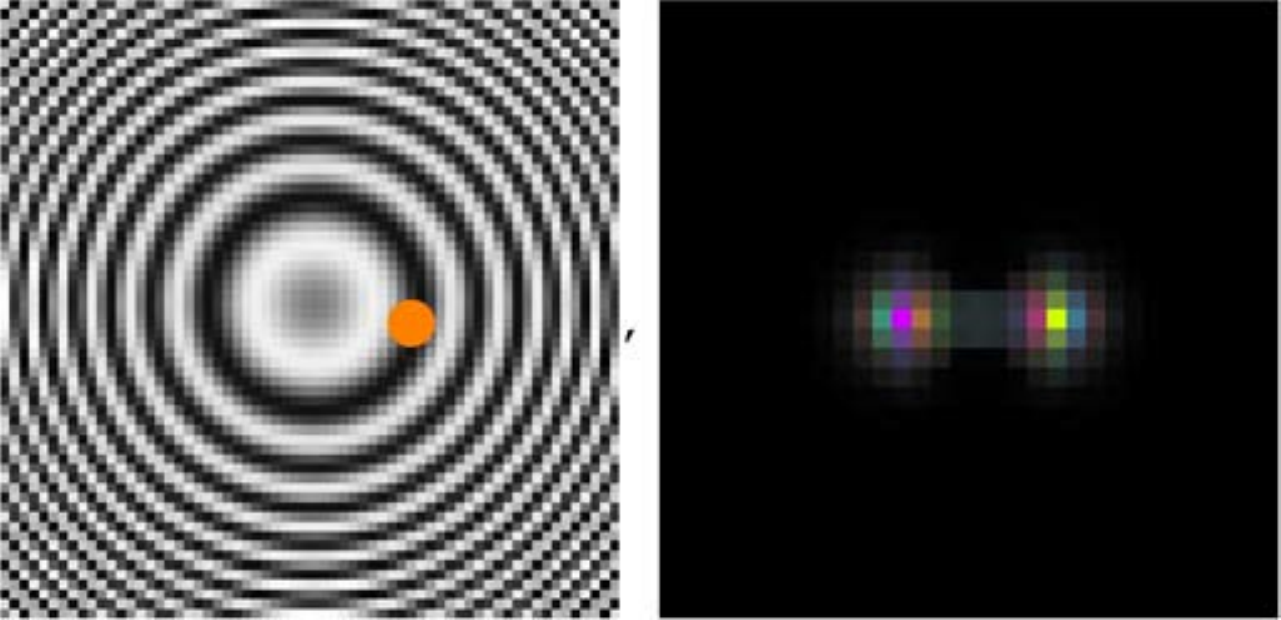}
}
\caption{Application of re-assignment of Gabor transforms of 2D-images. The dot
indicates the position $(p^{1},p^{2}) \in \R^{2}$ in the image (left)
where we depicted the local frequency distribution
$(q^{1},q^{2})\mapsto \mathcal{G}_{\psi}f(p^{1},p^{2},q^{1},q^{2})$ (right)
where we color-coded the phase. Top row input at two locations. Middle row output
for $t=0.1$. Bottom row output for $t=4$. Further parameter settings: $a=1$, $N=64$, $K=M=32$, $L=2$,
$D^{11}=D^{22}=1$.
}\label{fig:reass2D}
\end{figure}

%Short experimental section To be written...
In Section \ref{ch:reass} we applied differential reassignment to 1D-signals.
This technique can also be applied to Gabor transform of 2D-images. Here
we choose for the second option in (\ref{conchoice}) as the corresponding algorithm is faster.
In this case the reassigned Gabor transforms concentrate towards the lines
$\{(p^{1},p^{2},q_{1}, q_{2},s) \in H_{5} \; |\; \partial_{p^{i}} |\mathcal{G}_{\psi}f|(p,q,s)=\partial_{q_{i}} |\mathcal{G}_{\psi}f|(p,q,s)=0, i=1,2, p=(p^{1},p^{2}), q=(q_{1},q_{2}) \in \R^{2}\}$, which
coincides with the stationary solutions of the corresponding
Hamilton Jacobi equation on the modulus. See Figure \ref{fig:reass2D}, where the amount of large
Gabor coefficients is strongly reduced while maintaining the original image.
%One of the problems with differential reassignment

\subsection{Left-invariant Diffusion on Phase Space as a Pre-processing Step for Differential Reassignment}

Gabor transforms of noisy medical images often require smoothing as pre-processing before differential reassignment can be applied.
Linear left-invariant diffusion is a good choice for such a smoothing
as a pre-processing step for differential reassignment and/or local frequency extraction, since such smoothing does not affect the original signal (up to a scalar multiplication).
Here we aim for pre-processing before reassignment rather than signal enhancement.
%\begin{lemma}
%Let $\Phi:\mathcal{H}_{n} \to \mathcal{H}_{n}$ be a closed, linear and left-invariant operator in the Gabor domain, with $n=1$, and let $\tilde{\Phi}= \mathcal{S} \circ \Phi \circ \mathcal{S}^{1}$ be the corresponding operator on phase space. If $\tilde{\Phi}$ commutes with $\tilde{\mathcal{L}}_{g}$ given by
%\begin{equation}\label{Lgtilde}
%(\tilde{\mathcal{L}}_{(p,q,s)}\tilde{U})(p',q')=e^{2\pi i s}e^{\pi i (2 q p'-pq)}\tilde{U}(p'-p,q'-q),
%\end{equation}
%for all $\tilde{U} \in \mathbb{L}_{2}(\R^{2d})$ and all $(p,q,s) \in H_{r}$ and almost every $(p',q') \in \R^{2d}$ then
%the corresponding operator in the signal domain
%$
%\Upsilon_{\psi} = \mathcal{W}_{\psi}^{*} \Phi \mathcal{W}_{\psi} =  \mathcal{G}_{\psi}^{*} \tilde{\Phi} \mathcal{G}_{\psi}$
%is a scalar multiple of the identity.
%\end{lemma}
%\textbf{Proof }
%Recall from Theorem~\ref{th:one} that $\tilde{\mathcal{L}}_{(p,q,s)}= \mathcal{S} \mathcal{L}_{(p,q,s)} \mathcal{S}^{-1}$,
%so by means of
%Eq.~\!(\ref{S}) and Eq.~\!(\ref{Smin1}) and direct computation one obtains (\ref{Lgtilde}). Now $\Phi$ is linear iff $\tilde{\Phi}$ is linear and $\Phi$ commutes with
%$\mathcal{L}_{g}$ iff $\tilde{\Phi}$ commutes with $\tilde{L}_{g}$ which implies that $\Upsilon_{\psi}$ is linear and commutes with $\mathcal{U}^{n}_{g}$ for all $g \in H_{r}$. Now $\Upsilon_{\psi}$ is a closed linear operator that commutes with the irreducible representation
%$\mathcal{U}^{n}$ of $H_{r}$, cf.~\cite[Thm 9.2.1]{Grochenic}, and by an extension of Schur's lemma \cite{Dieudonne} it must be a multiple times the identity. $\hfill \Box$ \\
%\\
The Dunford-Pettis Theorem \cite{Bukhvalov} shows minor conditions for a linear operator on $\mathbb{L}_{2}(X)$, with $X$ a measurable space,
to be a kernel operator.
Furthermore, all linear left-invariant kernel operators on $\mathbb{L}_{2}(H_{2d+1})$ are convolution operators.
Next we classify all left-invariant kernel operators on phase space (i.e. all operators that commute with $\tilde{\mathcal{L}}_{(p,q,s)}:= \mathcal{S} \mathcal{L}_{(p,q,s)} \mathcal{S}^{-1}$).
\begin{lemma}
A kernel operator
$\mathcal{K}_{k}: \mathbb{L}_{2}(\R^{2d}) \to \mathbb{L}_{\infty}(\R^{2d}) \cap \mathbb{L}_{2}(\R^{2d})$ given by
\[
(\mathcal{K}_{k}\tilde{U})(p,q)= \int_{\R^{2d}} \tilde{k}(p,q,p',q') \tilde{U}(p',q') {\rm d}p'{\rm d}q',
\]
with $k \in \mathbb{L}_{1}(\R^{2d} \times \R^{2d})$
is left-invariant if
\begin{equation}\label{lips}
\tilde{k}(p,q,p',q')= e^{2\pi i \overline{q} (p'-p)} \tilde{k}(p+\overline{p},q+\overline{q}, p'+\overline{p},q' +\overline{q} )
\end{equation}
for almost every $(p,q), (p',q'), (\overline{p},\overline{q}) \in \R^{2d}$.
\end{lemma}
\textbf{Proof }
We have $\tilde{\mathcal{L}}_{g}\tilde{U}(p',q')=(\mathcal{S}\tilde{\mathcal{L}}_{g}\mathcal{S}^{-1}\tilde{U})(p',q')=e^{2\pi i s + \pi i (2qp'-pq)} \tilde{U}(p'-p,q'-q)$ so that
\[
\begin{array}{l}
(\mathcal{K}_{k} \mathcal{L}_{g}\tilde{U})(\tilde{p},\tilde{q})= \int_{\R^{2d}} \tilde{k}(\tilde{p},\tilde{q},p',q') e^{2\pi i s+ \pi i(2qp'-pq)} \tilde{U}(p'-p,q'-q) {\rm d}p'{\rm d}q' \\
=
e^{2\pi i s }\int_{\R^{2d}} e^{\pi i (2q\tilde{p}-pq)} \tilde{k}(\tilde{p}-p,\tilde{q}-q,p'',q'') \tilde{U}(p'',q'') {\rm d}p'' {\rm d}q''
= ( \mathcal{L}_{g}\mathcal{K}_{k}\tilde{U})(\tilde{p},\tilde{q})
\end{array}
\]
from which the result follows by substitution $p''=p'-p$, $q''=q'-q$.$\hfill \Box$
\begin{theorem} \label{th:diffusionphasespace}
Let $d \in \mathbb{N}$.
Let $\Phi_{t}$ be a left-invariant semi-group operator on $H_{r}$ given by
\[
(\Phi_{t}(U))(g)=(K_{t}*U)(g) := \int \limits_{H_{r}} K_{t}(h^{-1}g) U(h) \,{\rm d}\mu(h)
\]
with $\mu$ the left-invariant measure on $H_{r}$, $t>0$.
Then the corresponding operator on phase space $\tilde{\Phi}_{t} =\mathcal{S} \circ \Phi_{t} \circ \mathcal{S}^{-1}$ is given by
\begin{equation} \label{pso}
\begin{array}{l}
(\tilde{\Phi}_{t}(\mathcal{G}_{\psi}f))(p,q)=
\int \limits_{\R^{d}} \int \limits_{\R^{d}} \tilde{k}_{t}(p,q,p',q') \mathcal{G}_{\psi}f(p',q')\,
{\rm d}p' {\rm d}q', \\
\textrm{with }\tilde{k}_{t}(p,q,p',q')= \int \limits_{0}^{1} e^{-\pi i p' \cdot q'} K_{t}(p-p',q-q',\frac{-p \cdot q}{2}-s'-\frac{1}{2}
(p \cdot q' \!-\!q \cdot p')) e^{-2\pi i s'}\, {\rm d}s'.
\end{array}
\end{equation}
In particular we consider the diffusion kernel for horizontal diffusion generated by $D^{11}\sum \limits_{i=1}^{d}(\mathcal{A}_{i})^2 + D^{22}\sum \limits_{i=d+1}^{2d}(\mathcal{A}_{i})^2$ on the sub-Riemannian manifold $(H_{2d+1},{\rm d}s + \frac{1}{2}(p \cdot {\rm d}q -q \cdot {\rm d}p))$.
\end{theorem}
\textbf{Proof }Eq. (\ref{pso}) follows by direct computation where we note $h^{-1}=(p',q',s')^{-1}=(-p',-q',-s')$ and taking the phase inside, which is possible since $(\mathcal{S} \mathcal{W}_{\psi}f)(p',q')=(\mathcal{G}_{\psi}f)(p',q')$ is independent of $s'$.
The heat kernel on $H_{r}$ is obtained by the heat kernel on $H_{2d+1}$ via
$
K_{t}^{H_{r}}(p,q,s)= \sum \limits_{k \in \mathbb{Z}} K_{t}^{H_{3}}(p,q,s+k)$. $\hfill \Box$ \\
\\
As
the analytic closed form solution of this kernel is intangible, cf. \cite{Gaveau,DuitsR2006SS2},
we consider the well-known local analytic approximation
{\small
\begin{equation} \label{localapprox}
K_{t}(p,q,s)= \frac{c}{8 t \sqrt{D^{11}D^{22}}} \frac{1}{(4\pi t c^{-1}D^{11})^{\frac{d}{2}}(4\pi t c^{-1}D^{22})^{\frac{d}{2}}} e^{-\left( \frac{|s|}{\sqrt{D^{11}D^{22}}}+\frac{\|p\|^{2}}{D^{11}} + \frac{\|q\|^{2}}{D^{22}}\right) \frac{c}{4t}}.
\end{equation}
}
This approximation is due to the ball-box theorem \cite{Bellaiche} or theory of weighted sub-coercive operators on Lie groups \cite{TerElst3}, where we assign weights $w_{i}=1$, $i=1\ldots,2d$ and $w_{2d+1}=2$ to the left-invariant vector fields
$\{\mathcal{A}_{1},\ldots,\mathcal{A}_{2d+1}\}$.
\begin{lemma}
Let $K_{t}:H(2d+1) \to \R^{+}$ be given by Eq.~(\ref{localapprox}).
Then the corresponding kernel on phase space is given by
{\small
\begin{equation} \label{resultforpractice}
\tilde{k}_{t}(p,q,p',q')= e^{-\left(\frac{\|p-p'\|^2}{D^{11}}+\frac{\|q-q'\|^2}{D^{11}} \right)\frac{c}{4t}} e^{-\pi i (p'-p)\cdot (q'+q)} \frac{\left( \frac{c}{4\pi t \sqrt{D^{11}D^{22}}}\right)^{d}}{1+ 64 D^{11}D^{22} t^2 c^{-2}} .
\end{equation}
}
and satisfies Eq.~(\ref{lips}).
\end{lemma}
\textbf{Proof }
By Eq.~(\ref{localapprox}) and Eq.~(\ref{pso}) we find by substitution $v=s'-s-k+\frac{1}{2}(p \cdot q' - q\cdot p')$ that
\[
\begin{array}{l}
\tilde{k}_{t}(p,q,p',q')= e^{\pi i (p \cdot q' -q \cdot p')} e^{-2\pi s} e^{-i p'\cdot q'}  \frac{c}{8 t \sqrt{D^{11}D^{22}}} \frac{1}{(4\pi t c^{-1}D^{11})^{\frac{d}{2}}(4\pi t c^{-1}D^{22})^{\frac{d}{2}}} \cdot \\
 \sum \limits_{k \in \mathbb{Z}} \int \limits_{s'-k+\frac{1}{2}(p \cdot q' - q \cdot p')}^{1+s'-k+\frac{1}{2}(p \cdot q' - q \cdot p')} e^{-|v| \frac{c}{4t \sqrt{D^{11}D^{22}}}} e^{-2\pi v i}\, {\rm d}v.
 \end{array}
\]
The integral can be computed by means $\int_{\R} e^{-|v| \gamma} e^{-2\pi i v}\, {\rm d}v=\frac{2\gamma}{\gamma^2+4\pi^2} $ for
$\gamma>0$ . The kernel in Eq.\!(\ref{resultforpractice}) satisfies the left-invariance constraint
(\ref{lips}) since \mbox{$-(\tilde{p}-\overline{p})(\tilde{q}+\overline{q})=2q(\tilde{p}-\overline{p})-(\tilde{p}-\overline{p})(\tilde{q}+\overline{q}+2q)$}
$\hfill \Box$

\subsubsection{Possible Extension to Texture Enhancement in 2D images via Left-invariant Evolutions}

\begin{figure}
\centerline{
\includegraphics[width=\hsize]{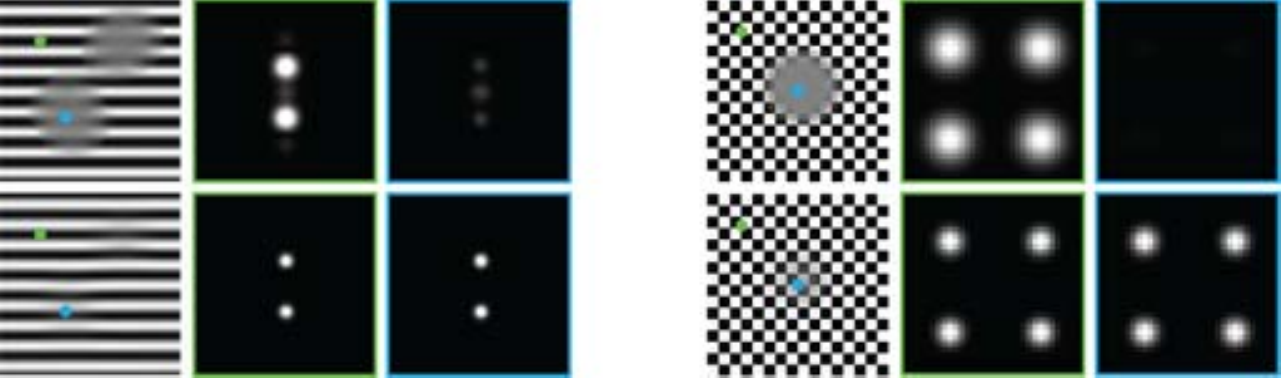}
}
\caption{Top row: Input image $f_{i}$ and two slices $|\mathcal{G}_{\psi}f_i(\ul{p}_{k},\cdot,\cdot)|$, $k=1,2$, in the Gabor domain centered around the indicated positions. Bottom row: output in spatial and Gabor domain
of linear diffusion combined with a squaring of the modulus. (Left: input image $f_{1}$, i.e. $i=1$, right: input image $f_{2}$, i.e. $i=2$). Akin to left-invariant diffusions on orientation scores where we can generically diffuse along crossing lines (Recall Figure \ref{fig:OS}) we can apply left-invariant diffusion on Gabor transforms and thereby diffuse along crossing frequencies.
}\label{fig:basic}
\end{figure}

The techniques signal enhancement by non-linear left-invariant diffusion on Gabor transforms in subsection~\ref{ch:covdiff} can be extended to enhancement of local 2D-frequency patterns and/or textures. This can have similar applications as the enhancement of lines via non-linear diffusion on invertible orientation scores, cf.~\cite{Fran2009,DuitsAMS2} and Figure~\ref{fig:OS}. However, such an extension would yield a technical and slow algorithm. Instead, akin to our earlier works on contour
enhancement via invertible orientation scores, cf.~\cite{Duits2005IJCV,DuitsAMS1}, we can use the (more basic) concatenation of linear left-invariant diffusion and monotonic transformations in the co-domain. Figure \ref{fig:basic} shows a basic experiment of such an approach.

\subsection{Local Frequency Estimation in Cardiac Tagged MRI Images \label{ch:pract}}
%\subsection{Left-invariant diffusion on phase space for enhancement}
%
%Short experimental section. To be written.
In the limiting case differential reassignment concentrates around local maxima of the absolute value of a Gabor transform. These local maxima produce per position $(p^{1},p^{2}) \in \R^{2}$ a local frequency $(q_1,q_2) \in \R^{2}$ estimation. The reliability of such local frequency estimations can be checked via differential reassignment. %, e.g. the larger the time where differential reassignment still preserves the image
%(up to an a priori set tolerance)
%the more reliable the corresponding frequency estimates.
In this section we briefly explain a 2D medical imaging application where local frequency estimation is important for measuring heart wall deformations during a systole. For more details, see \cite{Bruurmijnthesis}.

Quantification of cardiac wall motion may help in (early) diagnosis of cardiac abnormalities such as ischemia and myocardial infarction. To characterize the dynamic behavior of the cardiac muscle, non-invasive acquisition techniques such as MRI tagging can be applied. This allows to locally imprint brightness patterns in the muscle, which deform accordingly and allow detailed assessment of myocardial motion. Several optical flow techniques have been considered in this application. However, as the constant brightness assumption \cite{HandSpaper,BruhnetAl,Flor98a} in this application is invalid, these techniques end up in a concatenation of technical procedures,
\cite{Duits2011QAM,BecciuThesis,Florackharp}
yielding a complicated overall model and algorithm. For example, instead of tracking constant brightness one can track local extrema in scale space \cite{Jansssen2006b} (taking into account covariant derivatives and Helmholtz decomposition of the flow field \cite{Duits2011QAM}) or one can compute optical flow fields \cite{Florackharp,vanAssen2007CBfMII} that follow equiphase lines, where the phase is computed by Gabor filtering techniques known as the harmonic phase method (HARP), cf.~\!\cite{Osman1999}.
Gabor filtering techniques were also used in a recent applied approach \cite{Arts} where one obtains cardiac wall deformations directly from local scalar-valued frequency estimations in tagging directions.
Here we aim for a similar short-cut, but in contrast to the approach in \cite{Arts} we extract the maxima from all Gabor transform coefficients producing per tagging direction an accurate frequency covector field
$f=f_{1} {\rm d}x^{1}+f_{2} {\rm d}x^{2}$ (not necessarily aligned with the tagging direction).
From these covector fields one can deduce the required deformation gradient $D$ by duality as we briefly explain next.

\begin{figure}
\centerline{
\includegraphics[width=0.8\hsize]{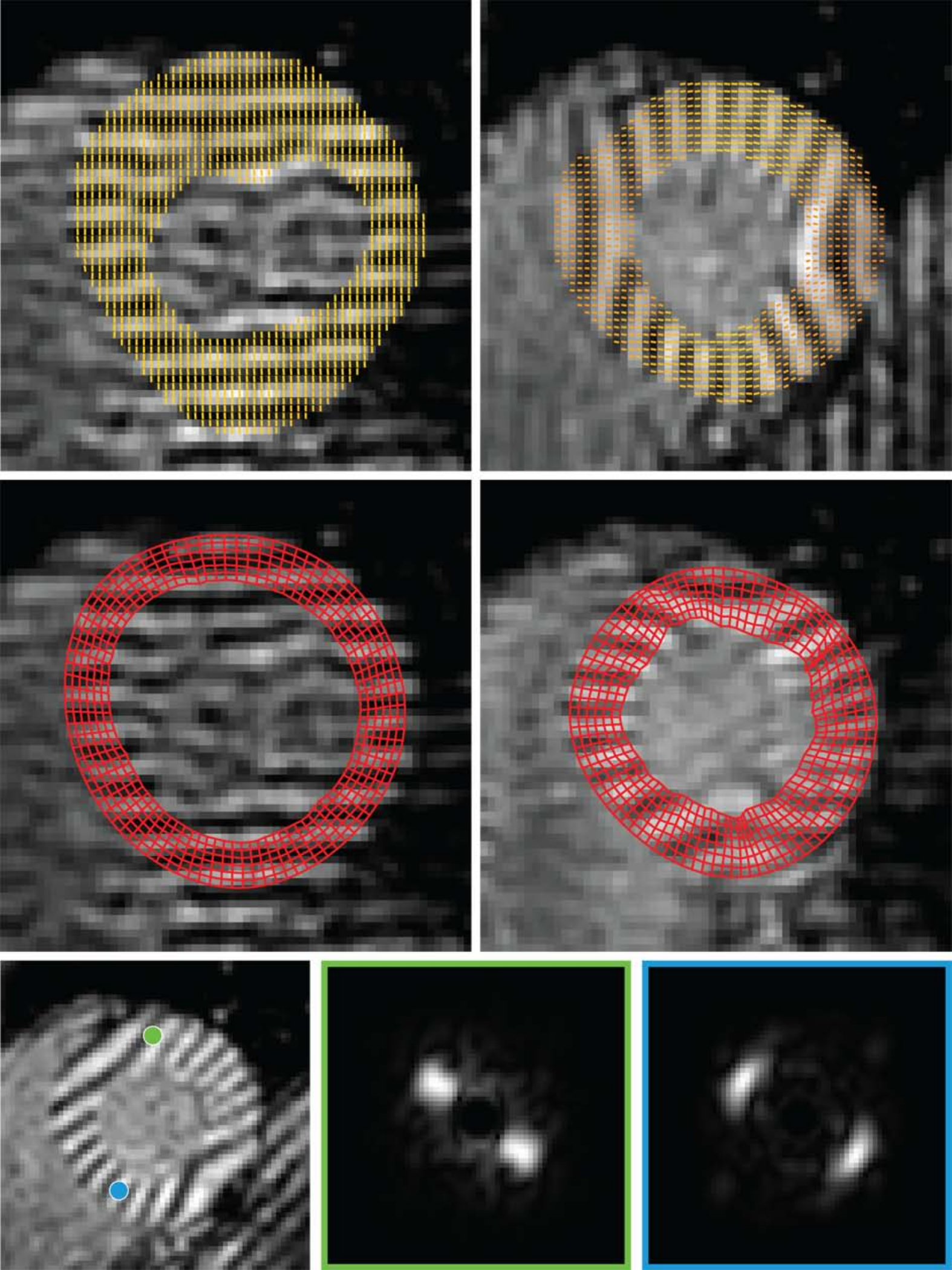}
}
%\end
\caption{Top row: MRI-tagging images (on a $64\times 64$-grid) at time frame $t=2$ and time frame $t=8$ (the systolic phase of the heart is sampled in 30 time frames), scanned in 4 tag-directions, i.e. $N=4$. Left: horizontal tag $\theta_{1}=0$ at $t=2$ with
estimated frequency field $\ul{q}^{2,1}$, right: vertical tag $\theta_{3}=\pi/2$ with estimated frequency field $\ul{q}^{8,3}$ plotted on top. In the frequency estimates we used maximum over-sampling $L=1$, $K=64$ in the discrete Gabor transform Eq.~(\ref{discretegabor}) with Gaussian window-size with a standard deviation of $\sigma=6$ (pixel lengths).
Middle row: Corresponding deformation nets computed from the frequency fields and outer boundary via Eq.~(\ref{Dtx}) and Eq.~(\ref{defnet})
Bottom row: Plots of the absolute value of the Gabor transform $(q_{1},q_{2}) \mapsto |\mathcal{G}_{\psi}f(p^{1},p^{2}, q_{1},q_{2})|$ of MRI-tagging image $f$ with orientation $\theta_{2}=\pi/4$ restricted to the green and blue point $(p^{1},p^{2})$ of interest. We have masked the low frequencies to show the 2 frequency blobs of the tagging lines.}\label{fig:tagswithGT}
\end{figure}

Let $f_{t}^{\theta_{i}}:\R^{2} \to \R$ be the images corresponding to tag-direction $\theta_i \in [0,2\pi)$, $i=1,\ldots, N$, with $N \geq 2$ the
total number of tags at time $t>0$ corresponding to (independent) tag-direction $(\cos(\theta_i),\sin (\theta_{i}))$ with $\theta_i \in [0,2\pi)$.
We compute the Gabor transforms $\mathcal{G}_{\psi}(f_{t}^{\theta_{i}}):\R^{4} \to \mathbb{C}$ and apply a linear left-invariant evolution and extract per position $\ul{p}$ the remaining maxima w.r.t. the frequency variable $\ul{q}$. For each
position $\ul{p}=(p^{1},p^{2})\in \R^{2}$ the Gabor transform $\mathcal{G}_{\psi}(f_{t}^{\theta_{i}})(\ul{p},\cdot)$ typically shows only two
dominant and noisy blobs (for $t$ small), cf.~\!Figure \ref{fig:tagswithGT}.
The 2 blobs relate to each other by reflection $\ul{q} \mapsto -\ul{q}$.

\begin{figure}
\centerline{\includegraphics[width=\hsize]{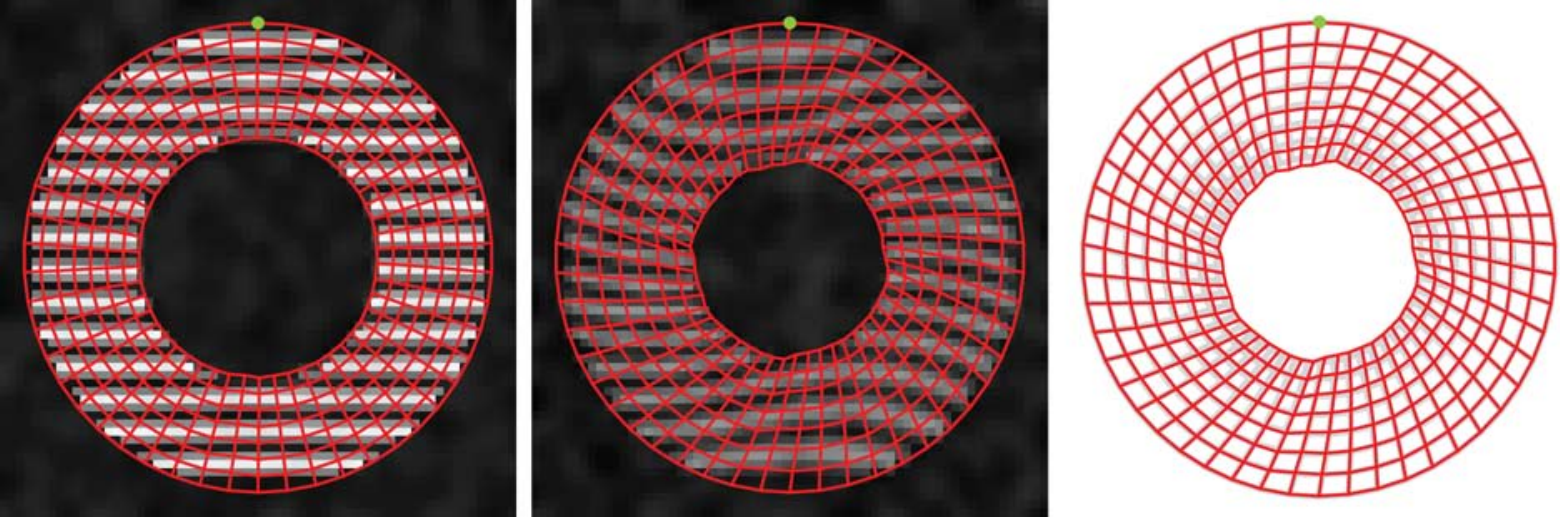}}
\caption{Ground truth phantom where the middle image at time $t>0$ is obtained from the left image $t=0$ via both non-uniform scaling and non-uniform rotation. % (both relatively large compared to two subsequent cardiac MRI-tags)
The phantom includes a strong fading of the tag-line intensities in this phantom. The deformation net depicted on top in red is computed according to Eq.~(\ref{defnet}), with $N=4$, $R=8$ and $J=50$. The deformation net is entirely computed from a single indicated material point at the boundary. The right image shows a comparison between
the computed deformation net (red) with the ground-truth deformation net (grey).
%as follows:
%After detecting the inner heart-wall boundary, we iteratively construct the next outward radial closed curve via the estimated %deformation gradient obtained by the local frequency field (as in Figure \ref{fig:tagswithGT}) via %\mbox{$\ul{x}_{t}^{R+1,i}=\ul{x}_{t}^{R,i}+ \ul{D}_{t}(\ul{x}_{0}^{R+1,i}) \;
%(\ul{x}_{0}^{R+1,i}-\ul{x}_{0}^{R,i})$}, where $R$ is a radial index and $i$ enumerates points on the $R$-th closed curve.
%This experiment shows that our key assumption Eq.\!~(\ref{111}) produces valid deformation estimates.
}\label{fig:exp}
\end{figure}

Best frequency estimates
are obtained by extracting maxima after applying a large time evolution. By the results in \cite{Loog2001a} this maximum converges to the center of mass, that
gives a sub-pixel accurate estimate. This yields per position $(p^{1},p^{2}) \in \R^{2}$
and per tagging direction $\theta_{i}$ the local frequency estimate $\ul{q}^{t,i}(p^{1},p^{2}):=\pm (q_{1}^{t,i}(p^{1},p^{2}),q_{2}^{t,i}(p^{1},p^{2}))$ that we store in a covector field
\[
\ul{q}^{t,i}=q_{1}^{t,i}{\rm d}x^{1}+q_{2}^{t,i}{\rm d}x^{2},
\]
where ${\rm d}x^{1}={\rm d}x$,
${\rm d}x^{2}={\rm d}y$. %\\
%\textbf{Question: How do we fix the signs of this field ?}
%\\
The obtained frequency field $\ul{q}^{t,i}$ can be related to the requested deformation gradient $\ul{D}_t=[\frac{\partial \ul{x}_t}{\partial \ul{x}_{t-1}}]$, where $\ul{x}_t$ denotes the position of a material point in the heart wall at time $t>0$. We assume duality between
frequencies and velocities by imposing
\begin{equation} \label{111}
% \left.\frac{d}{ds} \langle \ul{q}^{t,i}(\ul{x}_t),\ul{x}_{t}(s) \rangle \right|_{s=0}=\langle \ul{q}^{t,i}(\ul{x}_t)=
% \left.\frac{d}{ds} \langle \ul{q}^{t-1,i}(\ul{x}_{t-1}),\ul{x}_{t-1}(s) \rangle \right|_{s=0}=
\langle \ul{q}^{t,i}(\ul{x}_{t}),
 \left.\frac{d}{ds} \ul{x}_{t}(s)\right|_{s=0} \rangle = \langle \ul{q}^{t-1,i}(\ul{x}_{t-1}), \left.\frac{d}{ds} \ul{x}_{t-1}(s)\right|_{s=0} \rangle
\end{equation}
for all smooth parameterizations $(s,t) \mapsto \ul{x}_{t}(s) \in \R^{2}$ with $\ul{x}_{t}(0)=\ul{x}_t$ and all $i=1,\ldots, N$.
This yields
\[
\begin{array}{c}
\sum \limits_{k,j=1}^{2}q_{j}^{t,i}(\ul{x}_{t}) D^{j}_{k}(\ul{x}_{t-1})\; \dot{x}^{k}_{t-1}(0)
=\sum \limits_{j=1}^{2}q_{j}^{t,i}(\ul{x}_{t}) \; \dot{x}^{j}_{t}(0)=\sum \limits_{j=1}^{2}q_{j}^{t-1,i}(\ul{x}_{t-1}) \; \dot{x}^{j}_{t-1}(0) \desda \\
\ul{Q}_{t}(\ul{x}_{t}) \ul{D}_t(\ul{x}_{t-1}) \, \dot{\ul{x}}_{t-1}= \ul{Q}_{t-1}(\ul{x}_{t-1}) \, \dot{\ul{x}}_{t-1} \\
\end{array}
\]
%chainrule
with $\ul{D}_{t}=[D^{j}_{k}]=\left[\frac{\partial x_{t}^{j}}{\partial x_{t-1}^{k}}\right] \in \R^{2\times 2}$,
$\ul{Q}_{t}=\left[q_{j}^{t,i}\right] \in \R^{N\times 2}$ and
$\dot{\ul{x}}_{t}(0)=(\dot{x}^{1}_{t},\dot{x}^{2}_t)=\frac{d}{ds}\ul{x}_{t}(0)$ arbitrary so that
\[
\ul{Q}_{t}(\ul{x}_{t}) \ul{D}_t(\ul{x}_{t-1})=\ul{Q}_{t-1}(\ul{x}_{t-1}).
\]
Provided that the frequency estimates are smooth and slowly varying with respect to the position variable one can apply first order Taylor expansion on the left-hand side and evaluate $\ul{Q}_{t}(\ul{x}_{t-1})$ in stead of $\ul{Q}_{t}(\ul{x}_{t})$. This corresponds to linear deformation theory, where $\ul{D}_t(\ul{x}_{t})\approx \ul{D}_t(\ul{x}_{t-1})$ provided that displacements \mbox{$\|\ul{x}_{t}-\ul{x}_{t-1}\|$} are small.
%so that the Lagrangian and Eulerian approaches coincide
Consequently, we use per (fixed) spatial position $\ul{x} \in \R^{2}$, the least squares solution of $\ul{Q}_{t}(\ul{x}) \ul{D}_t(\ul{x})=\ul{Q}_{t-1}(\ul{x})$ given by
\begin{equation} \label{Dtx}
\ul{D}_t(\ul{x})= ((\ul{Q}_{t}(\ul{x}))^{T}\ul{Q}_{t}(\ul{x}))^{-1}(\ul{Q}_{t}(\ul{x}))^{T}\ul{Q}_{t-1}(\ul{x})
\end{equation}
as our deformation gradient estimate. %The corresponding Green-Lagrangian strain tensor equals %$\ul{E}_{t}(\ul{x})=\frac{1}{2}(\ul{D}_{t}(\ul{x})^{T}\ul{D}_{t}(\ul{x})-I)$.
For experiments on a phantom with both considerable non-uniform scaling, fading, and non-uniform rotation see Figure \ref{fig:exp}.
The iteratively computed deformation net is close to the ground truth deformation net.

The deformation nets are computed as follows. We fix a single material point at the boundary (the green point in Figure \ref{fig:exp}) and first compute the grid points on the outer contour in circumferential direction. From these points we compute the inner grid-points iteratively in inward radial direction. Let $\ul{x}_{t}^{r,j}$ denote the $j$-th grid-point in the deformation net at time $t>0$ on the $r$-th closed contour.
Let $J$ denote the number of points on one closed
circumferential contour and let $R$ denote the number of points on a circumferential contour, let $T$ be the total time (number of tags), then
\begin{equation}\label{defnet}
\begin{array}{l}
\textrm{Single material point } \ul{x}_{t}^{1,1} \textrm{ given for all }t \in \{0,\ldots T-1\}, \\
\textrm{Initial polar grid }
\ul{x}_{0}^{r,j} \textrm{ given for all } r \in \{1,\ldots R\}, j \in \{1,\ldots J\}, \\
\textrm{Deformation field } D_{t} \textrm{ computed for all }t \geq 0 \textrm{ from the local frequency fields }\{\ul{q}^{t,i}\}_{i=1}^{N} \textrm{ via Eq.~(\ref{Dtx})}, \\
%{\rm d}\ul{x}_{0}^{1,i}:= \ul{x}_{t}^{1,i+1}-\ul{x}_{t}^{1,i}, \textrm{ for }i \in \{1,\ldots N\}, \\
\textrm{Then for }t = 1 \textrm{ to } \textrm{T-1} \textrm{ do }\\
\textrm{for }j = 1 \textrm{ to } \textrm{J} \textrm{ compute }
\ul{x}_{t}^{1,j+1}:= \ul{x}_{t}^{1,j}+  \left.\ul{D}_{t}\right|_{\ul{x}_{t}^{1,j}} (\ul{x}_{t-1}^{1,j+1}-\ul{x}_{t-1}^{1,j}), \\
\textrm{for }r = 1 \textrm{ to } \textrm{R} \textrm{ compute } \ul{x}_{t}^{r+1,j}:= \ul{x}_{t}^{r,j}+  \left.\ul{D}_{t}\right|_{\ul{x}_{t}^{r,j}} (\ul{x}_{t-1}^{r+1,j}-\ul{x}_{t-1}^{r,j}). %, \\
%
%
%{\rm d} \ul{x}_{t}^{1,i}:=\ul{x}_{t}^{1,i+1}-\ul{x}_{t}^{1,i} \\
%\ul{x}_{0}^{1,i}=\ul{x}_{0}^{1,i-1}+ \ul{D}_{t}(\ul{x}_{0}^{1,i-1})({\rm d}\ul{x}_{0}^{1,i-1}) \textrm{ for }i=1, \ldots N \in \mathbb{N}.\\
%\ul{x}_{t}^{R+1,i}=\ul{x}_{t}^{R,i}+ \ul{D}_{t}(\ul{x}_{0}^{R+1,i}) \; \\
%(\ul{x}_{0}^{R+1,i}-\ul{x}_{0}^{R,i})
\end{array}
\end{equation}
In contrast to well-established optical flow methods, \cite{HandSpaper,LucasKanade,BruhnetAl,Florackharp,Flor98a,vanAssen2007CBfMII}
this method is robust with respect to inevitable fading of the tag-lines, and can be applied directly on the original MRI-tagging images
without harmonic phase or sine phase pre-processing \cite{Osman1999,Florackharp}. Some of the optic flow methods, such as \cite{Janssen2006b,BecciuCAIP,Dorst2009}, do allow direct computation on the original MRI-tagging images as well, but they are relatively expensive, complicated and technical, e.g.~\cite{Duits2011QAM}.

First experiments show that our method can handle relatively large nonlinear deformations in both
radial and circumferential direction. The frequency and deformation fields look very promising from the qualitative viewpoint, e.g. Figure \ref{fig:exp} and Figure \ref{fig:tagswithGT}.
Quantitative comparisons to other methods such as \cite{Arts,VanAssen2b,FlorackDEF} and investigation
of the reliability of the frequency estimates (and their connection to left-invariant evolutions on Gabor transforms) are interesting topics for future work.

Although Gabor transforms are widely used in 1D-signal processing, they are less common in 2D-imaging. %compared to, for example, wavelet transforms on the similitude group.
However, the results in this section indicate applications of both Gabor transforms and the left-invariant evolutions acting on them in 2D-image processing.

%, but fall without the scope of this article.

\ \\ %\ \\
\textbf{Acknowledgements }The authors gratefully acknowledge Jos Westenberg (Leiden University Medical Center) for providing us the MRI-tagging data sets that were acquired in 4 different tagging directions.
%dr.~H.C.~van Assen Eindhoven University of Technology (TU/e) for drawing our attention to the cardiac imaging application in %Section \ref{ch:pract}, for pointing us to \cite{Arts} and
%for providing us the tagged-MRI datasets. Furthermore, the authors gratefully acknowledge Prof.~L.M.J.~Florack (TU/e)
%for discussions and suggestions in the relation between frequency fields and deformation fields via tensor calculus.
%in Section \ref{ch:pract}.

\appendix

\section{Existence and Uniqueness of the Evolution Solutions \label{ch:unique}}
The convection diffusion systems (\ref{theeqs}) have unique solutions, since the coefficients $a_{i}$ and $D_{ij}$ depend smoothly on the modulus of the initial condition $|\mathcal{W}_{\psi}f|=|\mathcal{G}_{\psi}f|$. So for a given initial condition $\mathcal{W}_{\psi}f$ the left-invariant convection diffusion generator $Q(|\mathcal{W}_{\psi}f|,\mathcal{A}_{1},\ldots,\mathcal{A}_{2d})$ is of the type
\[
Q(|\mathcal{W}_{\psi}f|,\mathcal{A}_{1},\ldots,\mathcal{A}_{2d})=\sum \limits_{i=1}^{d}\alpha_{i}\mathcal{A}_{i} + \sum \limits_{i,j=1}^{d}\mathcal{A}_{i} \beta_{ij} \mathcal{A}_{j}.
\]
Such hypo-elliptic operators with almost everywhere smooth coefficients given by $\alpha_{i}(p,q)=a_{i}(|\mathcal{G}_{\psi}f|)(p,q)$ and $\beta_{ij}(p,q)=D_{ij}(|\mathcal{G}_{\psi}f|)(p,q)$ generate strongly continuous, %\cite{Yosida},
semigroups on $\mathbb{L}_{2}(\R^2)$, as long as we keep the functions $\alpha_{i}$ and $\beta_{ij}$ fixed, \cite{Taylor,Hoermander},
yielding unique solutions where at least formally we may write
\[
W(p,q,s,t)=\Phi_{t}(\mathcal{W}_{\psi}f)(p,q,s)= e^{t \left(\sum \limits_{i=1}^{d}\alpha_{i}\mathcal{A}_{i} + \sum \limits_{i,j=1}^{d}\mathcal{A}_{i}\, \beta_{ij} \, \mathcal{A}_{j}\right)} \mathcal{W}_{\psi}f(p,q,s)
\]
with $\lim_{t \downarrow 0} W(\cdot,t)=\mathcal{W}_{\psi}f$ in $\mathbb{L}_{2}$-sense.
Note that if $\psi$ is a Gaussian kernel and $f \neq 0$ the Gabor transform $\mathcal{G}_{\psi}f \neq 0$ is real analytic on $\R^{2d}$, so it can not vanish on a set with positive measure, so that $\alpha_{i}: \R^2 \to \R$ are almost everywhere smooth. %The $D_{ij}$
%Here we note that by \cite{Evans} $\mathcal{G}_{\psi}f \in \mathbb{L}_{2}$
This applies in particular to the first reassignment approach in (\ref{conchoice}) (mapping everything consistently into phase space using Theorem \ref{th:correspondence}), where we
have set
$D=0$, $a_{1}(|\mathcal{G}_{\psi}f|)=\mathcal{M}(|\mathcal{G}_{\psi}f|)|\mathcal{G}_{\psi}f|^{-1} \partial_{p}|\mathcal{G}_{\psi}f|$
and $a_{2}(|\mathcal{G}_{\psi}f|)=\mathcal{M}(|\mathcal{G}_{\psi}f|)|\mathcal{G}_{\psi}f|^{-1} \partial_{q}|\mathcal{G}_{\psi}f|$.
In the second approach in (\ref{conchoice}) the operator $\tilde{U} \mapsto \tilde{\mathcal{C}}(\tilde{U})$ is non-linear, left-invariant and maps the space $\mathbb{L}_{2}^{+}(\R^2)=\{f \in \mathbb{L}_{2}(\R^2)\; |\; f \geq 0\}$ into itself again. In these cases the erosion solutions (\ref{erosion}) are the unique viscosity solutions, of (\ref{conpde}), see \cite{Crandall}. \\
\begin{remark}For the diffusion case, \cite[ch:7]{DuitsFuehrJanssen}, \cite[ch:6]{Janssen2009b}, we have \mbox{$D=[D_{ij}]_{i,j=1,\ldots, 2d}>0$}, in which case the (horizontal) diffusion generator $Q(|\mathcal{W}_{\psi}f|,\mathcal{A}_{1},\ldots,\mathcal{A}_{2d})$ on the group is hypo-elliptic, whereas the corresponding generator $\tilde{Q}(|\mathcal{G}_{\psi}f|,\tilde{\mathcal{A}}_{1},\ldots,\tilde{\mathcal{A}}_{2d})$ on phase space is elliptic. By the results \cite[ch:7.1.1]{Evans} and \cite{Lions} there exists a unique weak solution $\tilde{W}=\mathcal{S}W \in \mathbb{L}_{2}(\R^{+}, \mathbb{H}_{1}(\R^2)) \cap \mathbb{H}_{1}(\R^{+}, \mathbb{L}_{2}(\R^2))$. %and point evaluation is continuous in time
%
%and operator $\mathbb{L}_{2}(\R^2) \ni \mathcal{G}_{\psi}f \mapsto %\tilde{\Phi}_{t}(\mathcal{G}_{\psi}f):=\tilde{W}(\cdot, \cdot,t) \in %\mathbb{L}_{2}(\R^2)$ is well-defined for all $t>0$.
By means of the combination of Theorem \ref{th:one} and Theorem \ref{th:correspondence} we can transfer the existence and uniqueness result for the elliptic diffusions on phase space
to the existence and uniqueness result for the hypo-elliptic diffusions on the group $H_{r}$ (that is via conjugation with $\mathcal{S}$).
\end{remark}
{\small
\bibliographystyle{elsarticle-num}

%\bibliography{lit_duits_fuehr_janssen}
}
%% Authors are advised to submit their bibtex database files. They are
%% requested to list a bibtex style file in the manuscript if they do
%% not want to use elsarticle-num.bst.

%% References without bibTeX database:

% \begin{thebibliography}{00}

%% \bibitem must have the following form:
%%   \bibitem{key}...
%%

% \bibitem{}

% \end{thebibliography}

\end{document}